\documentclass[12pt]{article}

\usepackage{latexsym}
\usepackage{a4wide}
\usepackage{amscd}
\usepackage{graphicx}
\usepackage{amsmath}
\usepackage{amssymb}
\usepackage{soul,xcolor}
\usepackage{esint}
\usepackage{mathrsfs}
\usepackage{amsthm}
\usepackage{bbm}
\usepackage{color}
\usepackage{comment}
\usepackage{enumerate}
\usepackage{ulem}
\input xy
\xyoption{all}

\usepackage{hyperref}

\newcommand{\N}{\mathbb{N}}                     
\newcommand{\Z}{\mathbb{Z}}                     
\newcommand{\R}{\mathbb{R}}                     
\newcommand{\C}{\mathbb{C}}                     
\newcommand{\U}{\mathbb{U}}                     
\newcommand{\T}{\mathbb{T}}                     
\newcommand{\re}{\mathrm{Re\,}}                 
\newcommand{\Ker}{\mathrm{Ker\,}}               

\newcommand{\di}{\mathrm{d}}                     
\newcommand{\tr}{\mathrm{tr}\,}			

\newtheorem{mainthm}{\sc Theorem}[subsection]           
\newtheorem{mainprop}[mainthm]{\sc Proposition}           
\newtheorem{mainrem}[mainthm]{\sc Remark}           
\newtheorem{mainque}[mainthm]{\sc Question}           
\newtheorem{maindefn}[mainthm]{\sc Definition}        
\newtheorem{mainex}[mainthm]{\sc Example}           

\newtheorem{thm}{\sc Theorem}[section]              
\newtheorem*{thm*}{\sc Theorem}               
\newtheorem*{cor*}{\sc Corollary}        
\newtheorem{lem}[thm]{\sc Lemma}            
\newtheorem{prop}[thm]{\sc Proposition}     
\newtheorem{rem}[thm]{\sc Remark}           
\newtheorem{ex}[thm]{\sc Example}          

\def\gray#1{{\textcolor{darkgray}{#1}}}

\title{Lorentz--Finsler metrics on symplectic and contact transformation groups}

\author{Alberto Abbondandolo, Gabriele Benedetti and Leonid Polterovich}

\begin{document}

\maketitle

\vspace{-.5cm}

\begin{abstract}
It has been noticed a while ago that several fundamental transformation groups
of symplectic and contact geometry carry natural causal structures, i.e., fields of tangent convex cones. 
Our starting point is that quite often the latter come together with Lorentz--Finsler metrics, a notion originated in relativity theory, which enable one to do geometric measurements with timelike curves. This includes finite-dimensional linear symplectic groups, where these metrics can be seen as Finsler generalizations of the classical anti-de Sitter spacetime, infinite-dimensional groups of contact transformations, with the simplest example being the group of circle diffeomorphisms, and symplectomorphism groups of convex domains. In the first two cases, the Lorentz--Finsler metrics we introduce are bi-invariant. A Lorentz--Finsler perspective on these transformation groups turns out to be unexpectedly rich: some basic questions about distance, geodesics and their conjugate points, and existence of time functions, are naturally related to the contact systolic problem, group quasi-morphisms, the Monge--Amp\`ere equation, and a subtle interplay between symplectic rigidity and flexibility.  We discuss these interrelations,
providing necessary preliminaries, albeit mostly focusing on new results which have not been published before. 
Along the way, we formulate a number of open questions.
\end{abstract} 

\vspace{-.5cm}

\tableofcontents

\section*{Introduction and main results} \addcontentsline{toc}{section}{Introduction and main results}

Endow the vector space $\R^{2n}$ with linear coordinates $x_1,y_1,\dots,x_n,y_n$ and with the standard symplectic form
\[
\omega_0 := \sum_{j=1}^n \di x_j \wedge \di y_j.
\]
The group of linear automorphisms of $\R^{2n}$ that preserve $\omega_0$ is the symplectic group $\mathrm{Sp}(2n)$. It is well known that $\mathrm{Sp}(2n)$ admits no bi-invariant distance function inducing the Lie group topology. Here is the simple argument for $n=1$: the symplectic automorphisms
\[
W_{\lambda}:= \left( \begin{array}{cc} 1 & \lambda  \\ 0 & 1 \end{array} \right), \qquad \lambda>0,
\]
are all pairwise symplectically conjugate and hence any bi-invariant distance function on $\mathrm{Sp}(2)$ assigns the same positive distance from the identity to each of them. But then the distance function cannot be continuous with respect to the Lie group topology, as $W_{\lambda}$ converges to the identity for $\lambda\rightarrow 0$. The same argument applies in higher dimension and shows, in particular, that $\mathrm{Sp}(2n)$ does not admit bi-invariant Riemannian or Finsler metrics. 

Similarly, the contactomorphism group $\mathrm{Cont}(M,\xi)$ of a closed contact manifold $(M,\xi)$ does not admit any bi-invariant distance function which is continuous with reasonable Lie group topologies. More precisely, any bi-invariant distance function on $\mathrm{Cont}(M,\xi)$ is discrete, meaning that the distance of any pair of distinct elements has a positive lower bound, see \cite[Theorem 3.1]{fpr18}.
\bigskip

In this monograph, 
 we show that $\mathrm{Sp}(2n)$ and $\mathrm{Cont}(M, \xi)$ admit natural bi-invariant Lorentz--Finsler structures and initiate a systematic study of their properties.  The former provides yet another multi-dimensional generalization of the classical 3-dimensional anti-de-Sitter space and yields a new viewpoint at the twist condition in Hamiltonian dynamics. The latter (which is related to the former) provides a natural geometric language for studying a  non-autonomous version of the contact systolic problem. While our motivation comes from symplectic and contact geometry and dynamics, we develop the subject along the lines
which are customary in Lorentzian geometry, and which are influenced by its physical interpretation. As we believe that keeping in mind this interpretation may facilitate the
understanding of the (otherwise, purely mathematical) material of the present monograph, we start with its very brief overview.  

Lorentzian (or, more generally, Lorentz--Finsler) metrics are sign-indefinite cousins of Riemannian (resp.\ Finsler) metrics. They originated in the relativity theory as a natural geometric structure on a space-time invariant under Lorentz transformations, modeling the change of an inertial coordinate system. The necessity to deal with anisotropies of the space-time \cite{js14} motivated a passage from sign-indefinite Lorentzian quadratic forms to more general Lorentz--Finsler functionals having similar convexity/concavity features.

In the Lorentzian world,  the space-time $M$ which is modeled by a smooth manifold,
gets equipped with a field of cones consisting of vectors of the real (as opposed to the imaginary) length. Physically admissible causal (resp.\ timelike) curves are characterized by the  fact that their tangent vectors point into these cones (resp.\ into their interior). When there are no closed causal curves, the existence of a causal curve starting
at $x$ and ending at $y$ introduces a partial order on $M$, defining the so called causal structure.
The natural parameter along causal curves is neither length nor time, but a so called proper time, which is modeled by the Lorentzian length. Local extremizers of the Lorentzian length, i.e., timelike geodesics, are of particular interest: they model the motion of a free particle in the spacetime.

The behaviour of proper time is far from being intuitive. According to the famous ``twin paradox", moving close to the light cone consisting of vectors of the zero Lorentzian length enables a spacetime traveller who takes off at a point $x$ to arrive at the destination $y$ within an arbitrary small proper time. His twin sibling, simultaneously starting at $x$, may pursue a different objective - to reach $y$ within the maximal possible proper time. This quantity is sometimes  infinite and sometimes finite, depending on the geometry and topology of the spacetime, and maybe also on the specific choice of the points $x$ and $y$. When finite, it defines the Lorentzian distance $\mathrm{dist}(x,y)$, a global geometric invariant.

\bigskip

The main results of the present monograph together with some open questions are presented in Sections A - N of the
Introduction. After recalling the definition of a Lorentz--Finsler structure (Section \ref{secA}), we introduce a Lorentz--Finsler structure in parallel on the finite dimensional Lie group $\text{Sp}(2n)$ (Sections \ref{secB} and \ref{secC}) and on the infinite dimensional one $\text{Cont}(M, \xi)$ (Section \ref{secD}), as the two structures are closely related (Section \ref{secE}). We then study the local properties of the induced length functionals (i.e., proper time) and of its geodesics (Sections \ref{secF} and \ref{secG}), which in the infinite dimensional case is related to a contact systolic question (Section \ref{secH}). The Lorentzian viewpoint enables us
to establish ``systolic freedom" for time-dependent contact forms which manifests an interplay
between dynamics and geometry.   Loosely speaking, the proper time of our Lorentz--Finsler structures can be described dynamically, as a certain ``magnitude of twisting" of the flow corresponding to a path on the group, and geometrically, via the contact volume. 

Furthermore, we discuss to which extent these structures can be used in order to produce global bi-invariant measurements on these groups. For $\mathrm{Sp}(2n)$,
the Lorentz--Finsler distance $\mathrm{dist}(x,y)$ between causally related points $x$ and $y$  can take both finite and infinite values, depending on the location of $x^{-1}y$ (Sections \ref{secI} and \ref{secJ}). In contrast to this,
for contactomorphisms of the projective space, $\mathrm{dist}(x,y)$ is always infinite whenever there is a timelike curve from $x$ to $y$ (Section \ref{secK}). This can be seen as a manifestation of the flexibility of contactomorphisms. However, flexibility is expensive: long paths connecting $x$ and $y$  necessarily possess a ``high complexity", properly understood. 
For the simplest contact manifold $S^1= \R \mathrm{P}^1$, our approach to this phenomenon involves a delicate $L^1$-version of Bernstein's classical inequality for positive trigonometric polynomials due to Nazarov (Section \ref{secL}). In higher dimensions, we use an ingredient from ``hard" contact topology, namely Givental's non-linear Maslov index, in combination with the analysis on $\mathrm{Sp}(2n)$ which turns out to be crucial (Section \ref{secM}).

Finally, in Section \ref{secN} we discuss Lorentz--Finsler phenomena on the group of symplectomorphisms of a uniformly convex domain, which turn out to be related to the Monge--Amp\`ere equation and to a variational problem which is linked to the maximization of the affine area functional.

\renewcommand{\thesubsection}{\Alph{subsection}}
\numberwithin{equation}{subsection}

\subsection{Lorentz--Finsler structures} 
\label{secA}

Let $M$ be a (possibly infinite dimensional) manifold. In this monograph, we shall use the following notion of {\it Lorentz--Finsler structure} on $M$: 

\begin{maindefn}
\label{deflorfin}
{\rm A {\it Lorentz--Finsler structure} $(K,F)$ on $M$ is given by the following data:
\begin{enumerate}[(i)]
\item An open subset $K\subset TM$ such that for every $p\in M$ the intersection $K\cap T_p M$ is a non-empty convex cone in the vector space $T_p M$, and $\overline{K} \cap \overline{-K}$ coincides with the zero-section of $TM$. The set $K$ is called {\it cone distribution} on $M$.
\item A smooth function $F: K \rightarrow (0,+\infty)$ which is fiberwise positively 1-homogeneous, fiberwise strongly concave in all directions other than the radial one, meaning that
\[
\mathrm{d}^2 F(v)\cdot (w,w) < 0 \qquad \forall v\in K \cap T_p M, \; \forall w\in T_p M \setminus \R v,\; \forall p\in M,
\]
and extends continuously to $\overline{K}$ by setting $F|_{\partial K}=0$. The function $F$ is called {\it Lorentz--Finsler metric} on $(M,K)$.
\end{enumerate}}
\end{maindefn}

If $M$ is infinite dimensional, the smoothness of $F$ can be understood in several ways, depending on the class of infinite dimensional objects one is working with. In this monograph, we will work with a Fr\'echet manifold which is modeled on the space of smooth real functions on a closed manifold, and smoothness is to be understood in the diffeological sense: the restriction of $F$ to any finite dimensional submanifold of the open set $K$ is smooth. 

The above definition generalizes the classical notion of a {\it time-oriented Lorentz structure}, in which the manifold $M$ is endowed with a non-degenerate symmetric bilinear form $g:TM\times TM\rightarrow \R$ of signature $(-,+,\dots,+)$ and there is a continuous vector field $X$ on $M$ such that $g(X,X)<0$: indeed, in this case one chooses as $K$ the connected component of the set $\{v\in TM \mid g(v,v) <0\}$ containing the image of $X$ and sets $F(v):= \sqrt{-g(v,v)}$. The assumption on the signature of $g$ implies that $K$ is convex and $F$ is fiberwise strongly concave, as required in Definition \ref{deflorfin}.

Apart from regularity and strong convexity issues on the boundary of $K$, the above definition of a Lorentz--Finsler structure agrees with Asanov's definition from \cite{asa85} and its later refinements, see \cite{min16}, \cite{js20}. In particular, it agrees with the idea that a Lorentz--Finsler metric needs to be defined only on the convex cone of {\it causal vectors}.

Vectors in $K$ are called {\it timelike}, non-vanishing vectors in $\partial K$ are called {\it lightlike}, and vectors which are either timelike or lightlike are called {\it causal}. A $C^1$ curve in $M$ is called {\it timelike} (resp.\  {\it lightlike}, resp.\ {\it causal}) if its derivative is everywhere timelike (resp.\ lightlike, resp.\ causal). The {\it Lorentz--Finsler length} of a $C^1$ causal curve $\gamma:[a,b]\rightarrow M$ is the non-negative number
\[
\mathrm{length}_F(\gamma) := \int_a^b F(\gamma'(t))\, \di t.
\]
This functional is invariant under orientation preserving reparametrizations and additive under juxtaposition of curves. Moreover, it is positive and has directional derivatives of every order at each timelike curve. 

\subsection{A bi-invariant Lorentz--Finsler structure on the linear symplectic group} 
\label{secB}

The Lie algebra of the linear symplectic group $\mathrm{Sp}(2n)$ is
\[
\mathrm{sp}(2n) := \{ X\in \mathrm{Hom}(\R^{2n},\R^{2n}) \mid \mbox{the bilinear form } (u,v)\mapsto \omega_0(u,Xv) \mbox{ is symmetric}\},
\]
and we consider the subset
\[
\mathrm{sp}^+(2n) := \{ X\in \mathrm{sp}(2n) \mid (u,v)\mapsto \omega_0(u,Xv) \mbox{ is positive definite}\},
\]
which  is an open convex cone. All the elements of $\mathrm{sp}^+(2n)$ have positive determinant, and the function
\[
G: \mathrm{sp}^+(2n) \rightarrow \R, \qquad G(X) = (\det X)^{\frac{1}{2n}},
\]
is positive, positively 1-homogeneous, smooth, strongly concave in every direction other than the radial one and extends continuously (but not smoothly) to the closure of $\mathrm{sp}^+(2n)$ by setting it to be zero on the boundary. 

An endomorphism $X\in \mathrm{Hom}(\R^{2n},\R^{2n})$ belongs to $\mathrm{sp}^+(2n)$ if and only if 
\begin{equation}
\label{isdiag}
X = \bigoplus_{j=1}^n \theta_j J_j,
\end{equation}
where the direct sum refers to a symplectic splitting of $\R^{2n}$ into $n$ pairwise $\omega_0$-orthogonal symplectic planes $V_j\subset\R^{2n}$, each $\theta_j$ is a positive number and each $J_j$ is an $\omega_0$-compatible complex structure on $(V_j,\omega_0|_{V_j})$ (recall that a complex structure $J$ on a symplectic vector space $(V,\omega)$ is said to be $\omega$-compatible if the bilinear form $(u,v) \mapsto \omega(u,Jv)$ is symmetric and positive definite on $V$). See Proposition \ref{propA1} in Appendix \ref{linear-app} for a proof of this characterization of the elements of $\mathrm{sp}^+(2n)$. If $X$ has the form (\ref{isdiag}), then $G(X)$ is the geometric mean of the positive numbers $\theta_1,\dots, \theta_n$:
\begin{equation}
\label{geomean}
G(X) = \sqrt[n]{\theta_1\cdots \theta_n}.
\end{equation}
The cone $\mathrm{sp}^+(2n)$ and the function $G$ are easily seen to be invariant under the adjoint action of $\mathrm{Sp}(2n)$ on $\mathrm{sp}(2n)$. Therefore, $\mathrm{sp}^+(2n)$ extends by translation to a bi-invariant cone distribution 
\begin{equation}
\label{cone}
\{\mathrm{sp}^+(2n)W \mid W\in \mathrm{Sp}(2n)\} = \{ W \mathrm{sp}^+(2n) \mid W\in \mathrm{Sp}(2n)\} \subset T\mathrm{Sp}(2n)
\end{equation}
in the tangent bundle of $\mathrm{Sp}(2n)$, and $G$ extends to a bi-invariant function on this set. 

\begin{mainprop}
\label{lorfinSp}
The pair $(\mathrm{sp}^+(2n),G)$ defines a bi-invariant Lorentz--Finsler structure on $\mathrm{Sp}(2n)$.
\end{mainprop}

The easy proof is contained in Section \ref{lorfinSpsec} below. We shall denote this bi-invariant Lorentz--Finsler structure on $\mathrm{Sp}(2n)$ simply by $(\mathrm{sp}^+(2n),G)$, without introducing a special name for the cone distribution (\ref{cone}).

Studying causality on $\mathrm{Sp}(2n)$ with the above bi-invariant cone distribution means understanding the behaviour of timelike and causal curves on $(\mathrm{Sp}(2n),\mathrm{sp}^+(2n))$. Although not under this terminology, the study of causality on $\mathrm{Sp}(2n)$ is a classical subject. Indeed, since the elements of $\mathrm{sp}^+(2n)$ have the form $J_0 S$, where $J_0$ is the standard $\omega_0$-compatible complex structure on $\R^{2n}$ satisfying
\[
\omega_0(u,v) = J_0 u \cdot v \qquad \forall u,v\in \R^{2n},
\]
and $S$ belongs to the cone $\mathrm{Sym}^+(2n)$ of positive definite symmetric endomorphisms of $\R^{2n}$, a continuously differentiable curve $W: [a,b] \rightarrow \mathrm{Sp}(2n)$ is timelike if and only if it solves the non-autonomous positive definite linear Hamiltonian system
\[
W'(t) = J_0 S(t) W(t),
\]
for some continuous path $S: [a,b] \rightarrow \mathrm{Sym}^+(2n)$. For this reason, timelike curves in $\mathrm{Sp}(2n)$ are also called {\it positive paths} of linear symplectomorphisms. Similarly, causal curves are solutions of a non-autonomous linear Hamiltonian system as above with $S(t)$ non-zero and positive semi-definite for every $t$.

Positive definite linear Hamiltonian system have been widely studied due to their special role in Krein's stability theory of linear Hamiltonian systems, see  e.g., \cite{kre50,kre51,kre55,gl58,kl62}. Comprehensive expositions of Krein's stability theory and of the theory of positive definite linear Hamiltonian systems can be found in \cite[Chapter III]{ys75} and \cite[Chapter I]{eke90}. More results about positive paths in $\mathrm{Sp}(2n)$ can be found in \cite{lm97}. More generally, the study of invariant convex cones in Lie algebras, such as $\mathrm{sp}^+(2n)$, is a classical topic in Lie theory, see e.g., \cite{vin80,pan81,ols81b,ols81,ols82}. 

The novelty here is the study of the Lorentz--Finsler metric $G$ which, albeit very natural, does not seem to have been received much attention, except for the special case $n=1$, which corresponds to a classical spacetime in general relativity.

\subsection{The anti-de Sitter case ${\mathbf n=1}$}
\label{secC}

In the special case $n=1$, $\mathrm{Sp}(2n) = \mathrm{Sp}(2)$ coincides with $\mathrm{SL}(2,\R)$ and the Lorentz--Finsler metric $G$ comes from a genuine Lorentz metric, corresponding to the three-dimensional anti-de Sitter spacetime $\mathrm{AdS}_3$. We recall that this time-orientable Lorentz manifold can be defined as the restriction $g$ of a symmetric bilinear form $b$ of signature $(2,2)$ on a 4-dimensional real vector space $V$ to the hypersurface
\[
\mathrm{AdS}_3:= \{v\in V \mid b(v,v)=-1\}.
\]
By choosing $V=\mathrm{Hom}(\R^2,\R^2)$ and $b$ to be the symmetric bilinear form on $V$ whose associated quadratic form is $-\det$, we see that 
\[
\mathrm{AdS}_3 = \mathrm{SL}(2,\R) = \mathrm{Sp}(2),
\]
and for every $W$ in this manifold, $\mathrm{sp}^+(2)W$ is precisely one component of the cone of timelike vectors at $W$, hence we can choose it to be the cone of future pointing timelike vectors. Finally,
\[
G(Y) = |g(Y,Y)|^{\frac{1}{2}} \qquad \forall Y\in T\mathrm{AdS}_3,
\]
is precisely the Lorentz norm induced by the Lorentz metric $g$.

When $n>1$, the determinant is not a quadratic form anymore and the Lorentz--Finsler metric $G$ is not induced by a Lorentz metric. 

\begin{mainrem}
\label{unique-1}
{\rm The cones $\mathrm{sp}^+(2n)$ and $-\mathrm{sp}^+(2n)$ are the unique invariant open convex cones that are proper subsets of $\mathrm{sp}(2n)$, see \cite{pan81}. In the case $n=1$, the Lorentz metric $G$ is, up to the multiplication by a positive number, the unique bi-invariant Lorentz--Finsler metric on $(\mathrm{Sp}(2),\mathrm{sp}^+(2))$. For $n\geq 2$, there are other bi-invariant Lorentz--Finsler metrics on $(\mathrm{Sp}(2n),\mathrm{sp}^+(2n))$. For instance, one can check that the function $H: \mathrm{sp}^+(2n) \rightarrow (0,+\infty)$ given by the quadratic harmonic mean of the numbers $\theta_j$ appearing in (\ref{isdiag}), i.e.,
\[
H(X) = \sqrt{2n} \bigl( - \mathrm{tr}\, (X^{-2}) \bigr)^{-\frac{1}{2}} = \left( \sum_{j=1}^n \frac{1}{\theta_j^2} \right)^{-\frac{1}{2}},
\]
defines a bi-invariant Lorentz--Finsler metric on $(\mathrm{Sp}(2n),\mathrm{sp}^+(2n))$ which for $n>1$ is not a multiple of $G$. More generally, any positively 1-homogeneous function $f:(0,+\infty)^n \rightarrow (0,+\infty)$ which is invariant under permutations of the coordinates induces via the formula
\[
F(X) = f(\theta_1,\dots,\theta_n) \qquad \mbox{for } X = \bigoplus_{j=1}^n \theta_j J_j
\]
a 1-homogeneous function on $\mathrm{sp}^+(2n)$ which is invariant under the adjoint action of $\mathrm{Sp}(2n)$. If the function $f$ is smooth, then so is $F$, thanks to Glaeser's differentiable version of Newton's theorem on the representation of symmetric functions, see \cite{gla63}. Moreover, if $f$ extends continuously to the closure of its domain by setting it to be zero on the boundary, the same is true for $F$. 
It is unclear to us whether the strong concavity of $f$ in all directions other than the radial one imply the corresponding property for $F$. Therefore, we raise the following question.} 
\end{mainrem}

\begin{mainque}
Which functions  $f:(0,+\infty)^n \rightarrow (0,+\infty)$ as above define a bi-invariant Lorentz--Finsler metric on $(\mathrm{Sp}(2n),\mathrm{sp}^+(2n))$?
\end{mainque}

See \cite{dav57} and \cite{lew00} for related results on the convexity of functions which are defined in terms of the eigenvalues.

\subsection{A bi-invariant Lorentz--Finsler metric on the contactomorphism group} 
\label{secD}

Let $\xi$ be a {\it co-oriented contact structure} on the closed $(2n-1)$-dimensional manifold $M$, where $n\geq 1$. We recall that this means that $\xi$ is the kernel of a {\it contact form} on $M$, i.e., a smooth 1-form $\alpha$ such that $\alpha\wedge \mathrm{d}\alpha^{n-1}$ is a volume form on $M$, and $\alpha$ is positive on each tangent vector that is positively transverse to $\xi$. A 1-form $\alpha$ as above is called a {\it defining contact form} for $\xi$. The {\it volume} of $M$ with respect to the contact form $\alpha$ is denoted by
\[
\mathrm{vol}(M,\alpha) := \int_M \alpha\wedge \mathrm{d} \alpha^{n-1},
\]
and the \textit{Reeb vector field of $\alpha$ is the vector field $R_{\alpha}$  which is defined by the identities
\[
\imath_{R_{\alpha}} \mathrm{d} \alpha = 0, \qquad \imath_{R_{\alpha}} \alpha = 1.
\] }
The group of smooth diffeomorphisms of $M$ that preserve the co-oriented contact structure $\xi$ is called {\it contactomorphism group} of $(M,\xi)$ and denoted by $\mathrm{Cont}(M,\xi)$. The connected component containing the identity is denoted by $\mathrm{Cont}_0(M,\xi)$. 

The Lie algebra of the infinite dimensional Lie group $\mathrm{Cont}(M,\xi)$ is the space $\mathrm{cont}(M,\xi)$ of {\it contact vector fields} on $(M,\xi)$, i.e., smooth vector fields on $M$  whose flow preserves $\xi$. Having fixed a defining contact form $\alpha$ for $\xi$, the space $\mathrm{cont}(M,\xi)$ can be identified with the space of real functions $C^{\infty}(M)$ by the map
\[
\mathrm{cont}(M,\xi) \rightarrow C^{\infty}(M), \qquad X\mapsto \imath_X \alpha,
\]
where the function $\imath_X \alpha$ is called {\it contact Hamiltonian} of the contact vector field $X$ with respect to the contact form $\alpha$. See Appendix \ref{contham} for some basic facts about this identification.

Denote by $\mathrm{cont}^+(M,\xi)$ the subset of $\mathrm{cont}(M,\xi)$ consisting of those contact vector fields $X$ that are positively transverse to $\xi$. If $\alpha$ is a defining contact form for $\xi$, $\mathrm{cont}^+(M,\xi)$ is the space of contact vector fields $X$ such that $\imath_X \alpha>0$, and in the above identification with $C^{\infty}(M)$ it corresponds to the set of positive Hamiltonians. It is easy to check that $\mathrm{cont}^+(M,\xi)$ is precisely the set of Reeb vector fields associated to all contact forms defining the co-oriented contact structure $\xi$ (see Appendix \ref{contham}).

Note that $\mathrm{cont}^+(M,\xi)$ is an open convex cone in $\mathrm{cont}(M,\xi)$ and
\[
\overline{\mathrm{cont}^+(M,\xi)} \cap \overline{-\mathrm{cont}^+(M,\xi)} = \{0\}.
\]
Here, $\mathrm{cont}(M,\xi)$ is equipped with an arbitrary metrizable vector space topology which, after the identification with $C^{\infty}(M)$, is not coarser than the $C^0$-topology of functions. For instance, we may use the $C^k$-topology on $\mathrm{cont}(M,\xi)$ for any $0\leq k \leq \infty$.

We define a real function $V: \mathrm{cont}^+(M,\xi) \rightarrow \R$ by
\begin{equation}
\label{LMV}
V(X) :=  \mathrm{vol}(M, \alpha)^{-\frac{1}{n}},
\end{equation}
where $\alpha$ is the unique contact form defining $\xi$ such that $X=R_{\alpha}$.

The adjoint action of $\mathrm{Cont}(M,\xi)$ on $\mathrm{cont}(M,\xi)$ is given by the push-forward:
\[
\mathrm{Ad}_{\phi} X = \phi_* X, \qquad \forall \phi\in \mathrm{Cont}(M,\xi), \; X\in \mathrm{cont}(M,\xi).
\]
The cone $\mathrm{cont}^+(M,\xi)$ is invariant under the adjoint action. Therefore, $\mathrm{cont}^+(M,\xi)$ extends to a bi-invariant cone distribution in the tangent bundle of $\mathrm{Cont}(M,\xi)$. Moreover, $V$ is invariant under the adjoint action and hence extends to a bi-invariant function on the bi-invariant cone distribution generated by $ \mathrm{cont}^+(M,\xi)$.

\begin{mainprop}
\label{lorfinCont}
The pair $(\mathrm{cont}^+(M,\xi),V)$ defines a bi-invariant Lorentz--Finsler structure on $\mathrm{Cont}(M,\xi)$.
\end{mainprop}

See Section \ref{lorfinContsec} for the proof.  As in the case of the linear symplectic group, we denote this bi-invariant Lorentz--Finsler structure simply by $(\mathrm{cont}^+(M,\xi),V)$.

It is instructive to look at the one-dimensional manifold $\T := \R/\Z$, which is a contact manifold with the trivial contact structure $\xi_0:=\{0\}$ that is co-oriented by the standard orientation of $\T$. In this case, $\mathrm{Cont}(\T,\xi_0)=\mathrm{Cont}_0(\T,\xi_0)$ coincides with $\mathrm{Diff}_0(\T)$, the group of orientation preserving diffeomorphisms of $\T$, $\mathrm{cont}(\T,\xi_0)$ is the space of all tangent vector fields on $\T$, and $\mathrm{cont}^+(\T,\xi_0)$ is the cone of vector fields of the form
\[
X(x) = H(x) \frac{\partial}{\partial x},
\]
where $H$ is a positive function on $\T$. Any lift $\phi: \R \rightarrow \R$ of a diffeomorphism in $\mathrm{Diff}_0(\T)$ has a well-defined {\it translation number}
\[
\rho(\phi) := \lim_{n\rightarrow \infty} \frac{\phi^n(x)-x}{n} \qquad \forall x\in \R,
\]
and the quantity $V(X)$ has the following interpretation, which we prove in Section \ref{rotnumbsec}.

\begin{mainprop}
\label{rotnumb}
Let $X$ be an element of $\mathrm{cont}^+(\T,\xi_0)$ and let $\phi^t:\R \rightarrow \R$ be the lift of its flow such that $\phi^0=\mathrm{id}$. Then
\[
V(X) =\rho(\phi^1).
\]
\end{mainprop}

\begin{mainrem}
\label{unique-2}
{\rm In the simple case of the one-dimensional contact manifold $(\T,\xi_0)$, $V$ is the only bi-invariant Lorentz--Finsler metric on $(\mathrm{Cont}(\T,\xi_0),\mathrm{cont}^+(\T,\xi_0))$ up to the multiplication by a positive number. Actually, more is true: Any positive function $W:  \mathrm{cont}^+(\T,\xi_0) \rightarrow \R$ which is positively 1-homogeneous and invariant under the adjoint action of $\mathrm{Cont}_0(\T,\xi_0)$ has the form $W=cV$ for some positive number $c$. This uniqueness statement does not need concavity or continuity assumptions on $W$ and is a simple consequence of the fact that $\mathrm{Cont}_0(\T,\xi_0)$ acts transitively on rays in $\mathrm{cont}^+(\T,\xi_0)$ (see Proposition \ref{uniqueness2} below). The situation is therefore similar to the case of $(\mathrm{Sp}(2),\mathrm{sp}^+(2))$, see Remark \ref{unique-1} above.  If $\dim M>1$, $\mathrm{Cont}(M,\xi)$ does not act transitively on rays in $\mathrm{cont}^+(M,\xi)$ and there is an infinite dimensional family of  positive functions $W: \mathrm{cont}^+(M,\xi)\rightarrow \R$ which are positively 1-homogeneous and invariant under the adjoint action of $\mathrm{Cont}(M,\xi)$. It would be interesting to understand how much this family gets reduced by imposing concavity and continuity conditions on $W$.}
\end{mainrem}

\begin{mainque}
In the case $\dim M>1$, are there other bi-invariant Lorentz--Finsler metrics on $(\mathrm{Cont}(M,\xi), \mathrm{cont}^+(M,\xi))$? Can one classify them?
\end{mainque}

Any smooth path $\phi^t$ in $\mathrm{Cont}(M,\xi)$ induces a smooth path of contact vector fields $X_t$, which is uniquely defined by the equation
\begin{equation}
\label{gencont}
\frac{\mathrm{d}}{\mathrm{d}t} \phi^t = X_t(\phi^t).
\end{equation}
The path of contactomorphisms $\phi^t$ is said to be {\it positive} if $X_t$ belongs to $\mathrm{cont}^+(M,\xi)$ for every $t$, or equivalently if $\imath_{X_t} \alpha>0$, where $\alpha$ is a defining contact form for $\xi$.
Positive paths of contactomorphisms are hence the timelike curves of $(\mathrm{Cont}(M,\xi), \mathrm{cont}^+(M,\xi))$. Causal curves are instead {\it non-negative} paths of contactomorphisms that are non-constant, where non-negative means that $\imath_{X_t} \alpha\geq 0$ for every $t$ (with this terminology, a constant path is non-negative but in accordance with the use in general relativity we do not consider it to be a causal path). 

The study of positive paths of contactomorphisms was initiated in \cite{ep00,bhu01} and has developed into an important topic in contact geometry, see e.g., \cite{ekp06,cn10a,cn10b,afm15,cn16,am18,cn20}. Its relationship with Lorentzian geometry, which is not limited to the fact that $\mathrm{cont}^+(M,\xi)$ defines a causal structure on $\mathrm{Cont}(M,\xi)$,  is explicitly noticed and discussed in the above mentioned papers of Chernov and Nemirovski.  What is new here is the bi-invariant Lorentz--Finsler metric $V$ on such a cone distribution.

\begin{mainrem}
{\rm As recalled above, any bi-invariant distance on $\mathrm{Cont}(M,\xi)$ is discrete, see \cite[Theorem 3.1]{fpr18}. This does not exclude the existence of bi-invariant Finsler metrics on $\mathrm{Cont}(M,\xi)$, that is, positively 1-homogeneous positive functions on $\mathrm{cont}(M,\xi)$ which are strongly convex in any direction other than the radial one and invariant under the adjoint action of $\mathrm{Cont}(M,\xi)$. Indeed, the bi-invariant pseudo-distance which is induced by a bi-invariant Finsler metric on $\mathrm{Cont}(M,\xi)$ could be identically zero. For instance, on the group of Hamiltonian diffeomorphisms of a closed symplectic manifold the $L^p$-norm on the space of normalized Hamiltonians induces a bi-invariant Finsler metric, whose induced pseudo-distance vanishes identically if $p<+\infty$. Bi-invariant Finsler metrics on the group of Hamiltonian diffeomorphisms have been studied in \cite{ow05,bo11,lem20}. This raises the following:}
\end{mainrem}
\begin{mainque}
\label{quenofinsler}
Can there be bi-invariant Finsler metrics on $\mathrm{Cont}(M,\xi)$ for some contact manifold $(M,\xi)$?
\end{mainque}
The answer to this question is negative if we require the Finsler metric to be $C^0$-continuous, after identifying the Lie algebra $\mathrm{cont}(M,\xi)$ with the space of contact Hamiltonians $C^{\infty}(M)$. Actually, any $C^0$-continuous function on $C^{\infty}(M)$ which vanishes at zero and is invariant under the adjoint action of $\mathrm{Cont}(M,\xi)$ must vanish on contact Hamiltonians which are supported in Darboux charts, see Remark \ref{nofinsler} below, preventing this function to be a Finsler metric. Since any function can be written as a convex combination of functions with support in Darboux charts, this shows that concavity, rather than convexity, is the right condition to require when looking for interesting invariant non-negative functions on the closure of $\mathrm{cont}^+(M,\xi)$. We suspect that the answer to Question \ref{quenofinsler} remains negative also for $C^{\infty}$-continuous Finsler metrics and we can confirm this for the standard contact structure of spheres and real projective spaces, see Remark \ref{nofinslerintro} further down in this Introduction. The techniques developed in \cite{ow05,bo11,lem20} might be helpful in settling the above question. 

\subsection{The contactomorphism groups of $(S^{2n-1},\xi_{\mathrm{st}})$ and $(\R \mathrm{P}^{2n-1},\xi_{\mathrm{st}})$} 
\label{secE}

The Liouville 1-form
\[
\lambda_0 := \frac{1}{2} \sum_{j=1}^n \bigl( x_j \, \mathrm{d} y_j - y_j \, \mathrm{d}x_j \bigr)
\]
of $\R^{2n}$ restricts to a contact form on the unit sphere $S^{2n-1}$, and the corresponding contact structure $\xi_{\mathrm{st}}$ is the standard contact structure of $S^{2n-1}$. Being invariant under the antipodal map $z\mapsto -z$, $\lambda_0|_{S^{2n-1}}$ descends to a contact form on the real projective space $\R \mathrm{P}^{2n-1}$, and the corresponding contact structure is also denoted by $\xi_{\mathrm{st}}$.

Any linear automorphism of $\R^{2n}$ acts on rays from the origin and on lines through the origin and hence induces a diffeomorphism of $S^{2n-1}$ and a diffeomorphism of $\R \mathrm{P}^{2n-1}$. In the case of a symplectic automorphism, the resulting diffeomorphism are contactomorphisms and we obtain the injective homomorphisms
\[
i: \mathrm{Sp}(2n) \rightarrow \mathrm{Cont}_0 (S^{2n-1},\xi_{\mathrm{st}}), \qquad j: \mathrm{PSp}(2n) \rightarrow \mathrm{Cont}_0 (\R \mathrm{P}^{2n-1},\xi_{\mathrm{st}}),
\]
where $\mathrm{PSp}(2n)$ denotes the quotient of $\mathrm{Sp}(2n)$ by the normal subgroup $\{\mathrm{id},-\mathrm{id}\}$. Note that the Lorentz--Finsler structure $(\mathrm{sp}^+(2n),G)$ descends to a Lorentz--Finsler structure on   $\mathrm{PSp}(2n)$, which we shall denote by the same notation.

In the case $n=1$, $(S^1,\xi_{\mathrm{st}})$ and $(\R \mathrm{P}^1,\xi_{\mathrm{st}})$ are clearly contactomorphic to each other and to $(\T,\xi_0)$. In this case, we actually have for every natural number $k$ an injective homomorphism
\[
j_k : \mathrm{PSp}_k(2) \rightarrow \mathrm{Cont}(\T,\xi_0) = \mathrm{Diff}_0(\T)
\]
whose image is a subgroup of diffeomorphisms of $\T$ commuting with the translation by $\frac{1}{k}$. Here, 
\[
p_k: \mathrm{PSp}_k(2) \rightarrow \mathrm{PSp}(2)
\]
is the connected $k$-th fold cover and $j_k(w)$ is defined by lifting the diffeomorphism
\[
j(p_k(w)) : \R \mathrm{P}^1 \cong \T  \rightarrow \R \mathrm{P}^1 \cong \T
\]
to the $k$-th fold cover $q_k: \T  \rightarrow \T$. The diffeomorphism $j(p_k(w))$ has $k$ distinct possible lifts, and the element $w\in  \mathrm{PSp}_k(2,\R)$ dictates which one we are choosing. See \cite[pp.\ 341-342]{ghy01}. The Lorentz structure of $\mathrm{Sp}(2)\cong \mathrm{AdS}_3$ induces bi-invariant Lorentz--Finsler structures on all the groups $\mathrm{PSp}_k(2)$, and we denote them still by $(\mathrm{sp}^+(2),G)$.

The following result shows that the Lorentz--Finsler structures $(\mathrm{cont}^+(S^{2n-1},\xi_{\mathrm{st}}),V)$ and $(\mathrm{cont}^+(\R \mathrm{P}^{2n-1},\xi_{\mathrm{st}}),V)$ on the contactomorphism groups of the sphere and the real projective space are, up to rescaling factors, infinite dimensional extensions of the Lorentz--Finsler structure $(\mathrm{sp}^+(2n),G)$ on $\mathrm{Sp}(2n)$ and $\mathrm{PSp}(2n)$.
  
\begin{mainprop}
\label{sottogruppi}
The homomorphisms $i$, $j$ and $j_k$ satisfy
\[
\begin{split}
\mathrm{d} i(\mathrm{id})^{-1} \bigl( \mathrm{cont}^+(S^{2n-1},\xi_{\mathrm{st}}) \bigr) = \mathrm{sp}^+(2n),  \qquad & \mathrm{d} j(\mathrm{id})^{-1} \bigl(  \mathrm{cont}^+(\R \mathrm{P}^{2n-1},\xi_{\mathrm{st}}) \bigr) = \mathrm{sp}^+(2n),\\
\mathrm{d}j_k(\mathrm{id})^{-1}  \bigl( \mathrm{cont}^+(\T,\xi_0) \bigr) &=  \mathrm{sp}^+(2),
\end{split}
\]
and
\[
V \circ  \mathrm{d} i(\mathrm{id}) = \frac{1}{2\pi} G, \qquad V \circ \mathrm{d} j(\mathrm{id})= \frac{2^{\frac{1}{n}}}{2\pi}  G, \qquad V \circ  \mathrm{d} j_k(\mathrm{id}) = \frac{1}{k\pi} G.
\]
\end{mainprop}

This proposition, whose proof is discussed in Section \ref{sottogruppisec},  allows us to deduce results about the Lorentz--Finsler structures on the contactomorphism groups of spheres and real projective spaces from finite-dimensional results about the Lorentz--Finsler structure on the linear symplectic group. 

\begin{mainrem}
\label{nofinslerintro}
{\rm The existence of the homomorphisms $i$ and $j$ allows us to give a negative answer to Question \ref{quenofinsler} above for the standard contact structures of spheres and real projective spaces:
on $\mathrm{Cont}(S^{2n-1},\xi_{\mathrm{st}})$ and $\mathrm{Cont}(\R \mathrm{P}^{2n-1},\xi_{\mathrm{st}})$ there are no bi-invariant $C^{\infty}$-continuous Finsler metrics. Indeed, the pull-back by $i$ of such a metric on $\mathrm{Cont}(S^{2n-1},\xi_{\mathrm{st}})$ would be a bi-invariant continuous Finsler metric on $\mathrm{Sp}(2n)$. By the finite dimensionality of $\mathrm{Sp}(2n)$, this metric would induce a bi-invariant distance function which is continuous with respect to the Lie group topology, and we have already noticed that such a distance cannot exist on $\mathrm{Sp}(2n)$. The same argument with the homomorphism $j$ works for $\mathrm{Cont}(\R \mathrm{P}^{2n-1},\xi_{\mathrm{st}})$.}
\end{mainrem}

\subsection{Timelike geodesics on $\mathrm{Sp}(2n)$} 
\label{secF}

{\it Timelike geodesics} on a manifold $M$ endowed with a Lorentz--Finsler structure $(K,F)$ can be defined as smooth timelike curves $\gamma$ that are extremal points of the functional $\mathrm{length}_F$, meaning that the first variation of $\mathrm{length}_F$ along any variation of $\gamma$ fixing the end-points vanishes. Moreover, we require timelike geodesics to be parametrized in such a way that $F(\gamma')$ is constant. 

In the case of the Lorentz--Finsler structure $(\mathrm{sp}^+(2n),G)$ on $\mathrm{Sp}(2n)$, timelike geodesics are precisely the solutions of {\it autonomous} positive definite linear Hamiltonian systems, i.e., the curves of the form
\[
W(t) = e^{tX} W_0,
\]
where $X\in \mathrm{sp}^+(2n)$ and $W_0\in \mathrm{Sp}(2n)$. In Appendix \ref{liegroups}, we discuss this fact in general for Lie groups that are endowed with a bi-invariant Lorentz--Finsler structure. Note that timelike geodesics depend on the bi-invariant cone distribution but not on the bi-invariant Lorentz--Finsler metric on it.

From the representation (\ref{isdiag}) for the elements of $\mathrm{sp}^+(2n)$, we deduce that timelike geodesics are, up to a right or left translation, direct sums of $n$ one-parameter groups of $\omega_0$-positive planar rotations:
\[
W(t) = \left( \bigoplus_{j=1}^n e^{\theta_j t J_j} \right) W_0.
\]
A timelike geodesic $W$ as above is periodic if and only if the numbers $\theta_j^{-1}$ are all integer multiples of the same real number, and in general is quasi-periodic. For $n=1$, we recover the well known fact that all timelike geodesics in the anti de-Sitter space are periodic and have the same length $2\pi$. For $n>1$, we obtain also quasiperiodic tori of non-closed timelike geodesics.

We shall compute the second variation of the Lorentz--Finsler length functional $\mathrm{length}_G$ at a timelike geodesic segment $W:[0,T] \rightarrow \mathrm{Sp}(2n)$. Due to the invariance under reparametrizations, this second variation has an infinite dimensional kernel, but modding out the reparametrizations we obtain a symmetric bilinear form which has a finite dimensional kernel and a finite {\it Morse co-index}, i.e., dimension of a maximal subspace on which the second variation is positive definite. 

As usual, $T$ is said to be a {\it conjugate instant} if this kernel is non-trivial, and in this case the dimension of this kernel is the multiplicity of the conjugate instant $T$. General facts about bi-invariant structures on Lie groups imply that the elements of this kernel are given by \textit{Jacobi fields}, i.e.~the paths $Y:[0,1] \rightarrow \mathrm{sp}(2n)$ such that
\[
Y''=[X,Y'], 
\]
where $X\in \mathrm{sp}^+(2n)$ is the generator of the timelike geodesic, which vanish for $t=0$ and $t=T$. Note that the equation for Jacobi fields, as the equation for geodesics of which this is the linearization, does not depend on the Lorentz--Finsler metric $G$. This is a consequence of the fact that we are working with a bi-invariant structure on a Lie group. See Appendix \ref{liegroups} for more about this.

Moreover, the \textit{Lorentzian Morse index theorem} holds: the Morse co-index of every timelike geodesic segment $W:[0,T] \rightarrow \mathrm{Sp}(2n)$ equals the sum of the multiplicities of the conjugate instants in the open interval $(0,T)$. In particular, timelike geodesics are locally length maximizing. These facts, which are well known for timelike geodesics on a Lorentzian manifold (see e.g., \cite{bee96}), still hold in the Lorentz--Finsler setting thanks to the strong concavity of the  Lorentz--Finsler metric. 

Instead of proving this in general, we content ourselves of checking these facts for bi-invariant Lorentz--Finsler metrics on Lie groups in Appendix \ref{liegroups} and to specialize them to $\mathrm{Sp}(2n)$ in Section \ref{varsec}.

In the special case of the periodic timelike geodesic $t\mapsto e^{t J}$ on $(\mathrm{Sp}(2n),\mathrm{sp}^+(2n))$, where $J$ is any $\omega_0$-compatible complex structure on $\R^{2n}$, we obtain the following result.

\begin{mainthm} 
\label{Morse}
Let $J$ be an $\omega_0$-compatible complex structure on $\R^{2n}$.
The timelike geodesic $W:\R \rightarrow \mathrm{Sp}(2n)$, $W(t) = e^{tJ}$, has a conjugate instant at $T>0$ if and only if $T\in \pi \N$ and each such conjugate instant $T$ has multiplicity $n^2+n$. The timelike geodesic segment $W|_{[0,T]}$ has finite Morse co-index, which equals the sum of conjugate instants in the interval $(0,T)$, counted with multiplicity, i.e.,
\[
\mbox{\rm co-ind} \, W|_{[0,T]} = \left( \Bigl\lceil \frac{T}{\pi} \Bigr\rceil - 1 \right) (n^2+n).
\]
\end{mainthm}

See Section \ref{jacsec} and Theorem \ref{Morsegen} below for a complete analysis of the conjugate instants and the Morse co-index of an arbitrary timelike geodesic on $(\mathrm{Sp}(2n),\mathrm{sp}^+(2n))$.

In the special case of $\mathrm{Sp}(2)=\mathrm{AdS}_3$, all the timelike geodesics are periodic and, after translation and affine reparametrization, are of the form considered in the theorem above, which hence recovers the familiar fact that every timelike geodesic on $\mathrm{AdS}_3$ has conjugate instants of multiplicity two at each semi-integer multiple of its period. In particular, simple closed timelike geodesics have co-index two and hence are not local maximizers of the Lorentzian length.

\subsection{Timelike geodesics on $\mathrm{Cont}(M,\xi)$}
\label{secG}
  
Timelike geodesics on $(\mathrm{Cont}(M,\xi), \mathrm{cont}^+(M,\xi),$ $V)$ can be defined as timelike curves (i.e.\ positive paths) which have constant speed and are extremal points of the functional $\mathrm{length}_V$ with respect to variations fixing the end-points. By the bi-invariance and strong concavity of $V$, we again obtain that a positive path $\phi^t$ of contactomorphisms of $(M,\xi)$ is a geodesic if and only if it is autonomous, i.e.\ satisfies
\[
\frac{\mathrm{d}}{\mathrm{d}t} \phi^t = X(\phi^t)
\]
for some time-independent contact vector field $X\in \mathrm{cont}^+(M,\xi)$. 

Since the Reeb vector fields induced by all contact forms defining $\xi$ are precisely the elements of $\mathrm{cont}^+(M,\xi)$, the timelike geodesics in $(\mathrm{Cont}(M,\xi), \mathrm{cont}^+(M,\xi), V)$ are, up to left or right translation by elements of $\mathrm{Cont}(M,\xi)$, precisely the Reeb flows induced by contact forms defining $\xi$.

If $\alpha$ is such a contact form and $\phi^t_{\alpha}$ is the flow of the corresponding Reeb vector field $R_{\alpha}$, then (\ref{LMV}) implies that the Lorentz--Finsler length of the geodesic segment $\{\phi^t_{\alpha}\}_{t\in [0,T]}$ is the quantity
\begin{equation}
\label{Vlengthgeo}
\mathrm{length}_V\bigl( \{\phi^t_{\alpha}\}_{t\in [0,T]} \bigr)= T \, \mathrm{vol}(M,\alpha)^{-\frac{1}{n}}.
\end{equation}

In Section \ref{secondvarsec}, we compute the second variation of the functional  $\mathrm{length}_V$ at a geodesic segment. Unlike in the finite dimensional case of $\mathrm{Sp}(2n)$, this quadratic form has always not only infinite Morse index but also infinite Morse co-index. Indeed, we shall prove the following result.

\begin{mainprop}
\label{infmorind}
Let $\phi^t$ be the Reeb flow of a contact form $\alpha$ defining $\xi$. Then for every $T>0$ the symmetric bilinear form
\[
\mathrm{d}^2 \mathrm{length}_V(\{\phi^t\}_{t\in [0,T]})
\]
is positive definite (resp.\ negative definite) on some infinite dimensional subspace $W^+$ (resp.\ $W^-$) of variations of $\{\phi^t\}_{t\in [0,T]}$ which vanish for $t=0$ and $t=T$. In particular, geodesics in $\mathrm{Cont}(M,\xi)$ are never locally length maximizing nor length minimizing.
\end{mainprop}

Conjugate instants along the Reeb flow $\phi^t$ of a contact form $\alpha$ defining $\xi$ can be defined as usual as the positive numbers $t^*$ such that, after modding out the invariance by reparametrizations, the second variation of $\mathrm{length}_V$ at $\{\phi^t\}_{t\in [0,t^*]}$ has a non-trivial kernel. The elements of this kernel are the Jacobi vector fields vanishing at $t=0$ and $t=t^*$. The equation for these time-dependent contact vector fields $Y$ reads exactly as in the finite dimensional case, i.e.
\[
\partial_{tt} Y = [R_{\alpha},\partial_t Y],
\]
provided that we define the Lie bracket of two vector fields $X,Y$ by the non-standard sign convention 
\[
[X,Y] = - \mathcal{L}_X Y.
\]
Although not standard, this sign convention is quite natural if one wishes to be consistent with the conventions from Lie group theory and is used by some authors, see \cite{arn78} and \cite[Remark 3.1.6]{ms95}. We shall adopt it also here.

The existence of conjugate instants along the geodesic which is determined by the Reeb vector field $R_{\alpha}$ depends on the dynamics of $R_{\alpha}$.  This is illustrated by the explicit computation of all conjugate instants in the following two examples, see Section \ref{secondvarsec}.

\begin{mainex}
\label{exintro1}
{\rm Consider the group $\mathrm{Cont}_0(\T,\xi_0) = \mathrm{Diff}_0(\T)$ of all orientation-preserving diffeomorphisms of $\T=\R/\Z$. Let $X$ be any positive vector field, i.e.\
\[
X(x) = H(x) \frac{\partial}{\partial x},
\]
where $H$ is a positive smooth function on $\T$, and let $\phi^t$ be its flow. Then $t^*>0$ is a conjugate instant along $\phi^t$ if and only if $t^*$ is a positive rational number times $\frac{1}{V(X)}$. Each conjugate instant has infinite multiplicity.}
\end{mainex}

In particular, in the above case conjugate points accumulate at zero. This fact could be used to give an alternative proof of the fact that the second variation of $\mathrm{length}_V$ at any geodesic segment in $\mathrm{Cont}_0(\T,\xi_0)$ has infinite Morse co-index. The latter fact, which as we have seen in the proposition above holds in general, does not require conjugate points accumulating at zero. Indeed, a timelike geodesic in $\mathrm{Cont}(M,\xi)$ may have no conjugate points at all, as shown by the next example.

\begin{mainex}
\label{exintro2}
{\rm Consider the group $\mathrm{Cont}(\T^3,\xi)$ where 
\[
\xi = \ker \alpha \qquad \mbox{and} \qquad \alpha (x,y,z)= \cos (2\pi z) \, \mathrm{d} x + \sin (2\pi z)\, \mathrm{d} y.
\]
The flow of the Reeb vector field of the contact form $\alpha$ is given by
\[
\phi^t(x,y,z) = \bigl( x + t \cos (2\pi z), y + t  \sin (2\pi z), z).
\]
If we identify $\T^3$ with the unit cotangent bundle of $\T^2$, the above flow is precisely the geodesic flow  induced by the flat Euclidean metric on $\T^2 = \R^2/\Z^2$. The geodesic $\phi^t$ in $\mathrm{Cont}(\T^3,\xi)$ has no conjugate instants.}
\end{mainex}

A study of conjugate instants and of the second variation of the $L^2$-energy functional in the context of the group of Hamiltonian diffeomorphisms of a closed symplectic manifold can be found in Vishnevsky's thesis \cite{vis21}. 

\subsection{A systolic question for non-autonomous Reeb flows}
\label{secH}

The {\it systolic ratio} of a contact form $\alpha$ on a closed $(2n-1)$-dimensional manifold $M$ is defined as 
\[
\rho_{\mathrm{sys}}(M,\alpha) := \frac{T_{\min}(\alpha)^n}{\mathrm{vol}(M,\alpha)},
\]
where $T_{\min}(\alpha)$ denotes the minimum over all the periods of closed orbits of the Reeb vector field $R_{\alpha}$. A contact form $\alpha$ is called {\it Zoll} if all the orbits of the corresponding Reeb flow are periodic and have the same minimal period. The main result of \cite{ab19} is that Zoll contact forms are local maximizers of the systolic ratio in the $C^3$-topology of contact forms: Any Zoll contact form $\alpha_0$ on the closed manifold $M$ has a $C^3$-neighborhood $\mathcal{U}$ such that
\begin{equation}
\label{sysin-intro}
\rho_{\mathrm{sys}}(M,\alpha) \leq \rho_{\mathrm{sys}}(M,\alpha_0) \qquad \forall \alpha\in \mathcal{U},
\end{equation}
with the equality holding if and only if $\alpha$ is Zoll. See \cite{apb14}, \cite{abhs18} and \cite{bk21} for previous results on the local systolic optimality of Zoll contact forms and for the relationship with metric systolic geometry. Note that this is a local phenomenon: The systolic ratio is always unbounded from above on the space of contact forms defining a given contact structure, as proven by Sa\u{g}lam in \cite{sag21} generalizing previous results from \cite{abhs18} and \cite{abhs19}.

Here, we would like to discuss whether the local systolic optimality of Zoll Reeb flows extends to non-autonomous Reeb flows. In order to formulate this precisely, recall that a discriminant point of a contactomorphism $\phi\in \mathrm{Cont}(M,\xi)$ is a fixed point $x_0$ of $\phi$ such that the endomorphism $\mathrm{d}\phi(x_0): T_{x_0} M \rightarrow T_{x_0} M$ has determinant one. Equivalently, $x_0$ is a fixed point such that
\[
(\phi_* \alpha) (x_0) = \alpha(x_0)
\]
for some, and hence any, contact form $\alpha$ defining $\xi$. Since the Reeb flow $\phi^t_{\alpha}$ of the contact form $\alpha$ preserves $\alpha$, any point on a $T$-periodic orbit of this flow is a discriminant point for $\phi^T_{\alpha}$. Thanks to the identity (\ref{Vlengthgeo}), the inequality (\ref{sysin-intro}) can then be restated in terms of the Lorentz--Finsler metric $V$ in the following way: Let $X_0=R_{\alpha_0}\in \mathrm{cont}^+(M,\xi)$ be a Zoll Reeb vector field with orbits of minimal period $T_0$. Then there exists $\epsilon>0$ such that for every $X\in \mathrm{cont}^+(M,\xi)$ with $\|X-X_0\|_{C^2} < \epsilon$ the following holds: If the autonomous positive path of contactomorphisms $\phi^t$ given by the flow of $X$ satisfies
\[
\mathrm{length}_V \bigl( \{\phi^t\}_{t\in [0,T_0]} \bigr) \geq \mathrm{length}_V \bigl( \{\phi^t_{\alpha_0}\}_{t\in [0,T_0]} \bigr),
\]
then there exists $t^*\in (0,T_0]$ such that $\phi^{t^*}$ has discriminant points. Our question here is whether this statement remains true for non-autonomous positive paths of contactomorphisms. 

More precisely: Let $\alpha_0$ be a Zoll contact form defining the contact structure $\xi$ on $M$, with Reeb flow $\phi^t_{\alpha_0}$ and minimal period $T_0$. Let $\{\phi^t\}_{t\in [0,T_0]}$ be a positive path in $\mathrm{Cont}(M,\xi)$ such that $\phi^0=\mathrm{id}$. Is it true that if $\{\phi^t\}_{t\in [0,T_0]}$ is suitably close to $\{\phi^t_{\alpha_0}\}_{t\in [0,T_0]}$ and 
\[
\mathrm{length}_V \bigl( \{\phi^t\}_{t\in [0,T_0]} \bigr) \geq \mathrm{length}_V \bigl( \{\phi^t_{\alpha_0}\}_{t\in [0,T_0]} \bigr),
\]
then there exists $t^*\in (0,T_0]$ such that $\phi^{t^*}$ has at least one discriminant point? Or even just a fixed point?

The fact that geodesic arcs in $\mathrm{Cont}(M,\xi)$ are never length maximizing implies that the answer to this question is negative. Indeed, in Section \ref{systolicsec} we shall deduce from Proposition \ref{infmorind} the following result. 

\begin{mainthm} 
\label{systolic}
Let $\alpha_0$ be a Zoll contact form defining the contact structure $\xi$ on $M$, with Reeb flow $\phi^t_{\alpha_0}$ and minimal period $T_0$. Then there exists a smooth 1-parameter family $\{\phi_s\}_{s\in (-\epsilon,\epsilon)}$ of positive paths
\[
\phi_s : [0,T_0] \rightarrow \mathrm{Cont}(M,\xi)
\]
such that $\phi_0(t) = \phi_{\alpha_0}^t$ for every $t\in [0,T_0]$, $\phi_s(0)=\mathrm{id}$ and 
\[
\mathrm{length}_V \bigl( \phi_s \bigr) = \mathrm{length}_V \bigl( \{\phi^t_{\alpha_0}\}_{t\in [0,T_0]} \bigr)
\]
for every $s\in (-\epsilon,\epsilon)$,  but $\phi_s(t)$ has no fixed points for every $s\neq 0$ and every $t\in (0,T_0]$.
\end{mainthm}

\subsection{A time function and a partial order on the universal cover of $\mathrm{Sp}(2n)$} 
\label{secI}

The space $(\mathrm{Sp}(2n),\mathrm{sp}^+(2n))$ is {\it totally vicious}, meaning that it admits closed timelike curves, such as for instance the curve $t\mapsto e^{tJ_0}$, $t\in [0,2\pi]$. Totally viciousness can be avoided if we pass to the universal cover $\widetilde{\mathrm{Sp}}(2n)$ of $\mathrm{Sp}(2n)$, which as usual we think of as the space of homotopy classes $[w]$ of paths $w:[0,1] \rightarrow \mathrm{Sp}(2n)$ starting at the identity; $\widetilde{\mathrm{Sp}}(2n)$ is a Lie group with the same Lie algebra $\mathrm{sp}(2n)$, and the covering map
\[
\pi: \widetilde{\mathrm{Sp}}(2n) \rightarrow \mathrm{Sp}(2n), \qquad [w] \mapsto w(1),
\]
is a homomorphism. The bi-invariant Lorentz--Finsler structure $(\mathrm{sp}^+(2n),G)$ on $\mathrm{Sp}(2n)$ lifts to a bi-invariant Lorentz--Finsler structure on the Lie group $\widetilde{\mathrm{Sp}}(2n)$, for which we keep the same notation. On $(\widetilde{\mathrm{Sp}}(2n),\mathrm{sp}^+(2n))$, there are no closed causal curves. Actually, more is true:

\begin{mainthm}
\label{time}
There exists a time function on $(\widetilde{\mathrm{Sp}}(2n),\mathrm{sp}^+(2n))$, namely a continuous function $f: \widetilde{\mathrm{Sp}}(2n) \rightarrow \R$ that is strictly increasing on every causal curve. Moreover, this time function $f$ can be chosen to be an unbounded quasimorphism.
\end{mainthm}

We recall that a real function $f$ on a group $G$ is said to be a quasimorphism if there is a global bound
\[
|f(v w) - f(v) - f(w)| \leq C \qquad \forall v, w\in G
\]
measuring the failure of $f$ from being a homomorphism. The interesting quasimorphisms are the unbounded ones (every bounded function is trivially a quasimorphism).

This time function is constructed in Section \ref{timesec} starting from a well known function on $\widetilde{\mathrm{Sp}}(2n)$, namely the {\it homogeneous Maslov quasimorphism}
\[
\mu : \widetilde{\mathrm{Sp}}(2n) \rightarrow \R.
\]
This function, which was first defined in \cite{gl58}, is the unique homogeneous real quasi-morphism on $\widetilde{\mathrm{Sp}}(2n)$ whose restriction to $\pi^{-1}(\mathrm{U}(n))$ agrees with the lift of the complex determinant. Here, homogeneity means $\mu(W^k) = k  \mu(W)$ for every integer $k$. Moreover, we are normalizing $\mu$ so that $\mu(\tau(\mathrm{id}))=1$, where $\tau$ is the positive generator of the group of deck transformations of the universal cover of $\mathrm{Sp}(2n)$, or equivalently $\mu(w)=n$ if $w$ is the homotopy class of the loop $\{e^{2\pi t J_0}\}_{t\in [0,1]}$. 

Furthermore, $\mu$ is conjugacy invariant, continuous, and non-decreasing on every causal curve. However, there are causal curves, and even timelike ones, on which $\mu$ is constant, see Lemma \ref{algebraic} below, so $\mu$ is not a time function. Nevertheless, the fact that $\mu$ is strictly increasing on causal curves which are contained in a suitable open subset of $\widetilde{\mathrm{Sp}}(2n)$ allows us to modify it and obtain a time function $f$ as in Theorem \ref{time}. Actually, $f$ can be chosen to be arbitrarily close to $\mu$ with respect to the supremum norm. See Section \ref{timesec} below.

It is worth noticing that this time function cannot be conjugacy invariant: Indeed, no continuous function on $\widetilde{\mathrm{Sp}}(2n)$ which strictly increases on timelike curves can be conjugacy invariant, see Proposition \ref{nonconju} below.

The existence of a time function implies that there are no closed causal curves on $(\widetilde{\mathrm{Sp}}(2n), \mathrm{sp}^+(2n))$. The latter fact is equivalent to the fact that the relation
\[
\mathcal{J} := \{ (w_0, w_1) \in \widetilde{\mathrm{Sp}}(2n)^2 \mid \mbox{either there is a causal curve from } w_0 \mbox{ to } w_1 \mbox{ or } w_0=w_1\}
\]
is a partial order on $\widetilde{\mathrm{Sp}}(2n)$. We shall use the notation $w_0 \leq w_1$ as shorthand for $(w_0, w_1)\in \mathcal{J}$, and $w_0\geq w_1$ as synonymous of $w_1\leq w_0$. In general relativity, cone structures satisfying the latter condition are called {\it causal}, while the existence of a time function is equivalent to a stronger condition called {\it stable causality}. See e.g., \cite[Chapter 3]{ms08} or  \cite[Chapter 4]{min19}. Thanks to the bi-invariance of the cone distribution, this partial order gives $\widetilde{\mathrm{Sp}}(2n)$ the structure of a partially ordered group in the sense of \cite{fuc63}.

\subsection{The Lorentz distance on the universal cover of $\mathrm{Sp}(2n)$} 
\label{secJ}

Let $(K,F)$ be a Lorentz--Finsler structure on the manifold $M$. We assume that the cone distribution $K$ is causal and denote by $\leq$ the corresponding partial order relation on $M$. The Lorentz--Finsler metric $F$ induces the {\it Lorentz distance}
\[
\mathrm{dist}_F: M \times M \rightarrow [0,+\infty], \qquad \mathrm{dist}_F(p,q) :=  \left\{ \begin{array}{ll} \sup  \mathrm{length}_F(\gamma) & \mbox{if }  p \leq  q, \\ 0 & \mbox{otherwise}, \end{array} \right.
\]
where the supremum is taken over all causal curves $\gamma:[0,1] \rightarrow M$ such that $\gamma(0)=p$ and $\gamma(1)=q$. Note that this function may be trivial, meaning that it takes only the values $0$ and $+\infty$. In general relativity, Lorentz distances are also called {\it time separation functions}. The Lorentz distance $\mathrm{dist}_F$  is lower semicontinuous and satisfies the reverse triangular inequality
\[
\mathrm{dist}_F (p_0,p_2) \geq \mathrm{dist}_F(p_0,p_1) + \mathrm{dist}_F(p_1,p_2),\qquad \mbox{if } p_0 \leq p_1 \leq p_2.
\]
See e.g., \cite[Section 2.9]{min19}. 

In this section, we discuss some properties of the Lorentz distance $\mathrm{dist}_G$ on the universal cover of $\mathrm{Sp}(2n)$, which as we have seen is causal. The bi-invariance of $G$ implies that $\mathrm{dist}_G$ is also bi-invariant, and hence it suffices to study $\mathrm{dist}_G(\mathrm{id},w)$ for $w\in \widetilde{\mathrm{Sp}}(2n)$ with $w\geq \mathrm{id}$.

In order to state our result, we need to recall some notions from Krein theory (see e.g., \cite{ys75} or \cite{eke90}). By extending the skew-symmetric bilinear form $\omega_0$ to a skew-Hermitian form on $\C^{2n} \times \C^{2n}$ and by multiplying it by $-i$, we obtain the Hermitian form 
\[
\kappa: \C^{2n} \times \C^{2n} \rightarrow \C, \qquad \kappa  := -i \omega_0,
\]
which has signature $(n,n)$ and is known as {\it Krein form}. The eigenvalues of $W\in \mathrm{Sp}(2n)$ that are either real or lie on the unit circle 
\[
\U:= \{z\in \C \mid |z|=1\}
\]
occur in pairs $\lambda,\lambda^{-1}$, while all other eigenvalues occur in quadruples $\lambda, \overline{\lambda}, \lambda^{-1}, {\overline{\lambda}}^{-1}$.
The restriction of $\kappa$ to the generalized eigenspace of an eigenvalue $\lambda$ on $\U$ is always non-degenerate, and $\lambda$ is said to be {\it Krein-positive} (resp.\ {\it Krein-negative}) if this restriction is positive (resp.\ negative) definite. If $\lambda\in \U$ is Krein-positive, then $\lambda^{-1}=\overline{\lambda}$ is Krein-negative. The {\it positively elliptic region}s is the set
\[
\begin{split}
\mathrm{Sp}_{\mathrm{ell}}^+ (2n) := \{ W\in \mathrm{Sp}(2n)  \mid & \mbox{ all the eigenvalues of } W \mbox{ are in } \U\setminus \{\pm 1\} \mbox{ and the ones}\\ & \mbox{ with positive imaginary part are Krein-positive} \}.
\end{split}
\]  
Equivalently, $\mathrm{Sp}_{\mathrm{ell}}^+ (2n)$ can be described as the set of linear symplectomorphisms of $\R^{2n}$ which split into $n$ rotations of angles in the interval $(0,\pi)$: more precisely, $W$ is in $\mathrm{Sp}_{\mathrm{ell}}^+ (2n)$ if and only if $W=e^X$ where $X\in \mathrm{sp}^+(2n)$ is as in (\ref{isdiag}) with $\theta_j\in (0,\pi)$ for every $j$ (see Proposition \ref{propA2} in Appendix \ref{linear-app}).

One of the fundamental results of Krein theory is that Krein-definite eigenvalues on $\U$ are stable, meaning that they cannot leave $\U$ after a perturbation. This implies that $ \mathrm{Sp}_{\mathrm{ell}}^+ (2n)$ is open in $\mathrm{Sp}(2n)$.

In the case $n=1$, elements $W$ of $\mathrm{Sp}(2)$ have the polar decomposition $W=UP$, where $U\in \mathrm{SO}(2)\cong S^1$ and $P$ is symmetric, positive definite and symplectic. Since any such $P$ is the exponential of a unique element in $\mathrm{sp}(2)\cap \mathrm{Sym}(2)\cong \R^2$, $\mathrm{Sp}(2)$ is homeomorphic to $\R^2 \times S^1$, or equivalently to $\mathbb{D}\times S^1$, where  $\mathbb{D}$ is the open disk in $\R^2$. 

\begin{figure}[h]
\centerline{\hspace{.2cm}\scalebox{.28}{\includegraphics{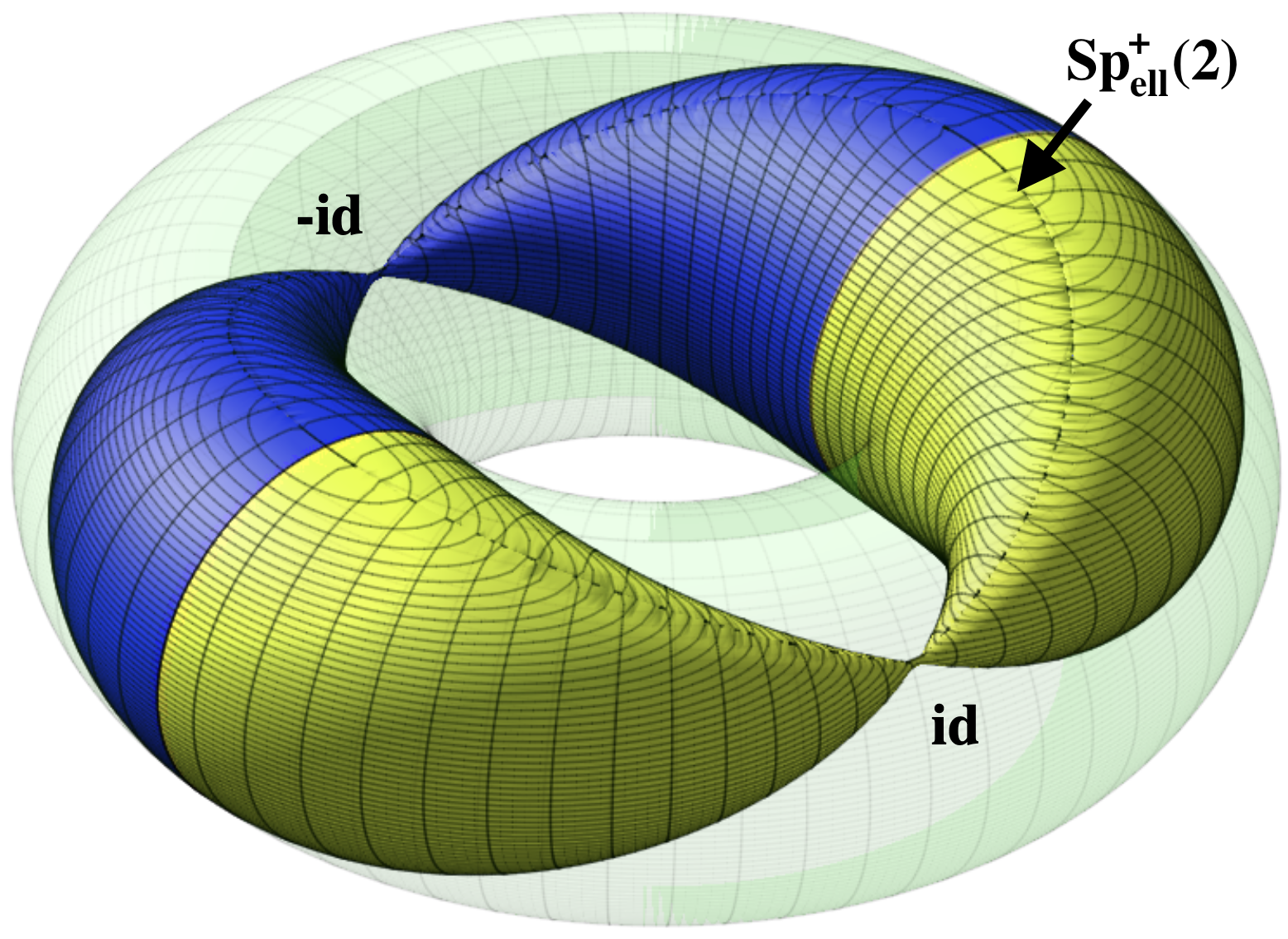}} \hspace{1cm} \scalebox{.45}{\includegraphics{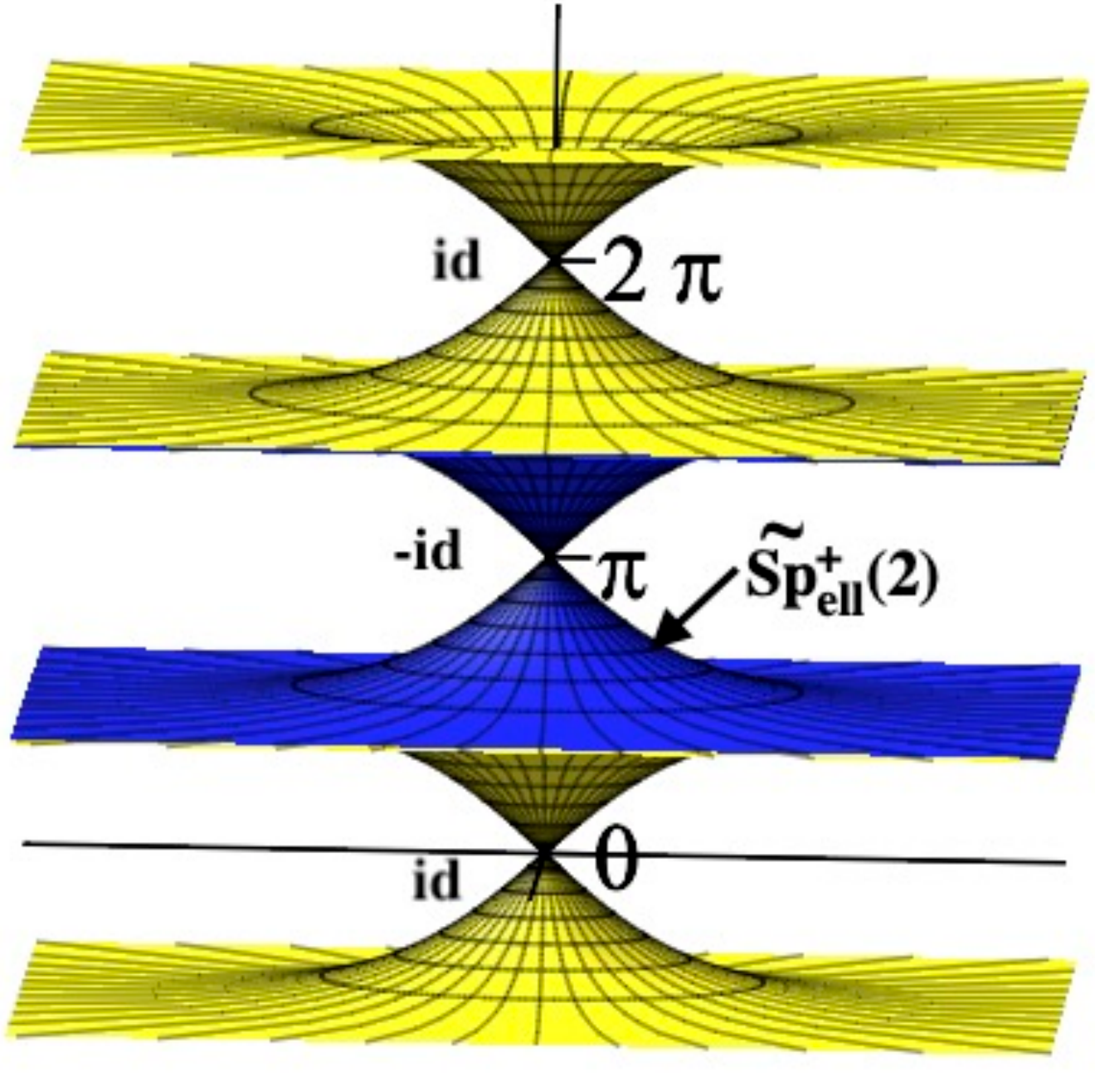}}}
\caption{On the left: The union of the sets $\Sigma_1$ and $\Sigma_{-1}$ decomposes $\mathrm{Sp}(2)$ into four open domains, one of which is $\mathrm{Sp}_{\mathrm{ell}}^+ (2)$. On the right: The lift of these sets to $\widetilde{\mathrm{Sp}}(2)\cong \R^2\times \R$.}
\label{figura}
\end{figure}

In the picture on the left in Figure \ref{figura}, we visualize $\mathrm{Sp}(2)\cong \mathbb{D} \times S^1$ as an open region in $\R^3$ bounded by a torus-like surface. The circle sitting at the core of this region (not represented in the picture) corresponds to the subgroup $\mathrm{SO}(2)=\mathrm{U}(1)$. The yellow double cone emanating from the identity represents the discriminant $\Sigma_1$, i.e., the set of $W$'s in $\mathrm{Sp}(2)$ having the eigenvalue 1. The blue double cone emanating from minus the identity is the set $\Sigma_{-1}$ of elements in $\mathrm{Sp}(2)$ having the eigenvalue $-1$ (these two surfaces seem to intersect in the picture, but their intersection is on the boundary of the region, which is not part of $\mathrm{Sp}(2)$). The open region bounded by the ``croissant'' on the upper right part is precisely the positively elliptic region $\mathrm{Sp}_{\mathrm{ell}}^+ (2)$, the region symmetric to it is the set consisting of $W$'s having both eigenvalues in $\U\setminus \{\pm 1\}$ with the Krein-positive one having negative imaginary part. The two outer regions correspond to $W$'s with either positive (region adherent to the identity) or negative (region adherent to minus the identity) real eigenvalues. In the companion picture on the right, we are visualizing a portion of the universal cover $\widetilde{\mathrm{Sp}}(2)$ as $\R^2 \times \R$. The lift of $\mathrm{SO}(2)$ is now the vertical axis, and the yellow and blue surfaces represent the lifts of the sets $\Sigma_1$ and $\Sigma_{-1}$. See \cite[Section 1.2.1 and 1.2.2]{abb01} for the explicit parametrizations leading to these pictures. 

A fundamental feature of timelike curves, or equivalently solutions of non-autonomous positive definite linear Hamiltonian systems, is that Krein-positive eigenvalues on $\U$ move counterclockwise, while Krein-negative ones move clockwise, see \cite[Proposition I.3.2 and Corollary I.3.3]{eke90}. Therefore, any timelike curve starting at the identity immediately enters $\mathrm{Sp}_{\mathrm{ell}}^+ (2n)$ and can leave it only if the eigenvalue $-1$ appears. In the case $n=1$, referring again to the picture on the left in Figure \ref{figura}, we have that timelike curves starting at the identity immediately enter the upper-right ``croissant'' and can leave it only through the blue surface. See also Section \ref{Sp(2)sec} below for another picture representing $\widetilde{\mathrm{Sp}}(2)$ and explicit coordinates on $\widetilde{\mathrm{Sp}}(2)$ which simplify the study of causality on this space. 

We denote by $\widetilde{\mathrm{Sp}}_{\mathrm{ell}}^+ (2n)$ the open subset of the universal cover of $\mathrm{Sp}(2n)$ consisting of homotopy classes of paths $w:[0,1] \rightarrow \mathrm{Sp}(2n)$ such that $w(0)=\mathrm{id}$ and $w((0,1]) \subset \mathrm{Sp}_{\mathrm{ell}}^+ (2n)$. Equivalently, $\widetilde{\mathrm{Sp}}_{\mathrm{ell}}^+ (2n)$ is the connected component of $\pi^{-1} ( \mathrm{Sp}_{\mathrm{ell}}^+ (2n))$ whose closure contains the identity. See again the right picture in Figure \ref{figura}.

The next result shows that the Lorentz distance $\mathrm{dist}_{G}$ on $\widetilde{\mathrm{Sp}}(2n)$ is non-trivial but is not everywhere finite either.

\begin{mainthm}
\label{Lordist}
Let $w$ be an element of $\widetilde{\mathrm{Sp}}(2n)$. Then we have:
\begin{enumerate}[(i)]
\item If $w$ is in the closure of $\widetilde{\mathrm{Sp}}_{\mathrm{ell}}^+ (2n)$, then $\pi(w)$ has spectrum $\{e^{\pm i \theta_1}, \dots, e^{\pm i \theta_n} \}$ with $\theta_j\in [0,\pi]$ for every $j=1,\dots,n$ and
\begin{equation}
\label{upperbound}
\mathrm{dist}_G(\mathrm{id},w) \leq \frac{2\pi}{n} \mu(w) = \frac{1}{n}   \sum_{j=1}^n \theta_j.
\end{equation}
\item If there is a timelike curve from the identity to $w$ and $w$ is not in the closure of $\widetilde{\mathrm{Sp}}_{\mathrm{ell}}^+ (2n)$, then $\mathrm{dist}_G(\mathrm{id},w)= +\infty$.
\end{enumerate}
\end{mainthm}

The function $\mu$ appearing in (i) is the homogeneous Maslov quasimorphism discussed in Section \ref{secH} above. This linear algebra result plays an important role in the length bounds of contactomorphism or symplectomorphism groups which we discuss below in Sections \ref{secM} and  \ref{secN}.

\begin{mainrem}{\rm Note that the length of the unique timelike geodesic from $\mathrm{id}$ to $w\in \widetilde{\mathrm{Sp}}_{\mathrm{ell}}^+ (2n)$ is the geometric mean of the numbers $\theta_j$ appearing in (i) above. This implies that in the special case in which $\theta_1=\dots=\theta_n\in (0,\pi)$, the inequality in (\ref{upperbound}) is actually an equality and the Lorentz distance $\mathrm{dist}_G(\mathrm{id},w)$ is achieved by the unique geodesic from $\mathrm{id}$ to $w$. We believe that the latter fact is true for any $w$ in $\widetilde{\mathrm{Sp}}_{\mathrm{ell}}^+ (2n)$. This is related to Question \ref{globhyp} below.}
\end{mainrem}

\begin{mainrem}{\rm (Long timelike curves not hitting the discriminant) Given $\lambda\in \C$, set
\[
\Sigma_{\lambda} := \{W\in \mathrm{Sp}(2n) \mid \lambda \mbox{ is an eigenvalue of } W \}.
\]
Statement (i) in the above theorem implies that any timelike curve $W:[0,1]\rightarrow \mathrm{Sp}(2n)$ with $W(0)= \mathrm{id}$ and $\mathrm{length}_G(W)\geq \pi$ must hit the set $\Sigma_{-1}$. On the other hand, statement (ii) implies that there are timelike curves $W:[0,1]\rightarrow \mathrm{Sp}(2n)$ with $W(0)= \mathrm{id}$ and arbitrarily large $\mathrm{length}_G(W)$ which never hit the discriminant $\Sigma_1$ after $t=0$. Note that this is a non-autonomous phenomenon: if $W(t)=e^{tX}$ with $X\in \mathrm{sp}^+(2n)$ satisfies $\mathrm{length}_G(W|_{[0,1]})\geq 2\pi$, then there exists $t^*\in (0,1]$ such that $W(t^*)\in \Sigma_1$. Indeed, since $G(X)=\mathrm{length}_G(W|_{[0,1]})\geq 2\pi$, identity (\ref{geomean}) implies that at least one of the numbers $\theta_j$ in (\ref{isdiag}) is at least $2\pi$, and we deduce that $W(t^*) = e^{t^*X}$ has the eigenvalue 1 for $t^*=\frac{2\pi}{\theta_j} \in (0,1]$. 
}
\end{mainrem}

\begin{mainrem}
\label{elli} {\rm
If we restrict the Lorentz--Finsler structure $(\mathrm{sp}^+(2n),G)$ of $\mathrm{Sp}(2n)$ to the open subset $\mathrm{Sp}^+_{\mathrm{ell}}(2n)$, we also obtain a stably causal space. The Lorentz distance on this space is everywhere finite and defines the structure of a Lorentzian length space in the sense of \cite{ks18}. The above theorem implies that this space has diameter $\pi$. In the case $n=1$, this space is {\it globally hyperbolic}, meaning that for every pair of points $w_0$, $w_1$ in it the set of all $w\in \mathrm{Sp}^+_{\mathrm{ell}}(2)$ with $w_0\leq w \leq w_1$ is compact, see Section \ref{Sp(2)sec} below. Global hyperbolicity is an important notion in general relativity. It has other equivalent characterizations, such as for instance the existence of a {\it Cauchy hypersurface}, i.e., a hypersurface which is met exactly once by every inextensible causal curve, and some striking consequences, such as the existence of a timelike geodesic between any two points which can be connected by a timelike curve and the well-posedness of the Cauchy problem for the wave equation. See  \cite[Section 3.11]{ms08}) and references therein. Therefore, we state the following:}
\end{mainrem}

\begin{mainque}
\label{globhyp}
Is $\mathrm{Sp}^+_{\mathrm{ell}}(2n)$ globally hyperbolic also for $n>1$?
\end{mainque}

Theorem \ref{Lordist} is proven in Section \ref{dimLordist} below. In the preceding Section \ref{Sp(2)sec}, we look more closely at the case $n=1$, i.e., at the case of the universal cover of the three-dimensional anti-de Sitter space $\mathrm{AdS}_3$. In this case, the Lorentz distance is completely described by Proposition \ref{lordisSp2}, which implies Theorem \ref{Lordist} for $n=1$. Moreover, the fact that elements $w\geq \mathrm{id}$ which are not in the closure of $\widetilde{\mathrm{Sp}}^+_{\mathrm{ell}}(2)$ have infinite distance from the identity is not specific of the Lorentz distance induced by $G$ and holds for any Lorentz distance which is conjugacy invariant, see Proposition \ref{trivlordisSp} below.

\subsection{The Lorentz distance on the universal cover of $\mathrm{Cont}_0 (\R \mathrm{P}^{2n-1},\xi_{\mathrm{st}})$}
\label{secK}

A general fact about causality in $(\mathrm{Cont}(M,\xi), \mathrm{cont}^+(M,\xi))$ is that the existence of a non-constant non-negative loop, i.e., a closed causal curve, implies the existence of a positive loop, i.e., a closed timelike curve, see \cite[Proposition 2.1.B]{ep00}. Besides this, causality depends on the contact manifold $(M,\xi)$ under consideration. There are contact manifolds $(M,\xi)$ such that $\mathrm{Cont}_0(M,\xi)$ admits no  positive loop whatsoever: This is the case of the cotangent sphere bundle $ST^*Q$ of any closed manifold $Q$ having infinite fundamental group, see \cite[Section 9]{cn10b}. At the opposite end of the spectrum, there are contact manifolds $(M,\xi)$ such that $\mathrm{Cont}_0(M,\xi)$ admits even contractible positive loops, such as $(S^{2n-1},\xi_{\mathrm{st}})$ for $n\geq 2$, see \cite{ekp06}. 

In the middle, there are contact manifolds such that $\mathrm{Cont}_0(M,\xi)$ admits positive loops but no contractible ones. This is the case of $\mathrm{Cont}_0(\R \mathrm{P}^{2n-1},\xi_{\mathrm{st}})$ for every $n\geq 1$. Indeed, the standard Reeb flow on $(\R \mathrm{P}^{2n-1},\xi_{\mathrm{st}})$ defines a non-contractible positive loop, but there are no contractible positive loops in $\mathrm{Cont}_0(\R \mathrm{P}^{2n-1},\xi_{\mathrm{st}})$. As shown in \cite{ep00}, this follows from the existence of Givental's  {\it asymptotic nonlinear Maslov index}
\[
\nu: \widetilde{\mathrm{Cont}}_0(\R \mathrm{P}^{2n-1},\xi_{\mathrm{st}}) \rightarrow \R
\]
on the universal cover of  $\mathrm{Cont}_0 (\R \mathrm{P}^{2n-1},\xi_{\mathrm{st}})$.
This is a conjugacy invariant homogeneous quasimorphism and is continuous with respect to the topology which is induced by the $C^0$-topology on the space of Hamiltonians. Moreover, if we endow $ \widetilde{\mathrm{Cont}}_0(\R \mathrm{P}^{2n-1},\xi_{\mathrm{st}})$ with the bi-invariant cone distribution that is induced by $\mathrm{cont}^+ (\R \mathrm{P}^{2n-1},\xi_{\mathrm{st}})$, we obtain that $\nu$ is non-decreasing along each non-positive path and it is strictly positive on each element of $\widetilde{\mathrm{Cont}}_0(\R \mathrm{P}^{2n-1},\xi_{\mathrm{st}})$ which is the end-point of a positive path starting at the identity. The asymptotic nonlinear Maslov index $\nu$ extends the homogeneous Maslov quasimorphism $\mu: \widetilde{\mathrm{Sp}} (2n) \rightarrow \R$, meaning that $\mu = \nu \circ \tilde{\jmath}$, where $\tilde{j}: \widetilde{\mathrm{Sp}}(2n) \rightarrow \widetilde{\mathrm{Cont}}_0(\R \mathrm{P}^{2n-1},\xi_{\mathrm{st}})$ is the lift of the homomorphism $j$ from  Section \ref{secE} above.

The situation of the contactomorphism group of the real projective space is then analogous to what we have encountered with the linear symplectic group and we get a genuine partial order $\leq$ on $\widetilde{\mathrm{Cont}}_0 (\R \mathrm{P}^{2n-1},\xi_{\mathrm{st}})$, where $\phi_0 \leq \phi_1$ means that there exists a non-negative path from $\phi_0$ to $\phi_1$. In the language of general relativity, $\widetilde{\mathrm{Cont}}_0 (\R \mathrm{P}^{2n-1},\xi_{\mathrm{st}})$ is then a causal space. A natural question is whether the stronger property of Theorem \ref{time} holds also for $\widetilde{\mathrm{Cont}}_0 (\R \mathrm{P}^{2n-1},\xi_{\mathrm{st}})$:

\begin{mainque}
Does $\widetilde{\mathrm{Cont}}_0 (\R \mathrm{P}^{2n-1},\xi_{\mathrm{st}})$ admit a time function, i.e., a real function that strictly increases along every non-negative path, which is continuous with respect to some reasonable topology?
\end{mainque} 

Being non-decreasing on non-negative paths, the nonlinear asymptotic Maslov index $\nu$ seems to be a good starting point to build a time function on $\widetilde{\mathrm{Cont}}_0 (\R \mathrm{P}^{2n-1},\xi_{\mathrm{st}})$. However, $\nu$ needs to be corrected, since it can be constant on some positive paths.
In the proof of Theorem \ref{time}, we can correct $\mu$ building on the fact that this function is strictly increasing along every non-negative path which is contained in a certain non-empty open subset of $\widetilde{\mathrm{Sp}}(2n)$, but it is not clear to us whether $\nu$ shares this property.

In the case $n=1$, the answer to the above question is positive. Indeed, since 
\[
\mathrm{Cont}(\R \mathrm{P}^1,\xi_{\mathrm{st}})\cong \mathrm{Cont}(\T,\xi_0) = \mathrm{Diff}_0(\T), 
\]
the universal cover $\widetilde{\mathrm{Cont}}_0 (\R \mathrm{P}^1,\xi_{\mathrm{st}})$ can be identified with the group $\mathrm{Diff}_1(\R)$ of diffeomorphisms $\phi: \R \rightarrow \R$ such that
\[
\phi(x+1) = 1 + \phi(x) \qquad \forall x\in \R,
\]
or, equivalently, diffeomorphisms $\phi: \R \rightarrow \R$ of the form $\phi=\mathrm{id} + \phi_0$ with $\phi_0:\R \rightarrow \R$ 1-periodic. The order $\leq$ on $\widetilde{\mathrm{Cont}}_0 (\R \mathrm{P}^1,\xi_{\mathrm{st}})$ corresponds to the standard order on real valued functions, and the function
\[
f: \mathrm{Diff}_1(\R) \rightarrow \R, \qquad f(\phi) := \int_0^1 (\phi(x)-x)\, dx,
\]
is readily seen to be a time function.

We now lift the Lorentz--Finsler metric $V$ to $\widetilde{\mathrm{Cont}}_0 (\R \mathrm{P}^{2n-1},\xi_{\mathrm{st}})$, and denote it by the same symbol.
This Lorentz--Finsler metric induces the Lorentz distance $\mathrm{dist}_V$ on $\widetilde{\mathrm{Cont}}_0 (\R \mathrm{P}^{2n-1},\xi_{\mathrm{st}})$. Unlike the Lorentz distance $\mathrm{dist}_{G}$ on $\widetilde{\mathrm{Sp}}(2n)$, the Lorentz distance $\mathrm{dist}_V$ is trivial on $\widetilde{\mathrm{Cont}}_0(\R \mathrm{P}^{2n-1},\xi_{\mathrm{st}})$. Actually, the following stronger result holds.

\begin{mainthm}
\label{nolordisCont}
Let $d: \widetilde{\mathrm{Cont}}_0 (\R \mathrm{P}^{2n-1},\xi_{\mathrm{st}}) \times \widetilde{\mathrm{Cont}}_0 (\R \mathrm{P}^{2n-1},\xi_{\mathrm{st}}) \rightarrow [0,+\infty]$ be a function such that:
\begin{enumerate}[(i)]
\item $d(\phi_0,\phi_1)>0$ if and only if there is a non-negative and somewhere positive path from $\phi_0$ to $\phi_1$;
\item $d(\phi_0,\phi_2) \geq d(\phi_0,\phi_1) + d(\phi_1,\phi_2)$ if $\phi_0\leq \phi_1 \leq \phi_2$;
\item $d$ is bi-invariant.
\end{enumerate}
Then $d(\phi_0,\phi_1)$ has the value $+\infty$ if there is a non-negative and somewhere positive path from $\phi_0$ to $\phi_1$, and $0$ otherwise.
\end{mainthm}

The result will come to no surprise to experts in contact geometry. Indeed, finding meaningful bi-invariant ``global measurements'' on contactomorphism groups is a notoriously difficult problem. As recalled at the beginning of this introduction, any bi-invariant distance function on the contactomorphism group $\mathrm{Cont}(M,\xi)$ is discrete. Therefore, $\mathrm{Cont}(M,\xi)$ does not admit a bi-invariant distance function coming from a Finsler metric, unlike the symplectomorphism group, whose Hofer metric is bi-invariant and is induced by a genuine Finsler metric (see \cite{hof93} and \cite{pol01}). As first shown by Sandon in \cite{san10}, non-trivial discrete bi-invariant distance functions do exist on some contactomorphism groups, see also \cite{zap13,cs15,san15,fpr18}. If one drops the requirement of being bi-invariant, there do exist interesting Lorentz distances on  orderable contactomorphisms groups, as recently shown by Hedicke in \cite{hed22}. 

Our proof of Theorem \ref{nolordisCont} is based on statement (ii) in Theorem \ref{Lordist} and will be carried out in Section \ref{nolordisContsec}. It is reasonable to believe that an analogous result holds for every orderable contactomorphisms group. Since we do not have a proof of this fact, we formulate the following:

\begin{mainque}
Does Theorem \ref{nolordisCont} extend to all orderable contactomorphisms groups?
\end{mainque}

Actually, we do not even know whether the Lorentz distance which is induced by the Lorentz--Finsler metric $V$ is trivial on every orderable contactomorphisms group.

\subsection{Length bounds in $\mathrm{Diff}_1(\R)$} 
\label{secL}

Consider again the 1-dimensional contact manifold $(\T=\R/\Z,\xi_0=\{0\})$. As discussed above, the universal cover of $\mathrm{Cont}(\T,\xi_0)=\mathrm{Diff}_0(\T)$ can be identified with the group $\mathrm{Diff}_1(\R)$ of diffeomorphisms $\phi: \R \rightarrow \R$ such that $\phi(x+1)=1+\phi(x)$ for every $x\in \R$, and the order $\leq$ on $\mathrm{Diff}_1(\R)$ is just the standard order on real valued functions.

Let $\phi\in \mathrm{Diff}_1(\R)$ be such that $\phi(x)>x$ for every $x\in \R$. Theorem \ref{nolordisCont} in the case $n=1$ tells us that  there are arbitrarily long positive paths $\{\phi^t\}_{t\in [0,1]}$ from $\mathrm{id}$ to $\phi$. See also Example \ref{leonid} for an explicit construction of arbitrarily long positive paths starting at the identity and staying below a fixed translation. Denote by $H$ the Hamiltonian which is associated to such a path: $H$ is a smooth positive function on $[0,1]\times \R$, 1-periodic in the second variable, and $\phi^t$ solves the ODE
\[
\frac{\mathrm{d}}{\mathrm{d}t} \phi^t = H(t,\phi^t), \qquad \phi^0=\mathrm{id}.
\]
A natural question is how ``complex'' $H$ must be, so that the path $\{\phi^t\}_{t\in [0,1]}$ connecting $\mathrm{id}$ to $\phi$ is long with respect to the Lorentz--Finsler metric $V$. 

A first observation is that $H$ must be non-autonomous. Indeed, if $H>0$ is autonomous and its flow $\phi^t$ satisfies $\phi^1(x) \leq \phi(x)$ for every $x\in \R$, then by Proposition \ref{rotnumb} we have 
\[
\mathrm{length}_V(\{\phi^t\}_{t\in [0,1]}) = \rho(\phi^1) \leq \rho(\phi),
\]
where 
\[
\rho: \mathrm{Diff}_1(\R) \rightarrow \R, \qquad \rho(\phi)= \lim_{n\rightarrow \infty} \frac{\phi^n(x) - x}{n}
\]
denotes the translation number quasi-morphism. 

In our next result, we consider the number of harmonics of the 1-periodic functions $x\mapsto H(t,x)$ as a measure of the complexity of $H$. 
Denoting by $\mathcal{P}_k$ the space of smooth functions on $[0,1]\times \R$ which are trigonometric polynomials of degree at most $k$ in the second variable, i.e., functions of the form
\[
p(t,x) = \sum_{j=0}^k \bigl(a_j(t) \cos (2\pi j x) + b_j(t) \sin (2\pi j x)  \bigr),
\]
for suitable smooth functions $a_j$, $b_j$, we have the following result.

\begin{mainthm}
\label{quantum}
For every $k\in \N$ and $\phi\in \mathrm{Diff}_1(\R)$  the following facts hold:
\begin{enumerate}[(i)]
\item If 
\[
\phi(x) \leq x + \frac{s}{4k} \qquad \forall x\in \R
\]
for some $s\in (0,1)$, then for every positive path $\{\phi^t\}_{t\in [0,1]}$ from $\mathrm{id}$ to $\phi$ which is generated by a time-dependent Hamiltonian in $\mathcal{P}_k$ we have
\[
\mathrm{length}_V (\{\phi^t\}_{t\in [0,1]}) \leq \frac{1}{1-s} \int_0^1 \bigl( \phi(x) - x \bigr)\, \mathrm{d}x.
\]
\item  If 
\[
\phi(x)> x + \frac{1}{k}  \qquad \forall x\in \R,
\]
then there exist $\psi\in \mathrm{Diff}_1(\R)$ with $\mathrm{id} < \psi < \phi$ and positive paths from $\mathrm{id}$ to $\psi$ which are generated by time-dependent Hamiltonians in $\mathcal{P}_k$ and have arbitrarily large $\mathrm{length}_V$.
\end{enumerate}
\end{mainthm}

\begin{mainrem}\label{quantumrem} 
{\rm Note the ``quantum'' nature of this result: If $\phi<\mathrm{id} + \frac{1}{4k}$ then positive paths from $\mathrm{id}$ to $\phi$ generated by Hamiltonians in $\mathcal{P}_k$ have uniformly bounded length, while if $\phi>\mathrm{id} + \frac{1}{k}$ is the end point of a positive path starting from $\mathrm{id}$ and generated by Hamiltonians in $\mathcal{P}_k$, then the length of such a path can be arbitrarily large. An interesting and presumably non-trivial question is how to close the gap between the thresholds $\frac{1}{4k}$ and $\frac{1}{k}$. }
\end{mainrem}

The proof of (i) uses a sharp Bernstein type inequality for non-negative periodic functions which is due to Nazarov, see Theorem \ref{bernstein} below. The proof of (ii) uses the embeddings $j_k$ of the classical Lorentzian spacetime $\mathrm{PSp}(2)$ into $\mathrm{Diff}_0(\T)$, see Proposition \ref{sottogruppi}.

In the context of the $L^2$-metric on the group of Hamiltonian diffeomorphisms of a compact symplectic manifold, a related phenomenon has been studied in the already mentioned \cite{vis21}.

\subsection{Length bounds on the universal cover of $\mathrm{Cont}_0(\R \mathrm{P}^{2n-1},\xi_{\mathrm{st}})$}
\label{secM}

We do not know whether the quantum phenomenon which is described by Theorem \ref{quantum} and Remark \ref{quantumrem} above holds also on contact manifolds of dimension larger than one. In particular, we would like to state the following:

\begin{mainque}
Does Statement (i) of Theorem \ref{quantum} generalize to $\widetilde{\mathrm{Cont}}_0(\R \mathrm{P}^{2n-1},\xi_{\mathrm{st}})$?
\end{mainque}

Here, it would be natural to replace the space $\mathcal{P}_k$ with the space of spherical harmonics of degree at most $k$ and the upper bound on $\phi$ with the condition that $\phi$ should satisfy $\phi\leq e_{c_k}$ for some $c_k>0$, where $e_t$ is the element of $\widetilde{\mathrm{Cont}}_0(\R \mathrm{P}^{2n-1},\xi_{\mathrm{st}})$ which is defined by the restriction to the interval $[0,t]$ of the standard periodic Reeb flow.

Another question which arises naturally (see also Remark \ref{condimpr} below) and we do not know how to answer concerns the lift
\[
\tilde{\jmath}: \widetilde{\mathrm{Sp}} (2n) \rightarrow  \widetilde{\mathrm{Cont}}_0(\R \mathrm{P}^{2n-1},\xi_{\mathrm{st}})
\]
of the homomorphism
\[
j: \mathrm{PSp} (2n) \rightarrow  \mathrm{Cont}_0(\R \mathrm{P}^{2n-1},\xi_{\mathrm{st}})
\]
from Section \ref{secE} above. The fact that the pull-back by $j$ of the cone $\mathrm{cont}^+(\R \mathrm{P}^{2n-1},\xi_{\mathrm{st}})$ is the cone $\mathrm{sp}^+(2n)$ implies that if $w_0\leq w_1$ in $\widetilde{\mathrm{Sp}} (2n)$, then $\tilde{\jmath}(w_0) \leq \tilde{\jmath}(w_1)$ in $\widetilde{\mathrm{Cont}}_0(\R \mathrm{P}^{2n-1},\xi_{\mathrm{st}})$. We do not know whether the converse is also true:

\begin{mainque}
\label{converse}
Assume that $w_0,w_1\in \widetilde{\mathrm{Sp}} (2n)$ satisfy $\tilde{\jmath}(w_0) \leq \tilde{\jmath}(w_1)$ in $\widetilde{\mathrm{Cont}}_0(\R \mathrm{P}^{2n-1},\xi_{\mathrm{st}})$. Is it true that $w_0\leq w_1$ in $\widetilde{\mathrm{Sp}} (2n)$?
\end{mainque}

In order to study this question, it is natural to consider the set of all $w$ in $\widetilde{\mathrm{Sp}} (2n)$ such that $\tilde{\jmath}(w) \geq \mathrm{id}$ in $\widetilde{\mathrm{Cont}}_0(\R \mathrm{P}^{2n-1},\xi_{\mathrm{st}})$, and the question is whether this set coincides with the set of $w$ in $\widetilde{\mathrm{Sp}} (2n)$ such that $w\geq \mathrm{id}$. Since both these sets are conjugacy invariant semigroups in $\widetilde{\mathrm{Sp}} (2n)$ and the second one is contained in the first one, Question \ref{converse} has a positive answer if one can show that the conjugacy invariant semigroup which is defined by the partial order $\leq$ on $\widetilde{\mathrm{Sp}} (2n)$ is a maximal proper conjugacy invariant semigroup. This fact is true for $n=1$, see \cite[Section 3.3]{bsh12}, leading to the positive answer to Question \ref{converse} in this case. See also Proposition \ref{converse1} below for a simpler argument.

To state the length bound that we can prove for positive paths in $\widetilde{\mathrm{Cont}}_0(\R \mathrm{P}^{2n-1},\xi_{\mathrm{st}})$, we need to introduce some notation. 
On $(\R \mathrm{P}^{2n-1},\xi_{\mathrm{st}})$, we fix the standard contact form $\alpha_0$, whose Reeb flow is Zoll with period $\frac{\pi}{2}$, and we use $\alpha_0$ to identify the space of Hamiltonians $C^{\infty}(\R \mathrm{P}^{2n-1})$ with the space of contact vector fields $\mathrm{cont}(\R \mathrm{P}^{2n-1},\xi_{\mathrm{st}})$. Note that there is a one-to-one correspondence between functions  $H\in C^{\infty}(\R \mathrm{P}^{2n-1})$ and 2-homogeneous even functions on $\R^{2n}$, which is obtained by identifying $\R \mathrm{P}^{2n-1}$ with the quotient of $S^{2n-1}$ by the antipodal $\Z_2$-action and by extending even functions on $S^{2n-1}$ to $\R^{2n}$ by 2-homogeneity. 

By a {\it positive quadratic function} on $\R \mathrm{P}^{2n-1}$ we mean a function which corresponds to a positive definite quadratic form on $\R^{2n}$ under this identification. Given a real number $c\geq 1$, we define $\mathcal{H}_c$ to be the space of time-dependent Hamiltonians $H\in C^{\infty}([0,1]\times \R \mathrm{P}^{2n-1})$ such that
\[
Q \leq H \leq c \, Q
\]
for some $Q\in C^{\infty}([0,1]\times \R \mathrm{P}^{2n-1})$ such that $Q(t,\cdot)$ is a positive quadratic function for every $t\in [0,1]$. Note that $\mathcal{H}_c$ is a convex cone and that the union of all $\mathcal{H}_c$ for $c\geq 1$ is the convex cone of all positive time-dependent Hamiltonians.

Note also that any $H\in C^{\infty}([0,1]\times \R \mathrm{P}^{2n-1})$ which is {\it convex}, meaning that its 2-homogeneous extension to $[0,1]\times \R^{2n}$ is convex in the second variable and positive on $[0,1]\times (\R^{2n}\setminus \{0\})$, belongs to $\mathcal{H}_{2n}$. This follows from John's theorem, stating that if $E$ is the ellipsoid of maximal volume which is contained in a centrally symmetric convex body $C\subset \R^N$, then $E\subset C \subset NE$. We can now state our next result.

\begin{mainthm}
\label{lengthboundn}
Let $\phi\in \widetilde{\mathrm{Cont}}_0(\R \mathrm{P}^{2n-1},\xi_{\mathrm{st}})$ be such that
\begin{equation}
\label{ipotesi}
\nu(\phi) \leq \frac{1}{2},
\end{equation}
with strict inequality in the case $n=1$. Then every positive path $\{\phi_t\}_{t\in [0,1]}$ from $\mathrm{id}$ to $\phi$ which is generated by a Hamiltonian in $\mathcal{H}_c$ satisfies
\[
\mathrm{length}_V\bigl(\{\phi_t\}_{t\in [0,1]} \bigr) \leq \frac{2^{\frac{1}{n}}}{n} c \, \nu (\phi).
\]
\end{mainthm}

This theorem is deduced in Section \ref{lengthboundnnsec} from Theorem \ref{Lordist} (i) above. 
The length bound of Theorem \ref{lengthboundn} is of a different nature than the one of Theorem \ref{quantum}: In Theorem \ref{quantum}, the paths of Hamiltonians are constrained to finite dimensional spaces but can take values in the whole cone of positive functions within these spaces, whereas in Theorem \ref{lengthboundn} there is no finite-dimensional constraint but the cone of positive functions is reduced by a pinching condition. 
In both results, the time dependence of the Hamiltonians can be arbitrary. 
 
\begin{mainrem}
{\rm
The upper bound $\frac{1}{2}$ in assumption (\ref{ipotesi}) of Theorem \ref{lengthboundn} is optimal: For every $s> \frac{1}{2}$ (and for $n=1$ also for $s=\frac{1}{2}$), there exists $\phi\in \widetilde{\mathrm{Cont}}_0(\R \mathrm{P}^{2n-1},\xi_{\mathrm{st}})$ with $\nu(\phi)=s$ such that there are positive paths from $\mathrm{id}$ to $\phi$ which are generated by Hamiltonians in $\mathcal{H}_1$ and have arbitrarily large length. See Remark \ref{optimal} below. }
\end{mainrem}

\begin{mainrem}
\label{condimpr}
{\rm
If the answer to Question \ref{converse} is positive, then the same conclusion of Theorem \ref{lengthboundn} holds also replacing (\ref{ipotesi}) by the assumption $\phi\leq e$, where $e$ denotes the element of $\widetilde{\mathrm{Cont}}_0(\R \mathrm{P}^{2n-1},\xi_{\mathrm{st}})$ generated by the constant Hamiltonian $H=1$, i.e., the homotopy class of the loop $\{\phi^{\frac{\pi}{2}t}_{\alpha_0}\}_{t\in [0,1]}$, where $\phi_{\alpha_0}$ denotes the Reeb flow of $\alpha_0$ (note that $\nu(e) = \frac{n}{2}$). See Remark \ref{alternative} below.

Since Question  \ref{converse} has a positive answer for $n=1$, in this case the length bound of Theorem \ref{lengthboundn} holds for every $\phi\leq e$. In the identification $\widetilde{\mathrm{Cont}}_0(\R \mathrm{P}^1,\xi_{\mathrm{st}})=\mathrm{Diff}_1(\R)$, we have $\mathcal{H}_1=\mathcal{P}_1$, $e(x) = x +1$ and the asymptotic nonlinear Maslov index $\nu$ coincides with $\frac{1}{2} \rho$, where $\rho(\phi)$ denotes the translation number of $\phi$. In this setting, the length bound of Theorem \ref{lengthboundn} holds for every $\phi\in \mathrm{Diff}_1(\R)$ satisfying $\phi(x)\leq x + 1$ for every $x\in \R$, which is a weaker condition than the assumption $\rho(\phi)<1$ from Theorem \ref{lengthboundn}. Note also that the optimality of the bound (\ref{ipotesi}) mentioned in the above remark now says that there  are diffeomorphisms $\phi\in \mathrm{Diff}_1(\R)$ with $\rho(\phi)=1$ which are the end-points of positive paths with arbitrarily large length which are generated by Hamiltonians in $\mathcal{P}_1$.
}
\end{mainrem}

We conclude this section by discussing length bounds for {\it autonomous positive paths} in the group $\widetilde{\mathrm{Cont}}_0(\R \mathrm{P}^{2n-1},\xi_{\mathrm{st}})$, i.e., solutions of
\[
\frac{\mathrm{d}}{\mathrm{d}t} \phi (t)= X(\phi(t)), \qquad \phi(0) = \mathrm{id},
\]
with $X\in \mathrm{cont}^+(M,\xi)$ independent of time. In the case $n=1$, Proposition \ref{rotnumb} and the identity $\rho=2\nu$ imply
\begin{equation}
\label{lengthnu}
\mathrm{length}_V(\phi|_{[0,1]}) = 2 \, \nu(\phi(1))
\end{equation}
for every autonomous positive path $\phi$ in $\widetilde{\mathrm{Cont}}_0(\R \mathrm{P}^1,\xi_{\mathrm{st}})$. In higher dimension, we certainly do not have this identity, but one may wonder whether the Lorentz--Finsler length of $\phi|_{[0,1]}$ for  $\phi$ an autonomous positive path in $\widetilde{\mathrm{Cont}}_0(\R \mathrm{P}^{2n-1},\xi_{\mathrm{st}})$ can be bounded from above in terms of $\nu(\phi(1))$. Note in fact that the Lorentz--Finsler length of the autonomous path 
\[
w(t)=e^{tX} \in \widetilde{\mathrm{Sp}}(2n)
\]
with $X$ as in (\ref{isdiag}) has the bound
\[
\mathrm{length}_G(w|_{[0,1]}) =  \left(\prod_{j=1}^n \theta_j \right)^{\frac{1}{n}} \leq  \frac{1}{n} \sum_{j=1}^n \theta_j = \frac{2\pi}{n} \, \mu(w(1)),
\]
thanks to the inequality between the geometric and arithmetic mean.
However, in the infinite dimensional group $\widetilde{\mathrm{Cont}}_0(\R \mathrm{P}^{2n-1},\xi_{\mathrm{st}})$ with $n>1$ the length of autonomous positive paths $\phi$ does not have an upper bound in terms of $\nu(\phi(1))$, as we now explain.

Indeed, let us consider the prequantization $S^1$-bundle $p:\R \mathrm{P}^{2n-1} \rightarrow \C\mathrm{P}^{n-1}$ and the moment map $m:\C \mathrm{P}^{n-1} \rightarrow \R^n$ which is associated to the standard Hamiltonian $\mathbb{T}^n$-action on $\C\mathrm{P}^{n-1}$ endowed with the Fubini--Study symplectic form. The image of $m$ is a simplex $\Delta$. Let $H\in C^{\infty}(\R \mathrm{P}^{2n-1})$ be a function which lifts a function $h: \Delta \rightarrow \R$ under the map $m\circ p$, and let $\phi: \R \rightarrow \widetilde{\mathrm{Cont}}_0(\R \mathrm{P}^{2n-1},\xi_{\mathrm{st}})$ be the path which is generated by $H$, seen as a contact Hamiltonian with respect to the standard contact form $\alpha_0$ of $\R \mathrm{P}^{2n-1}$. It can be proven that in this case 
\[
\nu(\phi(1)) = \frac{n}{2} \, h(b),
\]
where $b$ denotes the barycenter of $\Delta$. See \cite[Theorem 1.11]{ep09} and \cite{bs07}. By choosing a positive function $h$ with $h(b)$ small and $h$ equal to a large constant outside of a small neighborhood of $b$, we can make 
\[
\mathrm{length}_V(\phi|_{[0,1]}) = \left( \int_{\R \mathrm{P}^{2n-1}} H^{-n} \, \alpha_0 \wedge \mathrm{d} \alpha_0^{n-1} \right)^{-\frac{1}{n}}
\]
arbitrarily large and keep $\nu(\phi(1))$ arbitrarily small, as claimed above.

As we have already noted in  Section \ref{secI}, if we fix an element $\psi$ in $\widetilde{\mathrm{Cont}}_0(\R \mathrm{P}^1,\xi_{\mathrm{st}})$, then the Lorentz--Finsler length of an arbitrary autonomous positive path $\phi:[0,1] \rightarrow \widetilde{\mathrm{Cont}}_0(\R \mathrm{P}^1,\xi_{\mathrm{st}})$ such that $\phi(1) \leq \psi$ is uniformly bounded from above. The above example does not exclude that this holds true also in higher dimension. Therefore, we state the following question.

\begin{mainque}
Let $\psi\in \widetilde{\mathrm{Cont}}_0(\R \mathrm{P}^{2n-1},\xi_{\mathrm{st}})$. Is the Lorentz--Finsler length of autonomous positive paths $\phi: [0,1] \rightarrow \widetilde{\mathrm{Cont}}_0(\R \mathrm{P}^{2n-1},\xi_{\mathrm{st}})$ such that $\phi(1)\leq \psi$ uniformly bounded from above?
\end{mainque}

\subsection{Positive paths in the group of symplectomorphisms of uniformly convex domains} 
\label{secN}

Consider a bounded uniformly convex open set $\Omega$ with smooth boundary in the standard symplectic vector space $(\R^{2n},\omega_0)$. In other words, $\Omega$ is a non-empty open sublevel of a coercive smooth function on $\R^{2n}$ whose second differential is everywhere positive definite. Denote by $\mathrm{Symp}_0(\overline{\Omega})$ the identity component of the group of symplectomorphisms of the closure of $\Omega$. Let us emphasize that each element of $\mathrm{Symp}_0(\overline{\Omega})$ keeps the boundary $\partial \Omega$ invariant, but in general induces a non-trivial diffeomorphism of $\partial \Omega$. Any $\phi\in \mathrm{Symp}_0(\overline{\Omega})$ is the time-one map of a Hamiltonian vector field, i.e.\ $\phi=\phi_H^1$ where $\phi^t_H$ is the solution of the Cauchy problem
\begin{equation}
\label{hameq}
\frac{\mathrm{d}}{\mathrm{d}t} \phi_H^t = X_{H_t}(\phi_H^t), \qquad \phi_H^0 = \mathrm{id},
\end{equation}
where the Hamiltonian $H\in C^{\infty}([0,1]\times \overline{\Omega})$ is such that $H_t=H(t,\cdot)$ is constant on each leaf of the characteristic foliation of $\partial \Omega$. Here, $X_K$ denotes the Hamiltonian vector field of the function $K\in C^{\infty} (\overline{\Omega})$, which is defined by the identity
\[
\omega_0(X_K,\cdot) = -\di K,
\]
and the characteristic foliation of the hypersurface $\partial \Omega$ is the one-dimensional foliation which is tangent to the kernel of the restriction of $\omega_0$ to the tangent spaces of $\partial \Omega$. Both (\ref{hameq}) and the boundary condition are not affected by adding a function of $t$ to the Hamiltonian, and we obtain that the Lie algebra of $\mathrm{Symp}_0(\overline{\Omega})$ can be identified with the vector space
\[
\mathcal{H}(\overline{\Omega}) := \{ H\in C^{\infty} ( \overline{\Omega}) \mid H\; \mbox{\small is constant on each leaf of the characteristic foliation of } \partial \Omega \}/\R,
\]
where the quotient is with respect to the action of $\R$ which is given by adding constant functions. With a small abuse of notation, we shall see equivalence classes in $\mathcal{H}(\overline{\Omega})$ as functions on $\overline{\Omega}$.

It will be useful to use the standard Euclidean structure of $\R^{2n}$ and the standard complex structure $J_0$ and rewrite (\ref{hameq}) as
\begin{equation}
\label{hameq2}
\frac{\mathrm{d}}{\mathrm{d}t} \phi_H^t = J_0 \nabla H_t(\phi_H^t), \qquad \phi_H^0 = \mathrm{id},
\end{equation}
but all the notions we are introducing in this section depend only on the affine symplectic structure of $\R^{2n}$.

In this section, we discuss a structure on $\mathrm{Symp}_0(\overline{\Omega})$ which is induced from the Lorentz--Finsler metric $G$ on $\mathrm{Sp}(2n)$ in the following way. By linearizing (\ref{hameq2}) and using the identification $T\overline\Omega = \overline{\Omega} \times \R^{2n}$, we obtain for every $z\in \overline{\Omega}$ a path $t \mapsto \mathrm{d}\phi_H^t(z)$ in $\mathrm{Sp}(2n)$ which satisfies
\[
\frac{\mathrm{d}}{\mathrm{d}t} \mathrm{d}\phi_H^t (z) = J_0 \nabla^2 H_t(\phi_H^t(z)) \, \mathrm{d} \phi_H^t (z), \qquad \mathrm{d}\phi_H^0 (z)= \mathrm{id}.
\]
This path is timelike in $(\mathrm{Sp}(2n),\mathrm{sp}^+(2n))$ for every $z\in \overline{\Omega}$ if and only if $H_t$ is uniformly convex on $\overline{\Omega}$ for every $t\in [0,1]$, meaning that
$\nabla^2 H_t(z)$ is positive definite for every $z\in \overline{\Omega}$. This suggests to consider the following open convex cone in $\mathcal{H}(\overline{\Omega})$
\[
\mathcal{H}^+(\overline{\Omega}) := \{ H\in \mathcal{H}(\overline{\Omega}) \mid H \mbox{ uniformly convex on } \overline{\Omega} \},
\]
which is non-empty thanks to the uniform convexity of $\Omega$, 
and the smooth function $\mathcal{G}: \mathcal{H}^+(\overline{\Omega}) \rightarrow (0,+\infty)$ which is given by
\[
\mathcal{G}(H) := \frac{1}{\mathrm{vol}(\Omega)} \int_{\Omega} G(J_0 \nabla^2 H(z))\, \mathrm{d}z = \frac{1}{\mathrm{vol}(\Omega)} \int_{\Omega}\bigl( \det \nabla^2 H(z) \bigr)^{\frac{1}{2n}} \, \mathrm{d}z,
\]
where $\mathrm{vol}(\Omega)$ denotes the Euclidean volume of $\Omega$ and $\mathrm{d}z$ refers to integration with respect to the Euclidean volume form $\frac{1}{n!} \omega_0^n$ of $\R^{2n}$. 

The function $\mathcal{G}$ is smooth on $\mathcal{H}^+(\overline{\Omega})$ and
extends continuously to the closure of $\mathcal{H}^+(\overline{\Omega})$; notice that this extension is not identically zero on the boundary. Moreover, $\mathcal{G}$ satisfies the strong concavity condition
\[
\mathrm{d}^2\mathcal{G}(H) \cdot (K,K) < 0 \qquad \forall H\in \mathcal{H}^+(\overline{\Omega}), \; \forall K \in \mathcal{H}(\overline{\Omega})\setminus \R H,
\]
as shown in Proposition \ref{Mconcave} below. 

The cone $\mathcal{H}^+(\overline{\Omega})$ and the function $\mathcal{G}$ extend to a cone distribution on $\mathrm{Symp}_0(\overline{\Omega})$ and a function on it by using right-shifts on the Lie group. The resulting objects are not bi-invariant, because $\mathcal{H}^+(\overline{\Omega})$ and $\mathcal{G}$ are not invariant under the adjoint action of $\mathrm{Symp}_0(\overline{\Omega})$ on $\mathcal{H}(\overline{\Omega})$. This is due to the fact that we have used the affine structure of $\R^{2n}$ in order to identify tangent spaces at different points when linearizing (\ref{hameq}). These objects are nevertheless equivariant with respect to the affine symplectic group of $(\R^{2n},\omega_0)$.

The resulting structure on $\mathrm{Symp}_0(\overline{\Omega})$, which we still denote by $(\mathcal{H}^+(\overline{\Omega}),\mathcal{G})$, satisfies all the requirements of a Lorentz--Finsler structure as in Definition \ref{deflorfin}, except for the condition that $\mathcal{G}$ should vanish on the boundary of the cones. Therefore, we call $\mathcal{G}$ a \textit{weak} Lorentz--Finsler metric on the cone distribution determined by $\mathcal{H}^+(\overline{\Omega})$.

The Lorentz--Finsler length of any positive (i.e.\ timelike) path in $\mathrm{Symp}_0(\overline{\Omega})$ is still defined and denoted as usual by $\mathrm{length}_{\mathcal{G}}$.

As an example, take $\Omega$ to be the unit Euclidean ball in $\R^{2n}$ and consider the subalgebra $\mathfrak{u}\subset \mathcal{H}(\overline{\Omega})$ consisting of the quadratic Hamiltonians of the form
\[
H(z) = \frac{1}{2} Sz\cdot z,
\]
where $S\in \mathrm{Sym}(2n)$ commutes with $J_0$. The set of endomorphisms of $\R^{2n}\cong \C^n$ of the form $J_0S$ with $S$ as above is precisely the Lie algebra of $\mathrm{U}(n)$, so the Hamiltonian flow of $H$ is unitary. If $H$ is in $\mathfrak{u} \cap \mathcal{H}^+(\overline{\Omega})$, i.e.\ if $S$ as above is positive definite, then $\phi=\{e^{tJ_0S}\}_{t\in [0,1]}$ is an autonomous  positive path in both $\mathrm{Symp}_0(\overline{\Omega})$ and $\mathrm{Sp}(2n)$ and we have
\[
\mathrm{length}_{\mathcal{G}}(\phi) = \mathcal{G}(H) = G(J_0 S) = (\det S)^{\frac{1}{2n}} = \mathrm{length}_G(\phi).
\]
In general, it is easy to see that $\mathrm{length}_{\mathcal{G}}$ and $\mathrm{length}_G$ are related by the identity
\begin{equation}
\label{GGversusG}
\mathrm{length}_{\mathcal{G}}(\{t\mapsto \phi^t\}) = \frac{1}{\mathrm{vol}(\Omega)} \int_{\Omega} \mathrm{length}_G(\{t\mapsto \mathrm{d} \phi^t(z)\})\, \mathrm{d}z,
\end{equation}
for every positive path $\phi: [0,1] \rightarrow \mathrm{Symp}_0(\overline{\Omega})$, see Proposition \ref{propGGversusG} below.

We denote by $\widetilde{\mathrm{Symp}}_0(\overline{\Omega})$ the universal cover of $\mathrm{Symp}_0(\overline{\Omega})$. The homogeneous Maslov quasimorphism $\mu: \widetilde{\mathrm{Sp}}(2n)\rightarrow \R$ extends to a real valued quasimorphism on $\widetilde{\mathrm{Symp}}_0(\overline{\Omega})$ by setting
\[
\mathcal{M}(\tilde{\phi}) := \frac{1}{\mathrm{vol}(\Omega)} \int_{\Omega} \mu([\mathrm{d}\phi(z)])\,  \mathrm{d}z.
\]
Here, $\tilde{\phi}\in \widetilde{\mathrm{Symp}}_0(\overline{\Omega})$ is the homotopy class of a path $\{\phi^t\}_{t\in [0,1]}$ in $\mathrm{Symp}_0(\overline{\Omega})$ with $\phi^0=\mathrm{id}$ and $[\mathrm{d}\phi(z)]\in \widetilde{\mathrm{Sp}}(2n)$ denotes the homotopy class of the path $\{\mathrm{d}\phi^t(z)\}_{t\in [0,1]}$ in $\mathrm{Sp}(2n)$. This quasimorphism, which appears in Ruelle's work \cite{rue85}, was investigated by Barge and Ghys in \cite[Theorem 3.4]{bg92}. Thanks to (\ref{GGversusG}), Theorem \ref{Lordist} (i) has the following consequence, which is proven in Section \ref{Mproofs}.

\begin{mainthm}
\label{corti}
Let $\tilde{\phi}\in \widetilde{\mathrm{Symp}}_0(\overline{\Omega})$ be the homotopy class of a positive path $\phi=\{\phi^t\}_{t\in [0,1]}$ in $\mathrm{Symp}_0(\overline{\Omega})$ with $\phi^0=\mathrm{id}$ and such that $\mathrm{d}\phi^t(z)$ does not have the eigenvalue $-1$ for every $z\in \overline{\Omega}$ and $t\in [0,1]$. Then any positive path $\psi$ in $\mathrm{Symp}_0(\overline{\Omega})$ which is homotopic to $\phi$ with fixed ends satisfies the same condition and
\[
\mathrm{length}_{\mathcal{G}} (\psi) \leq \frac{2\pi}{n} \mathcal{M}(\tilde{\phi}).
\]
\end{mainthm}

The above results imply that the Lorentz distance $\mathrm{dist}_{\mathcal{G}}$ on $\widetilde{\mathrm{Symp}}_0(\overline{\Omega})$ which is induced by the lift of $\mathcal{G}$ is non-trivial, because
\[
0 < \mathrm{dist}_{\mathcal{G}} (\mathrm{id}, \tilde{\phi}) \leq \frac{2\pi}{n} \mathcal{M}(\tilde{\phi}) < +\infty
\]
for any $\tilde{\phi}$ which satisfies the assumptions of the above theorem. Due to Theorem \ref{Lordist} (ii), one may suspect that there are elements $\tilde{\phi}$ in $ \widetilde{\mathrm{Symp}}_0(\overline{\Omega})$ such that $\mathrm{dist}_{\mathcal{G}} (\mathrm{id}, \tilde{\phi})=+\infty$. However, we do not have a proof of this fact and hence we formulate this as a question.

\begin{mainque}
Are there elements $\tilde{\phi}$ in $ \widetilde{\mathrm{Symp}}_0(\overline{\Omega})$ for which the $\mathrm{length}_{\mathcal{G}}$ of positive paths representing the homotopy class $\tilde{\phi}$ has no upper bound?
\end{mainque}

Next, we focus on the following \textit{optimal extension problem}. A smooth path $\psi=\{\psi^t\}_{t\in [0,1]}$ of diffeomorphisms 
\[
\psi^t: \partial \Omega \rightarrow \partial \Omega
\]
is called \textit{extendable} if there exists a positive path $\{\phi^t_K\}_{t\in [0,1]}$ in $\mathrm{Symp}_0(\overline{\Omega})$, where $K\in C^{\infty}([0,1]\times \overline{\Omega})$ and $K_t\in \mathcal{H}^+(\overline{\Omega})$ for every $t\in [0,1]$, such that
\[
\psi^t:= \phi_K^t|_{\partial \Omega}.
\]
Given such an extendable path $\psi$, we consider the variational problem
\begin{equation}
\label{varprob}
 \sup \bigl\{ \mathrm{length}_{\mathcal{G}} (\phi) \mid \phi \mbox{ positive path in } \mathrm{Symp}_0(\overline{\Omega}) \mbox{ such that } \phi^t|_{\partial \Omega} =\psi^t \; \forall t\in [0,1] \bigr\}.
\end{equation}
Equivalently, the Hamiltonian $H\in C^{\infty}([0,1]\times \overline{\Omega})$ generating the path $\phi$ should be uniformly convex and satisfy
\begin{equation}
\label{firstjet}
\nabla H_t(z) = \nabla K_t(z) \qquad \forall (t,z)\in [0,1]\times \partial \Omega.
\end{equation}
Our next result is the finiteness of (\ref{varprob}). In fact, we are able to provide a constructive upper bound. In order to describe it, notice that by the convexity of $\Omega$ and the uniform convexity of $H_t$ the maps
\[
\nabla H_t : \overline{\Omega} \rightarrow \R^{2n}
\]
are embeddings and their image depends only on their restriction to $\partial \Omega$, so by  (\ref{firstjet}) only on the path of diffeomorphisms $\psi^t: \partial \Omega \rightarrow \partial \Omega$. The positive quantity
\[
\mathcal{V}(\psi) := \frac{1}{\mathrm{vol}(\Omega)^{\frac{1}{2n}}} \int_0^1  \mathrm{vol} (\nabla H_t(\Omega))^{\frac{1}{2n}} \, \di t = \frac{1}{\mathrm{vol}(\Omega)^{\frac{1}{2n}}} \int_0^1  \mathrm{vol} (X_{H_t}(\Omega))^{\frac{1}{2n}} \, \di t
\]
is then a function of the path $\psi$. The proof of the next result is discussed in Section \ref{Mproofs}.

\begin{mainthm}
\label{varthm}
For every extendable path of diffeomorphisms $\psi$ of $\partial \Omega$ and every positive path $\phi$ in $\mathrm{Symp}_0(\overline{\Omega})$ extending $\psi$, we have the upper bound
\[
\mathrm{length}_{\mathcal{G}} (\phi) \leq \mathcal{V}(\psi).
\]
The equality holds if and only $\phi$ is generated by a uniformly convex Hamiltonian $H$ satisfying the Monge--Amp\`ere equation
\[
\det \nabla^2 H_t(z) = c(t) \qquad \forall (t,z)\in [0,1]\times \overline{\Omega}
\]
for some positive numbers $c=c(t)$. 
\end{mainthm}

\begin{mainrem} 
{\rm Let us illustrate the quantity $\mathcal{V}(\psi)$ appearing in Theorem \ref{varthm} in the following situation. Assume that $\Omega$ is centrally symmetric, so in particular it contains the origin. Recall that the group $\R_+$ acts on $\R^{2n}$ by dilations $z\mapsto \sqrt{c} z$. If the Hamiltonians $H_t$ in $\mathcal{H}(\overline{\Omega})$ satisfy $H_t(\sqrt{c} z) = c H_t(z)$ near $\partial \Omega$, then the restriction of the corresponding Hamiltonian path $\{\phi_H^t\}_{t\in [0,1]}$ preserves the contact form $\alpha_{\Omega}$ on $\partial \Omega$ which is given by the restriction of the Liouville form $\lambda_0$ (see Section \ref{secE}). Let $H$ be the unique $\R^+$-equivariant function taking the value $\frac{1}{2}$ on $\partial \Omega$. Denote by $\psi=\{\psi^t\}_{t\in [0,1]}$ the path of diffeomorphisms of $\partial \Omega$ which is given by the restriction of the Hamiltonian path induced by $H$. The function $H$ is smooth and uniformly convex on $\R^{2n}\setminus \{0\}$, but in general not twice differentiable at the origin; however, $H$ can be modified near the origin to make it everywhere smooth and uniformly convex, so the path $\psi$ is extendable. The set $\nabla H(\Omega)$ is in this case the \textit{polar body} $\Omega^{\circ}$ of $\Omega$ and hence
\[
\mathcal{V}(\psi) = \mathrm{vol}(\Omega)^{-\frac{1}{2n}}\, \mathrm{vol} (\Omega^{\circ})^{\frac{1}{2n}}.
\]
Note now that $\psi$ is the Reeb flow of the contact form $\alpha_{\Omega}$ on the time interval $[0,1]$ and hence a positive path in $\mathrm{Cont}(\partial \Omega,\xi)$, where $\xi=\ker \alpha_{\Omega}$. Therefore, its length with respect to the Lorentz--Finsler metric $V$ from  Section \ref{secD} is
\[
\mathrm{length}_V(\psi) = \mathrm{vol}(\partial \Omega,\alpha_{\Omega} \wedge d\alpha_{\Omega}^{n-1})^{-\frac{1}{n}} = \mathrm{vol}(\Omega,\omega_0^n)^{-\frac{1}{n}} = \bigl( n! \, \mathrm{vol}(\Omega) \bigr)^{-\frac{1}{n}},
\]
by Stokes theorem, and we obtain the identity
\[
\mathcal{V}(\psi) \, \mathrm{length}_V(\psi)^{-1} = (n!)^{\frac{1}{n}} \bigl( \mathrm{vol}(\Omega)\, \mathrm{vol}(\Omega^0) \bigr)^{\frac{1}{2n}}.
\]
The quantity in brackets on the right-hand side of the above identity is the \textit{Mahler volume} of $\Omega$, a linear invariant of  centrally symmetric convex domains admitting the bounds
\[
\nu^{2n} \frac{4^{2n}}{(2n)!} \leq \mathrm{vol}(\Omega)\, \mathrm{vol}(\Omega^0) \leq \frac{\pi^{2n}}{(n!)^2}.
\]
Here, the upper bound is sharp and is given by the Blaschke--Santal\'o inequality, stating that the Mahler volume is maximized by ellipsoids (see \cite{bla17} and \cite{san49}). The value of the optimal positive number $\nu$ appearing in the above lower bound is not known, but is conjectured to be $\nu=1$. Indeed, the Mahler conjecture, which for now has been proven only in dimension at most three (see \cite{mah39} and \cite{is20}), states that the cube is a minimizer of the Mahler volume (the number $\frac{4^{2n}}{(2n!)}$ is precisely the Mahler volume of the cube in dimension $2n$). The best known bound for $\nu$ is due to Kuperberg, who in \cite{kup08} showed that $\nu\geq \frac{\pi}{4}$. Therefore, $\mathcal{V}(\psi)$ provides the following dimension-independent lower and upper bounds for the Lorentz--Finsler length of $\psi$:
\[
\frac{1}{\pi} \mathcal{V}(\psi) \leq \mathrm{length}_V(\psi) \leq \frac{1}{4\nu} \binom{2n}{n}^{\frac{1}{2n}} \mathcal{V}(\psi) \leq \frac{1}{2\nu}  \mathcal{V}(\psi) \leq \frac{2}{\pi}  \mathcal{V}(\psi),
\]
where we have used the inequality $\binom{2n}{n} \leq 2^{2n}$ and Kuperberg's bound for $\nu$. Together with Theorem \ref{varthm}, the left inequality implies that if $\phi: [0,1] \rightarrow \mathrm{Symp}(\overline{\Omega})$ is any positive path extending the path $\psi$ given by the restriction to the interval $[0,1]$ of the Reeb flow of $\alpha_{\Omega}$ on $\partial \Omega$, then
\[
\mathrm{length}_{\mathcal{G}}(\phi) \leq \pi \, \mathrm{length}_V(\psi).
\]
It would be interesting to explore further connections between the Lorentz--Finsler lengths on the group of symplectomorphisms of a convex domain and on the group of contactomorphisms of its boundary. }
\end{mainrem}

By the theory of the Dirichlet problem for the Monge--Amp\`ere equation, it is easy to produce examples of extendable paths of diffeomorphisms $\psi$ of $\partial \Omega$ admitting an extension $\phi$ such that $\mathrm{length}_{\mathcal{G}} (\phi) = \mathcal{V}(\psi)$. Indeed, let $h:\partial \Omega \rightarrow \R$ be an arbitrary smooth function which is constant on each leaf of the characteristic foliation of $\partial \Omega$. Then the Dirichlet problem
\[
\left\{ \begin{array}{ll} \det \nabla^2 H = 1 & \mbox{ on } \Omega, \\ H=h & \mbox{ on } \partial \Omega, \end{array} \right.
\]
has a unique uniformly convex solution $H\in C^{\infty}(\overline{\Omega})$, see \cite[Theorem 6.2.6 and Proposition 6.1.4]{han16}. Let $\psi=\{\psi^t\}_{t\in [0,1]}$ be the path of diffeomorphisms of $\partial \Omega$ which is given by the boundary restriction of the positive path $\phi = \{\phi_H^t\}_{t\in [0,1]}$. By Theorem \ref{varthm}, $\phi$ is the unique maximizer of the optimal extension problem (\ref{varprob}) and $ \mathrm{length}_{\mathcal{G}}(\phi)=  \mathcal{V}(\psi)$.

For a general extendable path $\psi$, we do not expect the existence of a positive path $\{\phi^t_H\}_{t\in [0,1]}$ in $\mathrm{Symp}_0(\overline{\Omega})$ whose generating Hamiltonian $H$ satisfies the above Monge--Amp\`ere equation. Indeed, the Monge--Amp\`ere equation
\[
\det \nabla^2 H = \mathrm{const} \;  \mbox{ on } \Omega,
\]
with the boundary condition (\ref{firstjet}) defines an overdetermined problem.

It is therefore natural to ask about existence and uniqueness of maximizers of the optimal extension problem (\ref{varprob}) for a general extendable path $\psi$. As proven in Section \ref{Mproofs} below, uniqueness is guaranteed by the strong concavity of $\mathcal{G}$.

\begin{mainprop}
\label{smooth-uniqueness}
Maximizers of the optimal extension problem (\ref{varprob}) are unique.
\end{mainprop}

Existence is a more difficult question.
As a non-essential simplification, consider the optimal extension problem for an autonomous path. Thus, we are given a function $K\in \mathcal{H}^+(\overline{\Omega})$ and we are looking for maximizers of the functional
\[
\mathcal{F}(H) := \int_{\Omega} \bigl( \det \nabla^2 H(z) \bigr)^{\frac{1}{2n}}\, \mathrm{d}z
\]
over the set of all uniformly convex functions $H\in C^{\infty}(\overline{\Omega})$ satisfying (\ref{firstjet}). This problem is invariant under the sum of constant functions, but we can mod this invariance out and obtain an equivalent problem by replacing the boundary condition (\ref{firstjet}) by
\begin{equation}
\label{firstjet2}
H(z) = K(z), \quad \nabla H(z) = \nabla K(z) \qquad \forall (t,z)\in [0,1]\times \partial \Omega.
\end{equation}
The analogous functional with exponent $\frac{1}{2n+2}$ instead of $\frac{1}{2n}$ (in dimension $2n$) is called \textit{affine area} of the hypersurface which is given by the graph of $H$, and the corresponding variational problem with boundary conditions (\ref{firstjet2}) is the \textit{first boundary value problem for affine maximal hypersurfaces}, which is discussed in detail by Trudinger and Wang in \cite{tw05} and \cite[Section 6.4]{tw08}.

As we explain in Section \ref{Mproofs}, some of the analysis of Trudinger and Wang goes through also for the functional $\mathcal{F}$ and we obtain the existence of convex but not necessarily smooth or uniformly convex maximizers of a suitable relaxation of the above problem to a space of less regular functions. Unlike in the case of the affine area, it is not clear to us whether maximizers of the relaxed problem are unique. See Section \ref{Mproofs} for more about this relaxation. Therefore, we state the following questions.  

\begin{mainque}
Assume that $\Omega\subset \R^{2n}$ is a uniformly convex smooth bounded open set. Does the functional $\mathcal{F}$ have maximizers in the set of all uniformly convex functions $H\in C^{\infty}(\overline{\Omega})$ satisfying (\ref{firstjet2})? Are the maximizers of the relaxed problem unique?
\end{mainque}

We conclude this section by indicating some other future directions. First, the construction and possibly the results presented in this section should extend to more general symplectic manifolds $(M,\omega)$ equipped with an affine structure and a flat affine symplectic connection. The matrix $\nabla^2 H$ of a uniformly convex function $H$ defines a so called Hessian Riemannian metric on $M$, see \cite{sy97}. When $\dim M = 2$, the quantity $\mathcal{G}(H)$ has a simple geometric interpretation as the ratio between the Riemannian and the symplectic areas of $M$.

Second, let $(W, \omega)$ be a $2n$-dimensional real symplectic vector space equipped with an $n$-form $\sigma$. Let $\Lambda$ be the subset of the oriented Lagrangian Grassmannian of $W$ consisting of those Lagrangian subspaces $L \subset W$ for which $\sigma|_L$ is a positive volume form, see \cite{Sol14}. The tangent space $T_L \Lambda$ is canonically identified with the space of bilinear symmetric forms on $L$. Denote by $K_L \subset T_L \Lambda$ the cone of positive forms. Every $g \in K_L$ defines a scalar product, and hence a positive volume form $\nu_g$ on $L$. For $L \in \Lambda$ and $g \in K_L$, we set 
\[
Z_{\sigma}(g) := \left( \frac{\nu_g}{(\sigma|_L)} \right)^{\frac{2}{n}}.
\]
This is a Lorentz--Finsler metric on $\Lambda$. It would be interesting to study the existence of time-functions, geodesics, and the Lorentzian distance on $(\Lambda,K,Z_{\sigma})$. Furthermore, $Z_{\sigma}(g)$ induces a weak Lorentz--Finsler metric on the cone of optical Hamiltonian diffeomorphisms of a compact symplectic manifold with boundary equipped with a Lagrangian distribution, see \cite{bp94b}. It would be interesting to explore the maximizers of the optimal extension problem in this context.

\paragraph{Acknowledgments.}  \addcontentsline{toc}{subsection}{Acknowledgments}
We would like to thank Fedor Nazarov for sharing with us his proof of the $L^1$-Bernstein inequality for non-negative trigonometric polynomials which appears here as Theorem \ref{bernstein}. We are also grateful to Stefan Nemirovski, Miguel S\'anchez and Stefan Suhr for discussions on Lorentzian geometry, to Alessio Figalli for pointing us to the literature on the affine area functional, to Marco Mazzucchelli for helping us with the graphics and to Yaron Ostrover for discussions on the motto ``flexibility is expensive''.

\medskip

{\small \noindent A.\ A.\ is partially supported by the DFG under the Collaborative Research Center SFB/TRR 191 - 281071066 (Symplectic Structures in Geometry, Algebra and Dynamics). 

\noindent
G.\ B.\ is partially supported by the
DFG under the Germany's Excellence Strategy EXC2181/1 - 390900948 (the Heidelberg STRUCTURES Excellence Cluster), the Collaborative Research Center SFB/TRR 191 - 281071066
(Symplectic Structures in Geometry, Algebra and Dynamics), and the Research Training Group RTG 2229
- 281869850 (Asymptotic Invariants and Limits of Groups and Spaces). 

\noindent
L.\ P.\ is partially supported by a Mercator Fellowship within the Collaborative Research Center SFB/TRR 191 - 281071066 (Symplectic Structures in Geometry, Algebra and Dynamics).}

\numberwithin{equation}{section}

\section{Proof of Proposition \ref{lorfinSp}}
\label{lorfinSpsec}

The subset
\[
\mathrm{sp}^+(2n) = \{ X\in \mathrm{sp}(2n) \mid (u,v)\mapsto \omega_0(u,Xv) \mbox{ is positive definite}\} = \{J_0 S \mid \mathrm{Sym}^+(2n)\}
\]
of the Lie algebra $\mathrm{sp}(2n)$ of $\mathrm{Sp}(2n)$ is clearly an open convex cone, invariant under the adjoint action of $\mathrm{Sp}(2n)$, i.e., under conjugacy by elements of $\mathrm{Sp}(2n)$, and satisfies
\[
\overline{\mathrm{sp}^+(2n)} \cap \overline{-\mathrm{sp}^+(2n)} = \{0\}.
\]
Therefore, it generates a bi-invariant cone distribution in the tangent bundle of $\mathrm{Sp}(2n)$ which satisfies the requirements of (i) in Definition \ref{deflorfin} from the Introduction.

The function
\[
G(X) = G(J_0 S) = (\det X)^{\frac{1}{2n}} = (\det S)^{\frac{1}{2n}}
\]
is smooth on $\mathrm{sp}^+(2n)$, positively 1-homogeneous and extends continuously to the closure of $\mathrm{sp}^+(2n)$ by setting it to be zero on the boundary. Moreover, it is strongly concave in all directions other than the radial one, meaning that
\[
\mathrm{d}^2G(X)\cdot (Y,Y)<0 \qquad \forall X\in \mathrm{sp}^+(2n), \; Y\in \mathrm{sp}(2n) \setminus \R X.
\]
This follows immediately from the following well known concavity property of the $N$-th root of the determinant on the cone $\mathrm{Sym}^+(N)$ of positive symmetric endomorphisms of $\R^N$, of which for sake of completeness we give a proof.

\begin{lem}\label{l:hessdet}
	Let $f: \mathrm{Sym}^+(N) \rightarrow \R$ be the smooth function $f(S) := (\det S)^{\frac{1}{N}}$. Then, for all $S\in \mathrm{Sym}^+(N)$ and all $H,H_1,H_2\in \mathrm{Sym}^+(N)$
	\[
	\begin{split}
		\mathrm{d}f(S)\cdot H &= \frac{1}{N} (\det S)^{\frac{1}{N}} \tr (S^{-1} H), \\
		\mathrm{d}^2 f(S)\cdot (H_1,H_2) &= \frac{1}{N^2} (\det S)^{\frac{1}{N}} \Bigl(\tr (S^{-1} H_1)\tr (S^{-1} H_1) - N \,\tr \bigl( (S^{-1}H_1)(S^{-1}H_2) \bigr) \Bigr).
	\end{split}
	\]
	It follows that
	\[  
	\mathrm{d}^2 f(S)\cdot (H,H) \leq 0 \qquad \forall S\in \mathrm{Sym}^+(N), \; \forall H \in \mathrm{Sym}(N),
	\]
	where the equality holds if and only if $H\in \R S$.
\end{lem}

\begin{proof}
	The first two equalities in the statement are readily obtained using the fact that the differential of the determinant has the form
	\[
	\mathrm{d} (\det) (A)\cdot H = \det A \cdot \tr (A^{-1} H) \qquad \forall A\in \mathrm{GL}(N), \; \forall A\in \mathrm{Hom}(\R^N,\R^N).
	\]
	Plugging $H=H_1=H_2\in \mathrm{Sym}(N)$ in the formula for the second differential of $f$ at $S\in \mathrm{Sym}^+(N)$, we find
	\[
	\mathrm{d}^2 f(S)\cdot (H,H) = \frac{1}{N^2} (\det S)^{\frac{1}{N}} \Bigl( \bigl( \tr (S^{-1} H)\bigr)^2 - N \,\tr \bigl( (S^{-1}H)^2 \bigr) \Bigr).
	\]
	Since the symmetric endomorphism $S^{-1}$ is positive, it has a square root $S^{-\frac{1}{2}}\in \mathrm{Sym}^+(N)$ and by the conjugacy invariance of the trace we can rewrite the last identity as
	\[
	\mathrm{d}^2 f(S)\cdot(H,H) = \frac{1}{N^2} (\det S)^{\frac{1}{N}} \Bigl( \bigl( \tr A \bigr)^2 - N\, \tr \bigl( A^2 \bigr) \Bigr),
	\]
	where $A:= S^{-\frac{1}{2}} H S^{-\frac{1}{2}}$ belongs to $\mathrm{Sym}(N)$. Since $A$ is diagonalizable over $\R$, the Cauchy--Schwarz inequality implies that the above quantity is not larger than zero, and equal zero if and only if $A= \alpha I$ for some $\alpha\in \R$, i.e., if and only if $H= \alpha S$.
\end{proof}

Being invariant under the adjoint action of $\mathrm{Sp}(2n)$, $G$ extends to a bi-invariant function on the bi-invariant cone distribution in $T\mathrm{Sp}(2n)$ which is induced by $\mathrm{sp}^+(2n)$. 
Actually, since every $W\in \mathrm{Sp}(2n)$ has determinant 1, this extension is still the $2n$-th root of the determinant:
\[
G(Y)= (\det Y)^{\frac{1}{2n}} \qquad \forall\, Y=XW\in \mathrm{sp}^+(2n) W \subset T_W \mathrm{Sp}(2n), \; \forall\, W\in \mathrm{Sp}(2n).
\]
This extended function $G$ satisfies the requirements of (ii) in Definition \ref{deflorfin}. This concludes the proof of Proposition \ref{lorfinSp} from the Introduction.

\begin{rem}
{\rm The cone distribution $\{\mathrm{sp}^+(2n) W\}_{W\in \mathrm{Sp}(2n)}$ fits into the definition of a Lipschitz cone structure from \cite{fs12} and of the (more general) cone field from \cite{bs18}. Note that when $n>1$, the boundary of $\mathrm{sp}^+(2n)$ is not a smooth hypersurface, even after removing the origin, as singularities occur at each $X\in \partial \mathrm{sp}^+(2n)$ having zero as an eigenvalue with multiplicity  larger than one. For the same reason, this boundary is not strongly convex. Due to these facts, this cone structure satisfies neither the smoothness requirement of a weak cone structure nor the strong convexity requirement of a strong cone structure from \cite{js20}. The function $G$ satisfies the conditions of a Lorentz--Finsler metric on the cone structure $\{\mathrm{sp}^+(2n) W\}_{W\in \mathrm{Sp}(2n)}$, as defined in \cite{js20}, except for the fact that $G^2$ should be smooth up to the boundary of the cone minus the zero section. We refer to \cite{min16} and \cite{js20} for a discussion on the various definitions of a Lorentz--Finsler structure and their relationships.}
\end{rem}

\section{Proof of Proposition \ref{lorfinCont}}
\label{lorfinContsec} 

Let $\xi$ be a co-oriented contact structure on the closed $(2n-1)$-dimensional manifold $M$, where $n\geq 1$. In Section \ref{secD} of the Introduction, we have defined the cone $\mathrm{cont}^+(M,\xi)$ in the Lie algebra $\mathrm{cont}(M,\xi)$ of contact vector fields as the set of vector fields $X$ that are positively transverse to $\xi$, and the function
\[
V: \mathrm{cont}^+(M,\xi) \rightarrow \R
\]
as
\[
V(X) := \mathrm{vol}(M,\alpha)^{-\frac{1}{n}},
\]
where $\alpha$ is the unique contact form defining $\xi$ whose Reeb vector field $R_{\alpha}$ coincides with $X$. The aim of this section is to show that $\mathrm{cont}^+(M,\xi)$ and $V$ define a bi-invariant Lorentz--Finsler structure on $\mathrm{Cont}(M,\xi)$, hence proving Proposition \ref{lorfinCont} from the Introduction.

The convexity of $\mathrm{cont}^+(M,\xi)$ and its invariance under the adjoint action of $\mathrm{Cont}(M,\xi)$, i.e., under the push-forward of contact vector fields by contactomorphisms, are clear, and hence $\mathrm{cont}^+(M,\xi)$ defines a bi-invariant cone distribution in the tangent bundle of $\mathrm{Cont}(M,\xi)$. The intersection 
\[
\overline{\mathrm{cont}^+(M,\xi)} \cap \overline{- \mathrm{cont}^+(M,\xi)}
\]
consists of contact vector fields that are sections of the contact structure $\xi$, but the only contact vector field with this property is the zero vector field, see Remark \ref{no-tangency} in Appendix \ref{contham}.  This proves that the cone distribution generated by $\mathrm{cont}^+(M,\xi)$ satisfies the requirements of Definition \ref{deflorfin} (i) from the Introduction.

The invariance of $V$ under the adjoint action of $\mathrm{Cont}(M,\xi)$ on $\mathrm{cont}^+(M,\xi)$ is also clear: If $\phi\in \mathrm{Cont}(M,\xi)$ and $X=R_{\alpha} \in \mathrm{cont}^+(M,\xi)$, then
\[
\phi_* X = \phi_* R_{\alpha} = R_{\phi_* \alpha}
\]
and hence
\[
V(\phi_* X) = \mathrm{vol}(M,\phi_* \alpha)^{-\frac{1}{n}} = \mathrm{vol}(M,\alpha)^{-\frac{1}{n}} = V(X).
\]
In order to study the concavity of $V$, it is useful to fix a contact form $\alpha$ defining $\xi$ and use the identification
\begin{equation}
\label{the-identification}
\mathrm{cont}(M,\xi) \rightarrow C^{\infty}(M), \qquad X \mapsto \imath_X \alpha,
\end{equation}
whose inverse is denoted by
\begin{equation}
\label{the-identification2}
C^{\infty}(M) \rightarrow \mathrm{cont}(M,\xi), \qquad H \mapsto X_H.
\end{equation}
See Appendix \ref{contham} for the properties of this identification that we use here. When $H$ is positive, the Reeb vector field of the contact form $H^{-1} \alpha$ is precisely $X_H$ (see identity (\ref{Reeb-contact}) in Appendix \ref{contham}) and we obtain the formula
\begin{equation}
\label{Valpha}
V(X_H)  = \mathrm{vol}(M, H^{-1} \alpha)^{-\frac{1}{n}} = \left( \int_M H^{-n}  \alpha\wedge \mathrm{d} \alpha^{n-1} \right)^{-\frac{1}{n}}.
\end{equation}

The function $V$ is clearly positive and positively 1-homogeneous on $\mathrm{cont}^+(M,\xi)$.  Moreover, it is smooth, in the sense that it admits directional derivatives of every order, and a short computation leads to the following formulas.
\begin{lem}
Let $H\in C^{\infty}(M)$ be positive and set
\[
\mu:=H^{-n-2}\alpha\wedge\di\alpha^{n-1}.
\]
Then for every $K\in C^\infty(M)$, we have
\begin{equation}\label{e:dV}
	\di V(X_H)\cdot X_K =\frac{\mathrm{d}}{\mathrm{d}t}\Big|_{t=0}V(X_H+t  X_K)=c_0\int_MHK\mu\, ,
\end{equation}
where $c_0=\mathrm{vol}(M,H^{-1}\alpha)^{-\frac1n-1}$. Moreover, for every $K_1,K_2\in C^\infty(M)$ we have
\begin{equation}\label{e:ddV}
	\begin{aligned}
		\di^2V(X_H) & \cdot(X_{K_1},X_{K_2})=\frac{\di}{\di t}\Big|_{t=0}\di V(X_H+t X_{K_1})\cdot X_{K_2}\\
		&=c_1\left[\left(\int_MHK_1\mu\right)\left(\int_MHK_2\mu\right)-\left(\int_M H^2\mu\right)\left(\int_MK_1K_2\mu\right)\right],
	\end{aligned}
\end{equation}
where $c_1=(n+1)\mathrm{vol}(M,H^{-1}\alpha)^{-\frac1n-2}$. 
\end{lem}

Using these formulas we can check that $V$ is strongly concave in every direction other than the radial one. Indeed, if $X=X_H$ with $H\in C^{\infty}(M)$ positive and $Y=X_K\in \mathrm{cont}(M,\xi)$ with $K\in C^\infty(M)$, then we have
\[
\frac{\mathrm{d}^2}{\mathrm{d}t^2}\Big|_{t=0}V(X+t Y)=\di^2V(X)\cdot(Y,Y)=c_1\left[\left(\int_MHK\mu\right)^2-\left(\int_MH^2\mu\right)\left(\int_MK^2\right) \right]
\]
By the Cauchy--Schwarz inequality, the above quantity is negative when $K$ is not a multiple of $H$, i.e., when $Y$ is not a multiple of $X$ and we conclude that
\[
\frac{\mathrm{d}^2}{\mathrm{d}t^2}\Big|_{t=0}  V(X+t Y) < 0 \qquad \forall X\in \mathrm{cont}^+(M,\xi), \; \forall Y\in \mathrm{cont}(M,\xi) \setminus \R X,
\]
so $V$ has the claimed strong concavity property.

There remains to check that $V$ extends continuously to $\overline{\mathrm{cont}^+(M,\xi)}$ by setting it to be zero on the boundary, where the closure refers to any vector space topology on $\mathrm{cont}(M,\xi)$ which is finer than the $C^0$-topology on $C^{\infty}(M)$, after the identification (\ref{the-identification}). Indeed, the vector field $X$ belongs to the boundary of $\mathrm{cont}^+(M,\xi)$ if and only if $X=X_H$ for some non-negative contact Hamiltonian $H$ that vanishes somewhere. If $x_0$ is a point of $M$ at which $H$ vanishes, then $x_0$ is a minimum for $H$ and we have 
\[
H(x) \leq c\, \mathrm{dist}(x,x_0)^2 \qquad \forall x\in M,
\]
for a suitably large constant $c$, where $\mathrm{dist}$ is the distance function that is induced by an auxiliary Riemannian metric on $M$. Then
\[
H(x)^{-n} \geq c^{-n} \mathrm{dist}(x,x_0)^{-2n} \qquad \forall x\in M,
\]
and since $2n$ is larger than the dimension of $M$ we have
\[
\int_M H^{-n} \, \alpha\wedge \mathrm{d}\alpha^{n-1} = +\infty.
\]
If the sequence $(X_{H_j}) \subset  \mathrm{cont}^+(M,\xi)$ converges to $X=X_H$ then $H_j$ converges to $H$ pointwise (and even uniformly), so Fatou's lemma implies that
\[
\lim_{j\rightarrow \infty} \int_M H^{-n}_j \, \alpha\wedge \mathrm{d}\alpha^{n-1} = +\infty,
\]
and hence $V(X_{H_j})$ converges to zero, as we wished to show. This concludes the proof of Proposition \ref{lorfinCont} from the Introduction.

\begin{rem}
\label{nofinsler}
{\rm Identify $\mathrm{cont}(M,\xi)$ with $C^{\infty}(M)$ by means of a contact form $\alpha$ defining $\xi$.
As mentioned in the Introduction, any $C^0$-continuous function $F: C^{\infty}(M) \rightarrow \R$ which is invariant under the adjoint action of $\mathrm{Cont}(M,\xi)$ and vanishes at zero must vanish on all functions which are supported in a Darboux chart, i.e., in the image of a diffeomorphism $\varphi:B_r \rightarrow M$ which pulls $\alpha$ back to the standard 1-form
\[
\frac{1}{2} \sum_{j=1}^{n-1} ( x_j \, \mathrm{d} y_j - y_j \, \mathrm{d} x_j ) + \di z.
\]
Here, $B_r$ denotes the ball of radius $r$ centered at $0$ in $\R^{2n-1}$, whose points are denoted by $(x,y,z)= (x_1,y_1,\dots,x_{n-1},y_{n-1},z)$. Indeed, let $H\in C^{\infty}(M)$ be supported in $\varphi(B_r)$. Choose $s<r$ such that $\varphi(B_s)$ contains the support of $H$ and consider a group of contactomorphisms $\{\phi^t: M \rightarrow M\}_{t\in \R}$ which are supported in $\varphi(B_r)$ and satisfy  
\[
\phi^t\circ \varphi (x,y,z) = \varphi (e^t x, e^t y, e^{2t} z),
\]
for all $(x,y,z)\in B_s$ such that $(e^t x, e^t y, e^{2t} z)\in B_s$. Such a family of contactomorphisms can be defined by integrating the contact vector field $X_K$ corresponding to a contact Hamiltonian $K\in C^{\infty}(M)$ which is supported in $\varphi(B_r)$ and satisfies
\[
K\circ \varphi  (x_1,y_1,\dots,x_{n-1},y_{n-1},z) = 2z \qquad \mbox{on } B_s.
\]
By identity \eqref{adjoint-action} from Appendix \ref{contham}, we have $\phi^t_*(X_H) = X_{H_t}$, where $H_t\in C^{\infty}(M)$ is supported in $\varphi(B_r)$ and satisfies
\[
H_t\circ \varphi(x,y,z) = e^{2t} H \circ \varphi (e^{-t} x, e^{-t} y, e^{-2t} z)
\]
for all $(x,y,z)\in B_s$ such that $(e^{-t} x, e^{-t} y, e^{-2t} z)\in B_s$. The above identity and the fact that $H$ is supported in $\varphi(B_s)$ imply that $H_t$ converges to zero uniformly on $M$ for $t\rightarrow -\infty$. By invariance, we have $F(H)=F(H_t)$ for every $t\in \R$ and hence a limit for $t\rightarrow -\infty$ and the $C^0$-continuity of $F$ yield $F(H)=F(0)=0$.
}
\end{rem}

\section{Proof of Proposition \ref{rotnumb} and uniqueness}
\label{rotnumbsec}

Recall that $\T$ denotes the 1-torus $\R/\Z$, whose trivial contact structure $\xi_0:=\{0\}$ is co-oriented by the standard orientation of $\T$.  The 1-form $\mathrm{d}x$, $x$ being the standard coordinate on $\R$, is a defining contact form for $\xi_0$. In this case,
\[
\mathrm{Cont}(\T,\xi_0)= \mathrm{Cont}_0(\T,\xi_0) = \mathrm{Diff}_0(\T)
\]
is the group of orientation preserving diffeomorphisms of $\T$, which is connected. Let $X\in \mathrm{cont}^+(\T,\xi_0)$ and let $H:= \imath_X \mathrm{d}x$ be the corresponding contact Hamiltonian, i.e.,
\[
X = H \frac{\partial}{\partial x}.
\]
We denote by $\phi^t: \R \rightarrow \R$ the canonical lift of the flow of $X$, that is, the solution of
\[
\frac{\mathrm{d}}{\mathrm{d}t} \phi^t(x) = H(\phi^t(x)), \qquad \phi^0(x) = x, \qquad \forall x\in \R,
\]
where we are seeing $H$ as a 1-periodic function on $\R$. By dividing by the term on the right-hand side and integrating on $[0,t]$ we obtain
\[
\int_0^t \frac{1}{H(\phi^s(x))} \frac{\mathrm{d}}{\mathrm{d}s} \phi^s(x)\, \mathrm{d}s = t,
\]
which thanks to the change of variable $y=\phi^t(x)$ can be rewritten as
\[
\int_{x}^{\phi^t(x)} \frac{\mathrm{d}y}{H(y)} = t.
\] 
Taking the inverse of both sides we find
\[
\frac{\phi^t(x) - x}{t} = \left( \frac{1}{\phi^t(x) - x} \int_{x}^{\phi^t(x)} \frac{\mathrm{d}y}{H(y)} \right)^{-1}.
\]
By the periodicity of $H$, the term in brackets on the right-hand side converges to
\[
\int_0^1 \frac{\mathrm{d}y}{H(y)},
\]
for $t\rightarrow +\infty$, and we conclude that
\[ 
\lim_{t\rightarrow +\infty} \frac{\phi^t(x) - x}{t} = \left( \int_{\T} \frac{\mathrm{d}y}{H(y)} \right)^{-1} \qquad \forall x\in \R.
\]
The term on the left-hand side is the translation number $\rho(\phi^1)$ of the diffeomorphism $\phi^1$, while the term on the right-hand side is precisely $V(X)$, thanks to (\ref{Valpha}). This proves Proposition \ref{rotnumb} from the Introduction.

We now discuss the uniqueness property of the bi-invariant Lorentz--Finsler metric $V$ on $(\mathrm{Cont}(\T,\xi_0),\mathrm{cont}^+(\T,\xi_0))$ that is mentioned in Remark \ref{unique-2} in the Introduction.

\begin{prop}
\label{uniqueness2} 
Let $W: \mathrm{cont}^+(\T,\xi_0)\rightarrow \R$ be a positive function that is positively 1-homogeneous and invariant under the adjoint action of $\mathrm{Cont}(\T,\xi_0) = \mathrm{Diff}_0(\T)$. Then there exists $c>0$ such that $W=cV$. In particular, $V$ is, up to rescaling, the unique bi-invariant Lorentz--Finsler metric on $(\mathrm{Cont}(\T,\xi_0),\mathrm{cont}^+(\T,\xi_0))$. 
\end{prop}

\begin{proof}
Since both $V$ and $W$ are positive, positively 1-homogeneous and invariant under the adjoint action of $\mathrm{Diff}_0(\T)$, it is enough to show that $\mathrm{Diff}_0(\T)$ acts transitively on rays in $\mathrm{cont}^+(\T,\xi_0)$. Let 
\[
X(x) = H(x) \frac{\partial}{\partial x}
\]
be any element of $\mathrm{cont}^+(\T,\xi_0)$, where $H$ is positive smooth function on $\T$. Set
\[
h:= \int_0^1 H(y)\, \mathrm{d}y.
\]
Then the diffeomorphism $\phi\in \mathrm{Diff}_0(\T)$ which is defined by
\[
\phi(x) := \frac{1}{h} \int_0^x H(y)\, \mathrm{d}y
\]
satisfies
\[
\phi_* \left( h \frac{\partial}{\partial x} \right) = X.\qedhere
\]
\end{proof}

\section{Proof of Proposition \ref{sottogruppi}}
\label{sottogruppisec}

Recall that
\[
\lambda_0= \frac{1}{2} \sum_{j=1}^n ( x_j \, \mathrm{d}y_j -  y_j \, \mathrm{d}x_j), \qquad \mbox{i.e.,} \qquad \lambda_0(z)[v] = \frac{1}{2}\, \omega_0 (z,v)\quad \forall\,z,v\in\R^{2n},
\]
denotes the standard radial primitive of the symplectic form $\omega_0$ on $\R^{2n}$. The restriction of $\lambda_0$ to the sphere
\[
S^{2n-1} := \{ z\in \R^{2n} \mid |z|=1 \}
\]
is a contact form defining the standard (co-oriented) contact structure $\xi_{\mathrm{st}}$. We denote this restriction by $\alpha_0$. Being invariant under the antipodal symmetry $z\mapsto -z$, the contact form $\alpha_0$ induces a contact form, which we still denote by $\alpha_0$, on the real projective space $\R \mathrm{P}^{2n-1}$. The corresponding contact structure is the standard contact structure of $\R \mathrm{P}^{2n-1}$ and is also denoted by $\xi_{\mathrm{st}}$.

It is well known and easy to check that a diffeomorphism $\phi$ of $\R^{2n}\setminus \{0\}$ preserves $\lambda_0$, i.e., $\phi^* \lambda_0 = \lambda_0$, if and only if it is a positively 1-homogeneous symplectomorphism: $\phi( r z) = r \phi(z)$ for every $z\in \R^{2n}\setminus \{0\}$ and $r>0$ and $\phi^* \omega_0 = \omega_0$. If $\phi$ is such a positively 1-homogeneous symplectomorphism, then the diffeomorphism
\[
\psi: S^{2n-1} \rightarrow S^{2n-1}, \qquad \psi(z) := \frac{\phi(z)}{|\phi(z)|},
\]
is readily seen to satisfy
\[
\psi^* \alpha_0 = \frac{1}{|\phi|^2} \alpha_0,
\]
and hence is a contactomorphism of $(S^{2n-1},\xi_{\mathrm{st}})$. Actually, all contactomorphisms of $(S^{2n-1},\xi_{\mathrm{st}})$ arise from this construction.

The elements of the symplectic group $\mathrm{Sp}(2n)$ are 1-homogeneous symplectomorphisms of $\R^{2n}$ and we obtain a map
\[
i: \mathrm{Sp}(2n) \rightarrow \mathrm{Cont}(S^{2n-1},\xi_{\mathrm{st}}), \qquad i(W) : z \mapsto \frac{Wz}{|Wz|},
\]
which is readily seen to be an injective homomorphism. Let $X$ be a tangent vector to $\mathrm{Sp}(2n)$ at the identity, i.e., an element of the Lie algebra $\mathrm{sp}(2n)$. Note that
\[
\mathrm{d}i(\mathrm{id}) \cdot X = \frac{\mathrm{d}}{\mathrm{d}t} \Big|_{t=0} i\bigl(e^{tX} \bigr) = P X,
\]
where
\[
(PX)(z) = Xz - (Xz\cdot z) z \qquad \forall z\in S^{2n-1}
\]
denotes the tangent vector field on $S^{2n-1}$ that is obtained by projecting the restriction of $X$ orthogonally onto the tangent spaces of the sphere. Since the radial direction is in the kernel of $\lambda_0$, for every $z\in S^{2n-1}$ we find
\[
H(z):=\bigl( \imath_{\mathrm{d}i(\mathrm{id})\cdot X} \alpha_0 \bigr) (z) = \bigl( \imath_{Xz} \lambda_0 \bigr) (z) = \frac{1}{2} \omega_0(z,X z),
\]
and hence $\mathrm{d}i(\mathrm{id})\cdot X$ is in $\mathrm{cont}^+(S^{2n-1},\xi_{\mathrm{st}})$ if and only if $X$ is in $\mathrm{sp}^+(2n)$. This proves the identity
\[
\mathrm{d}i(\mathrm{id})^{-1} \bigl( \mathrm{cont}^+(S^{2n-1},\xi_{\mathrm{st}}) \bigr) = \mathrm{sp}^+(2n).
\]
If $X$ is in $\mathrm{sp}^+(2n)$, writing $X=J_0 S$ with $S\in \mathrm{Sym}^+(2n)$, we get
\[
H(z)= \frac{1}{2} \omega_0(z,X z) = \frac{1}{2} J_0 z \cdot J_0 S z = \frac{1}{2} S z\cdot z.
\]
Thanks to the positivity of $S$, the set $\{H=1\}$ is an ellipsoid, and the radial projection
\[
p: S^{2n-1} \rightarrow \{H=1\}
\]
is easily seen to satisfy
\[
p^*\bigl( \lambda_0|_{\{H=1\}} \bigr) = H^{-1} \alpha_0,
\]
thanks to the 2-homogeneity of $H$. Therefore, (\ref{Valpha}) gives us
\[
V(\mathrm{d} i(\mathrm{id})\cdot X)^{-n} = \mathrm{vol}  (S^{2n-1},H^{-1} \alpha_0) =  \mathrm{vol}( \{H=1\},  \lambda_0|_{\{H=1\}}) = \mathrm{vol}( \{H<1\}, \omega_0^n),
\]
where in the last step we have used Stokes theorem. Since the ellipsoid $\{H<1\}$ is the image of the unit ball $B^{2n}\subset \R^{2n}$ by the linear map $\sqrt{2} S^{-\frac{1}{2}}$, we have
\[
 \mathrm{vol} ( \{H< 1\}, \omega_0^n) = \det( \sqrt{2} S^{-\frac{1}{2}} \bigr)\, \mathrm{vol} (B^{2n}, \omega_0^n) = (2\pi)^n \bigl( \det S \bigr)^{-\frac{1}{2}},
 \]
 and hence
 \[
 V(\mathrm{d}i(\mathrm{id})\cdot X) = \frac{1}{2\pi} (\det S)^{\frac{1}{2n}} = \frac{1}{2\pi} G(X).
 \]
This establishes the part of Proposition \ref{sottogruppi} concerning the homomorphism $i$ into the contactomorphism group of the sphere.

The argument for the injective homomorphism
\[
j: \mathrm{PSp}(2n) \rightarrow \mathrm{Cont}(\R \mathrm{P}^{2n-1},\xi_{\mathrm{st}})
\]
is analogous. In order to determine the scaling factor in that case, note that if
\[
p: \mathrm{Sp}(2n) \rightarrow \mathrm{PSp}(2n), \qquad q: S^{2n-1} \rightarrow \R \mathrm{P}^{2n-1}
\]
are the quotient projections, then we have for every $W\in \mathrm{Sp}(2n)$ the commutative diagram
\[
\begin{CD}
S^{2n-1} & @>{i(W)}>> & S^{2n-1} \\
@V{q}VV & & @VV{q}V \\
 \R \mathrm{P}^{2n-1} & @>{j(p(W))}>>  &  \R \mathrm{P}^{2n-1}.
 \end{CD}
 \]
By identifying the Lie algebras of $\mathrm{Sp}(2n)$ and $\mathrm{PSp}(2n)$ by the differential of $p$ at the identity, we deduce that for every $X\in \mathrm{sp}(2n)$ the vector fields
$Z=\mathrm{d}j(\mathrm{id})\cdot X$ and $Y = \mathrm{d}i(\mathrm{id})\cdot X$ are related by the identity $Z=q_* Y$. Since $q$ is a double cover intertwining the standard contact forms $\alpha_0$ of $S^{2n-1}$ and $ \R \mathrm{P}^{2n-1}$, we find by (\ref{Valpha})
\[
V(Y) = \left( \int_{S^{2n-1}}\!\! ( \imath_Y \alpha_0)^{-n} \alpha_0 \wedge \mathrm{d}\alpha_0^{n-1} \right)^{-\frac{1}{n}} = \left( 2 \int_{\R \mathrm{P}^{2n-1}}\!\! ( \imath_Z \alpha_0)^{-n} \alpha_0 \wedge \mathrm{d}\alpha_0^{n-1} \right)^{-\frac{1}{n}} = 2^{-\frac{1}{n}} V(Z),
\]
and hence
\[
 V(\mathrm{d} j(I)\cdot X) =   2^{\frac{1}{n}} \, V(\mathrm{d}(\mathrm{id}) \cdot X) = \frac{2^{\frac{1}{n}} }{2\pi} G(X).
\]
This proves the part of Proposition \ref{sottogruppi} concerning the homomorphism $j$. The part concerning the homomorphism $j_k$ follows from the case of the homomorphism 
\[
j: \mathrm{PSp}(2) \rightarrow \mathrm{Cont}(\R \mathrm{P}^1,\xi_{\mathrm{st}}) = \mathrm{Diff}_0(\T) 
\]
by the observation that small neighborhoods of the identity in $\mathrm{PSp}(2)$ and $\mathrm{PSp}_k(2)$ can be identified and the restriction of $j_k$ to such a neighborhood is the composition of $j$ with the following map from a neighborhood of the identity in $\mathrm{Diff}_0(\T)$ to $\mathrm{Diff}_0(\T)$:
\[
\phi \mapsto \psi, \qquad \mbox{where} \quad \psi(x) := \frac{1}{k} \tilde{\phi}(kx)  \mod 1,
\]
where $\tilde{\phi}: \R \rightarrow \R$ denotes the unique lift of $\phi$ with $\tilde{\phi}(0)$ close to 0. Indeed, this yields the implication 
\[
\mathrm{d}j(\mathrm{id}) \cdot X = H \, \frac{\partial}{\partial x} \qquad \Rightarrow \qquad \mathrm{d}j_k(\mathrm{id}) \cdot X = H_k \, \frac{\partial}{\partial x} \quad \mbox{with} \quad H_k(x) = \frac{1}{k} H(kx),
\]
and the identity 
\[
V(\mathrm{d}j_k(\mathrm{id}) \cdot X) = \frac{1}{\pi k} G(X)
\]
follows. This concludes the proof of Proposition \ref{sottogruppi} from the Introduction.

\begin{rem}
\label{poltrig}
{\rm The above argument also shows that the image of the linear map
\[
\mathrm{d}j_k(\mathrm{id}) : \mathrm{sp}(2) \rightarrow \mathrm{cont}(\T,\xi_{0}) = C^{\infty}(\T)
\]
is the space of functions of the form
\[
x\mapsto q\bigl(\cos (\pi k x), \sin (\pi k x) \bigr),
\]
where $q: \R^2 \rightarrow \R$ is a quadratic form. This is the space of trigonometric polynomials of the form
\[
x\mapsto a \cos (2\pi k x) + b \sin (2\pi k x) + c,
\]
for $a,b,c\in \R$. This observation will be useful later on.}
\end{rem}

\section{The Morse co-index theorem for timelike geodesics on $\mathrm{Sp}(2n)$}
\label{varsec}

In Appendix \ref{liegroups}, we compute the first and second variation of the length functional which is associated to any bi-invariant Lorentz--Finsler metric on a Lie group. To this purpose, we need the first and second variation of $G$ which we computed in Lemma \ref{l:hessdet}:
\begin{align}
\label{diffg}
\mathrm{d}G(X) \cdot Y &= \frac{1}{2n} (\det X)^{\frac{1}{2n}} \tr (X^{-1} Y),  \\
\label{diff2g}
\mathrm{d}^2G(X) \cdot (Y_1,Y_2) &= \frac{1}{4n^2} (\det X)^{\frac{1}{2n}} \left( \tr (X^{-1} Y_1)\, \tr (X^{-1} Y_2) - 2n \, \tr(X^{-1}Y_1X^{-1} Y_2) \right), 
\end{align}
for every $X\in \mathrm{sp}^+(2n)$ and every $Y,Y_1,Y_2\in \mathrm{sp}(2n)$. Therefore, Proposition \ref{lemvar1b} and \eqref{diffg} give us the following formula for the first variation of the length functional
\[
\mathrm{length}_G(W) = \int_0^1 G(W'(t))\, \mathrm{d}t
\]
 at a timelike curve $W:[0,1] \rightarrow \mathrm{Sp}(2n)$ in the direction of a curve $Y: [0,1] \rightarrow \mathrm{sp}(2n)$:
\[
\mathrm{d} \, \mathrm{length}_G(W) \cdot Y = \frac{1}{2n} \int_0^1 (\det X)^{\frac{1}{2n}}\, \tr  (X^{-1} Y' ) \, \mathrm{d}t,
\]
where the curve $X: [0,1] \rightarrow \mathrm{sp}^+(2n)$ is defined by $W'=XW$. By the same proposition, timelike curves $W$ such that $\mathrm{d} \, \mathrm{length}_G(W) \cdot Y$ vanishes for every curve $Y: [0,1] \rightarrow \mathrm{sp}(2n)$ with compact support in $(0,1)$ are precisely time reparametrizations of autonomous curves: $W(t) = e^{\tau(t) X} W_0$ for some $X\in \mathrm{sp}^+(2n)$ and $W_0\in \mathrm{Sp}(2n)$.

Therefore, the timelike geodesics on the Lorentz--Finsler manifold $(\mathrm{Sp}(2n), \mathrm{sp}^+(2n),G)$, i.e., timelike extremals of $\mathrm{length}_G$ that are parametrized with constant speed, are precisely the autonomous positive paths, i.e., the curves of the form
\[
W(t) = e^{tX} W_0,
\] 
where $X\in \mathrm{sp}^+(2n)$ and $W_0\in \mathrm{Sp}(2n)$. Again, this is a general fact and holds true for any bi-invariant Lorentz--Finsler metric on a Lie group (see Appendix \ref{liegroups}). In particular, a different bi-invariant Lorentz--Finsler metric on $(\mathrm{Sp}(2n), \mathrm{sp}^+ (2n))$ would give us the same geodesics, but would measure their length differently.

\begin{rem}
{\rm 
Lightlike geodesics on manifolds equipped with a Lorentz--Finsler structure are defined in \cite[Definition 2.9]{js20} as ``cone geodesics'', a notion that uses only the boundary of the causal cones. The well-posedness of the initial value problem for them uses the smoothness and strong convexity of these boundary cones. As already mentioned, the boundary of the cone $\mathrm{sp}^+(2n)$ is neither smooth nor strongly convex, even after removing the zero section. However, in our situation it seems natural to define lightlike geodesics just as curves of the form $W(t) = e^{tX} W_0$ where $W_0\in \mathrm{Sp}(2n)$ and $X$ is a non-zero element of the boundary of $\mathrm{sp}^+(2n)$, i.e., a non-zero endomorphism of $\R^{2n}$ such that the bilinear form $(u,v) \mapsto \omega_0(u,Xv)$ is symmetric, positive semidefinite but not positive definite. Equivalently, $X=J_0 S$ with $S\in \mathrm{Sym}(2n)\setminus \{0\}$ such that $S\geq 0$ and $\ker S\neq 0$. Lightlike geodesics starting at the identity are confined to the discriminant, i.e., to the singular hypersurface $\Sigma_1 \subset \mathrm{Sp}(2n)$ consisting of symplectic automorphisms having the eigenvalue 1.}
\end{rem} 

Proposition \ref{lemvar2b} and \eqref{diff2g} give us the following formula for the second variation of $\mathrm{length}_G$ at the timelike geodesic $W: [0,1] \rightarrow \mathrm{Sp}(2n)$, $W(t) = e^{tX} W_0$, with $X\in \mathrm{sp}^+(2n)$ and $W_0\in \mathrm{Sp}(2n)$: for every pair of curves $Y_1,Y_2:[0,1] \rightarrow \mathrm{sp}(2n)$ vanishing at $t=0$ and $t=1$ we have
\[
\begin{split}
\mathrm{d}^2 \, & \mathrm{length}_G(W) \cdot (Y_1,Y_2)  =  - \frac{(\det X)^{\frac{1}{2n}}}{2n} \int_0^1  \Bigl(  \tr (X^{-1} Y_1'  X^{-1} Y_2' ) + \\ & -  \frac{1}{2n}\bigl( \tr  (X^{-1} Y_1') \bigr)\bigl( \tr  (X^{-1} Y_2') \bigr)  +   \tr  (Y_1 X^{-1} Y_2')- \tr (X^{-1} Y_1 Y_2'   ) \Bigr)  \, \mathrm{d}t,
\end{split}
\]
This symmetric bilinear form has infinite dimensional kernel, because the length functional $\mathrm{length}_G$ is invariant under reparametrizations. As explained in Appendix \ref{liegroups}, this invariance can be killed by restricting the second variation to curves $Y:[0,1] \rightarrow \mathrm{sp}(2n)$ taking values in the kernel of $\mathrm{d} G(X)|_{\mathrm{sp}(2n)}$, i.e., in the hyperplane
\[
\mathrm{sp}_X(2n) := \{ Y\in \mathrm{sp}(2n) \mid \mathrm{tr} (X^{-1} Y ) = 0\}.
\]
Restricting the second variation to the Sobolev space $H^1_0((0,1),\mathrm{sp}_X(2n))$ of absolutely continuous curves in $\mathrm{sp}_X(2n)$ which vanish at the end-points of the interval $[0,1]$ and have square integrable derivative, yields a continuous symmetric bilinear form
which has a finite dimensional kernel and a finite co-index. 

The {\it Morse co-index} of the timelike geodesic segment $W: [0,1] \rightarrow \mathrm{Sp}(2n)$ is the co-index of $\mathrm{d}^2\, \mathrm{length}_G(W)$, or equivalently of its restriction to $H^1_0((0,1),\mathrm{sp}_X(2n))$:
\[
\begin{split}
\mbox{co-ind}(W) := & \mbox{co-ind} (\mathrm{d}^2\, \mathrm{length}(W)) =  \max \{ \dim V \mid  V \mbox{ linear subspace of } \\ & H^1_0((0,1),\mathrm{sp}(2n)),
 \mathrm{d}^2 \, \mathrm{length}_G(W) \mbox{ is positive definite on } V\}.
\end{split}
\]
The kernel of the restriction of $\mathrm{d}^2 \, \mathrm{length}_G(W)$ to $H^1_0((0,1),\mathrm{sp}_X(2n))$ consists of the {\it Jacobi vector fields} along $W$, i.e., the solutions $Y: [0,1] \rightarrow \mathrm{sp}(2n)$ of the equation
\[
Y'' = [X,Y'],
\]
such that $Y(0)=Y(1)=0$. A number $t^*\in (0,1]$ for which there exist non-trivial Jacobi vector fields $Y$ such that $Y(0)=Y(t^*)=0$ is called {\it conjugate instant} and its {\it multiplicity} $m(t^*)$ is the dimension of the space of Jacobi vector fields with this property. The Morse co-index theorem holds:
\begin{equation}
\label{morse}
\mbox{co-ind}(W) = \sum_{t^*\in (0,1)} m(t^*).
\end{equation}
See Appendix \ref{liegroups} for the proof of these facts for general bi-invariant Lorentz--Finsler metrics on Lie groups.

\section{Jacobi fields along timelike geodesics in $\mathrm{Sp}(2n)$ and proof of Theorem \ref{Morse}}
\label{jacsec}

In this section, we wish to determine the Jacobi vector fields, the conjugate instants and the co-index of timelike geodesic segments in $\mathrm{Sp}(2n)$.  We start with the special case of the $2\pi$-periodic geodesic
\[
W:\R \rightarrow \mathrm{Sp}(2n), \qquad W(t) :=  e^{tJ},
\]
where $J$ is an $\omega_0$-compatible complex structure on $\R^{2n}$. The Jacobi vector fields along $W$ are the solutions $Y: \R \rightarrow \mathrm{sp}(2n)$ of the equation
\begin{equation}
\label{jacobiJ}
Y''=[J,Y'].
\end{equation}
We recall that the solutions of the commutator equation
\[
Z'=[X,Z]
\]
are given by
\[
Z(t) = e^{tX} Z_0 e^{-tX},
\]
where $Z_0=Z(0)$. We deduce that the solutions of \eqref{jacobiJ} vanishing at $t=0$ are of the form
\begin{equation}
\label{Yint}
Y(t) = \int_0^t e^{sJ} Y_0 e^{-sJ}\, \di s
\end{equation}
with $Y_0=Y'(0)\in \mathrm{sp}(2n)$. The vector space $\mathrm{sp}(2n)$ has the linear splitting
\[
\mathrm{sp}(2n) = \mathrm{sp}^c(2n)  \oplus \mathrm{sp}^a(2n), \quad X = X^c + X^a, \; X^c:= \frac{1}{2} (X-JXJ) , \;X^a:= \frac{1}{2} (X+JXJ),
\]
where the subspaces
\[
 \mathrm{sp}^c(2n) := \{X\in  \mathrm{sp}(2n) \mid XJ=JX\}, \qquad   \mathrm{sp}^a(2n) := \{X\in  \mathrm{sp}(2n) \mid XJ=-JX\},
 \]
have dimension
\[
\dim  \mathrm{sp}^c(2n) = n^2, \qquad \dim  \mathrm{sp}^a(2n) = n^2+n.
\]
Every $X$ in $\mathrm{sp}^c(2n)$ commutes with $e^{sJ}$, whereas if $X$ is in $\mathrm{sp}^a(2n)$ we have $e^{sJ} X=X e^{-sJ}$, for every $s\in \R$. Therefore, (\ref{Yint}) can be rewritten as
\[
Y(t) = \int_0^t e^{sJ} (Y_0^c + Y_0^a) e^{-sJ}\, \di s = \int_0^t Y_0^c\, \di s +  \int_0^t e^{2sJ} Y_0^a\, \di s = t Y_0^c - \frac{J}{2}  (e^{2tJ}-I) Y_0^a.
\]
Since $t Y_0^c$ is in $\mathrm{sp}^c(2n)$ and $\frac{J}{2}  (e^{2tJ}- \mathrm{id} ) Y_0^a$ is in $\mathrm{sp}^a(2n)$, we deduce that $Y(t)$ vanishes for a given $t\neq 0$ if and only if $Y_0^c=0$ and $(e^{2tJ}- \mathrm{id} ) Y_0^a=0$. This shows that the instant $t^*>0$ is conjugate to $t=0$ for the timelike geodesic $W(t)=e^{tJ}$ if and only if $t^*\in \N \pi$, and in this case its multiplicity is $n^2+n$. Together with the Morse co-index formula (\ref{morse}), we deduce that for every $T>0$ the geodesic segment $W|_{[0,T]}$ has co-index 
\[
\mbox{co-ind}(W|_{[0,T]}) = \left( \Bigl\lceil \frac{T}{\pi} \Bigr\rceil - 1 \right) (n^2+n).
\]
This proves Theorem \ref{Morse} from the Introduction.

\bigskip

We conclude this section by determining the conjugate instants of an arbitrary timelike geodesic 
\[
W(t) := e^{tX} W_0, \qquad X\in \mathrm{sp}^+(2n), \; W_0\in \mathrm{Sp}(2n),
\] 
on $\mathrm{Sp}(2n)$. By the representation of Appendix \ref{linear-app}, we have, collecting together the identical eigenvalues,
\begin{equation}
\label{lambdaform}
\R^{2n} = \bigoplus_{\lambda \in \Lambda} V_{\lambda}, \qquad
X = \bigoplus_{\lambda \in \Lambda} \lambda J_{\lambda}.
\end{equation}
Here, $\Lambda$ is a finite set of positive numbers, the splitting of $\R^{2n}$ is symplectic, and $J_{\lambda}$ is a $\omega_0$-compatible complex structure on the symplectic subspace $V_{\lambda}$, which has dimension $2n_{\lambda}$.

If $Z \in\mathrm{sp}(2n)$, the Jacobi field $Y$ along $W$ with $Y(0)=0$ and $Y'(0)=Z$ has the form
\[
Y(t)=\int_0^t e^{s X} Z e^{-s X}\, \mathrm{d} s.
\]
The multiplicity $m(t)$ of the possible conjugacy instant $t>0$ is the dimension of the kernel of the linear mapping
\[
\mathrm{sp}(2n) \rightarrow \mathrm{sp}(2n), \qquad Z \mapsto Y(t) = \int_0^t e^{s X} Z e^{-s X}\, \mathrm{d} s.
\]
In order to determine $m(t)$, we decompose
\[
\mathrm{sp}(2n)=\bigoplus_{\substack{ (\lambda_1,\lambda_2)\in \Lambda^2 \\ \lambda_1\leq \lambda_2}} \mathrm{sp}_{\lambda_1\lambda_2},
\]
where the elements of $\mathrm{sp}_{\lambda \lambda}$ are the maps $A_{\lambda}:V_\lambda\to V_\lambda$ in $\mathrm{sp}(2n)$ and the elements of $\mathrm{sp}_{\lambda_1\lambda_2}$ for $\lambda_1 < \lambda_2$ are the maps $A_{\lambda_1\lambda_2}\oplus A_{\lambda_2\lambda_1}:V_{\lambda_1}\oplus V_{\lambda_2}\to V_{\lambda_2}\oplus V_{\lambda_1}$ in $\mathrm{sp}(2n)$. The map
\[
\mathrm{sp}_{\lambda_1 \lambda_2}\to \mathrm{Hom}(V_{\lambda_1},V_{\lambda_2}),\qquad A_{\lambda_1\lambda_2}\oplus A_{\lambda_2\lambda_1}\mapsto A_{\lambda_1\lambda_2}
\]
is an isomorphism. We can now decompose $\mathrm{sp}_{\lambda \lambda} =\mathrm{sp}^c_{\lambda \lambda} \oplus\mathrm{sp}^a_{\lambda \lambda}$, where $\mathrm{sp}^c_{\lambda \lambda}$ are the elements that commute with $J_\lambda$ and $\mathrm{sp}^a_{\lambda \lambda}$ are the elements that anti-commute with $J_\lambda$. If $2n_\lambda=\dim V_\lambda$, then $\dim\mathrm{sp}^c_{\lambda\lambda} =n_\lambda^2$ and $\dim\mathrm{sp}^a_{\lambda \lambda} =n_\lambda^2+n_\lambda$.

Similarly, for $\lambda_1 < \lambda_2$ we decompose $\mathrm{sp}_{\lambda_1\lambda_2}=\mathrm{sp}^c_{\lambda_1\lambda_2}\oplus\mathrm{sp}^a_{\lambda_1\lambda_2}$. Here $\mathrm{sp}^c_{\lambda_1\lambda_2}$ is made of the elements such that $A_{\lambda_1\lambda_2}$ (equivalently $A_{\lambda_2\lambda_1}$) intertwines $J_{\lambda_1}$ and $J_{\lambda_2}$, i.e.,
\[
J_{\lambda_2}A_{\lambda_1\lambda_2}=A_{\lambda_1\lambda_2}J_{\lambda_1},
\]
whereas $\mathrm{sp}^a_{\lambda_1\lambda_2}$ is made of the elements such that $A_{\lambda_1\lambda_2}$ (equivalently $A_{\lambda_2\lambda_1}$) anti-intertwines $J_{\lambda_1}$ and $J_{\lambda_2}$, i.e.,
\[
J_{\lambda_2}A_{\lambda_1\lambda_2}=-A_{\lambda_1\lambda_2}J_{\lambda_1}.
\]
There holds $\dim\mathrm{sp}_{\lambda_1\lambda_2}^c=\dim\mathrm{sp}_{\lambda_1\lambda_2}^a=2n_{\lambda_1}n_{\lambda_2}$.

If we decompose $Z\in\mathrm{sp}(2n)$ as
\[
Z=\bigoplus_{(\lambda_1,\lambda_2) \in\Lambda^2}(Z_{\lambda_1\lambda_2}^c\oplus Z_{\lambda_1\lambda_2}^a),
\]
then the corresponding path $Y(t)$ can be written with respect to the splitting as
\[
Y(t)=\bigoplus_{(\lambda_1,\lambda_2) \in\Lambda^2}\int_0^t e^{ s \lambda_2 J_{\lambda_2}}Z_{\lambda_1\lambda_2}^ce^{-s \lambda_1J_{\lambda_1}}\, \mathrm{d} s + \bigoplus_{(\lambda_1,\lambda_2)\in\Lambda^2} \int_0^t e^{s \lambda_2 J_{\lambda_2}}Z_{\lambda_1\lambda_2}^ae^{- s\lambda_1J_{\lambda_1}}\, \mathrm{d} s.
\]
Let us compute the integrals in the first direct sum. We distinguish two cases. For $\lambda_1=\lambda=\lambda_2$ we get
\[
\int_0^t e^{s \lambda J_{\lambda}}Z_{\lambda\lambda}^ce^{- s\lambda J_{\lambda}}\, \mathrm{d}s=\int_0^t Z_{\lambda\lambda}^c \, \mathrm{d}s=t Z_{\lambda\lambda}^c,
\]
which vanishes only at $t=0$ when $Z^c_{\lambda\lambda}\neq0$. On the other hand, for $\lambda_1\neq\lambda_2$ we get
\[
\int_0^t e^{s \lambda_2 J_{\lambda_2}}Z_{\lambda_1\lambda_2}^ce^{-s \lambda_1J_{\lambda_1}}\, \mathrm{d} s=\int_0^t e^{s(\lambda_2-\lambda_1)J_{\lambda_2}}Z_{\lambda_1\lambda_2}^c \, \mathrm{d} s=\frac{1}{\lambda_1-\lambda_2}J_{\lambda_2}(e^{t(\lambda_2-\lambda_1)J_{\lambda_2}}-\mathrm{Id})Z_{\lambda_1\lambda_2}^c,
\]
which still intertwines $J_{\lambda_1}$ and $J_{\lambda_2}$ and vanishes exactly for $t\in \frac{2\pi}{|\lambda_2-\lambda_1|}\Z$ when $Z_{\lambda_1\lambda_2}^c\neq0$. Therefore, each pair $\lambda_1 < \lambda_2$ such that $t$ is an integer multiple of $2\pi/(\lambda_2 - \lambda_1)$ gives a contribution of $2n_{\lambda_1} n_{\lambda_2}$ to the multiplicity of $t$.

For the integrals in the second direct sum, we get without distinguishing cases
\[
\int_0^t e^{s\lambda_2 J_{\lambda_2}}Z_{\lambda_1\lambda_2}^ae^{-s\lambda_1J_{\lambda_1}}\, \mathrm{d} s=\int_0^t e^{s(\lambda_1+\lambda_2)J_{\lambda_2}}Z_{\lambda_1\lambda_2}^a\, \mathrm{d} s=-\frac{1}{\lambda_1+\lambda_2}J_{\lambda_2}(e^{t(\lambda_1+\lambda_2)J_{\lambda_2}}-\mathrm{id})Z_{\lambda_1\lambda_2}^a,
\]
which still anti-intertwines $J_{\lambda_1}$ and $J_{\lambda_2}$, and vanishes exactly for $t\in\frac{2\pi}{\lambda_1+\lambda_2}\Z$ when $Z_{\lambda_1\lambda_2}^a\neq0$. Therefore, for $t=\pi/\lambda$ we have a contribution of $n_{\lambda^2} + n_{\lambda}$ to the multiplicity of $t$, whereas every pair $\lambda_1 < \lambda_2$ such that $t$ is an integer multiple of $2\pi/(\lambda_1 + \lambda_2)$ gives a contribution of $2n_{\lambda_1} n_{\lambda_2}$ to the multiplicity of $t$. 

These considerations imply the following result, which generalizes Theorem \ref{Morse} from the Introduction.

\begin{thm}
\label{Morsegen}
Let $X\in \mathrm{sp}^+(2n)$ be of the form (\ref{lambdaform}). Then the multiplicity of $t^*>0$ as a conjugate instant to $t=0$ for the timelike geodesic $W(t) = e^{tX}$ is given by the formula
\[
m(t^*) = \sum_{\lambda\in \Lambda} (n_{\lambda}^2+n_{\lambda}) \mathbbm{1}_{\frac{\pi}{\lambda} \N} (t^*) + \sum_{\substack{(\lambda_1,\lambda_2) \in \Lambda^2 \\ \lambda_1 < \lambda_2}} 2 n_{\lambda_1} n_{\lambda_2} \mathbbm{1}_{\frac{2\pi}{\lambda_1+\lambda_2} \N \cup \frac{2\pi}{\lambda_2 - \lambda_1} \N} (t^*).
\]
In particular, setting $\lambda_1 := \max \Lambda$, the first conjugate instant is $t_1=\frac{\pi}{\lambda_1}$ with multiplicity $n^2_{\lambda_1}+n_{\lambda_1}$.
\end{thm}

\begin{rem}
The value $t_1$ is the first positive exit time of $W(t)=e^{tX}$, $t\geq0$ from the positively elliptic region $\mathrm{Sp}_{\mathrm{ell}}^+(2n)$. Therefore, all timelike geodesics originating form $\mathrm{id}$ are local maximizers of the Lorentz--Finsler length as long as they remain in $\mathrm{Sp}_{\mathrm{ell}}^+(2n)$. 
\end{rem}

\section{The second variation of the Lorentz--Finsler length on $\mathrm{Cont}(M,\xi)$} 
\label{secondvarsec}

As explained in Remark \ref{extensiontoDiff}, Propositions \ref{lemvar1b} and \ref{lemvar2b} from Appendix \ref{liegroups} hold also for the length functional which is induced by a bi-invariant Lorentz--Finsler metric on a group of diffeomorphisms. By applying Proposition \ref{lemvar1b} to $(\mathrm{Cont}(M,\xi), \mathrm{cont}^+(M,\xi),V)$, we obtain that the extremal curves of $\mathrm{length}_V$ are the autonomous positive paths of contactomorphisms, i.e.\ the Reeb flows which are induced by contact forms defining $\xi$.

Let $\alpha$ be a contact form defining $\xi$, denote by $\phi^t$ the flow of the corresponding Reeb vector field $R_{\alpha}$ and set
\[
\mu:= \alpha \wedge \mathrm{d} \alpha^{n-1}.
\]
Let $Y=\{Y_t\}_{t\in [0,T]}$ be a time-dependent contact vector field vanishing for $t=0$ and $t=T$ and let $K=\imath_Y \alpha$ be the corresponding time-dependent contact Hamiltonian. Using the identification (\ref{the-identification})-(\ref{the-identification2}) between $\mathrm{cont}(M,\xi)$ and $C^{\infty}(M)$ which is induced by the contact form $\alpha$, we then have $R_{\alpha}=X_1$ and $Y=X_K$. As discussed in Appendix \ref{contham}, the Lie bracket $[R_{\alpha},Y]$ is the contact vector field corresponding to the contact Hamiltonian
\begin{equation}\label{e:K1}
\{K,1\}=\di K(R_\alpha)\, ,
\end{equation}
see (\ref{e:bracketHK}). Since $\phi^t$ preserves the volume form $\mu$, we have
\begin{equation}\label{e:averagezero}
\int_M\{K,1\}\mu=0\,.
\end{equation}
Denote by $\{\theta_{s,t}\}_{(s,t)\in \R \times [0,T]}$ the smooth family of contactomorphisms of $(M,\xi)$ which is defined by
\[
\frac{\partial}{\partial s} \theta_{s,t} = Y_t (\theta_{s,t}), \qquad \theta_{0,t} = \mathrm{id}, \qquad \forall  (s,t)\in \R \times [0,T].
\]
Then the smooth family 
\[
\psi_{s,t} = \theta_{s,t} \circ \phi^t
\]
satisfies
\[
\psi_{0,t} = \phi^t, \qquad \psi_{s,0} = \mathrm{id}, \qquad \psi_{s,T} = \phi^T, \qquad  \forall  (s,t)\in \R \times [0,T],
\]
so for $|s|$ small the path of contactomorphisms $\{\psi_{s,t}\}_{t\in [0,T]}$ is positive and joins $\mathrm{id}$ with $\phi^T$.
By Proposition \ref{lemvar2b} we have
\begin{equation}\label{ildiff2}
\begin{split}
\frac{\mathrm{d}^2}{\mathrm{d} s^2} \Big|_{s=0} & \mathrm{length}_V \bigl( \{\psi_{s,t}\}_{t\in [0,T]} ) =
\mathrm{d}^2 \mathrm{length}_V \bigl( \{\phi^t\}_{t\in [0,T]} \bigr) (Y,Y) = \\ &=c_2\int_0^T\left[\frac{1}{\mathrm{vol}(M,\alpha)} \Bigl( \int_M \partial_t K\, \mu \Bigr)^2 - \int_M (\partial_t K)^2  \, \mu -\int_M   \mathrm{d} K(R_{\alpha}) \partial_t K\, \mu \right]\,  \mathrm{d} t\,,
\end{split}
\end{equation}
where $c_2:=(n+1)\mathrm{vol}(M,\alpha)^{-\frac1n-1}$ and we used \eqref{e:K1}, \eqref{e:averagezero}, and the formula for $\di^2V$ given in \eqref{e:ddV}.
 
By the Cauchy--Schwarz inequality, the first two terms in the square bracket define a quadratic form which is negative semidefinite. However, the presence of the third integral introduces an infinite dimensional subspace of contact vector fields on which the second variation of $\mathrm{length}_V$ is positive definite. In other words the second variation of the Lorentz--Finsler length $\mathrm{length}_V$ on the contactomorphism group has always infinite Morse index and infinite Morse co-index, as we stated Proposition \ref{infmorind} from the Introduction, which we now prove.

\begin{proof}[Proof of Proposition \ref{infmorind}]
Up to multiplying $\alpha$ by a positive number, we may rescale time and assume that $T=\pi$. Let $U$ be an open subset of $M$ which is diffeomorphic to the cube $(0,\epsilon)^{2n-1}$ and such that, using the coordinate system
\[
(r,z) = (r,x_1,y_1, \dots, x_{n-1},y_{n-1}), 
\]
which is induced by the identification $U\cong (0,\epsilon)^{2n-1}$, we have
\[
R_{\alpha}|_U = \frac{\partial}{\partial r} \qquad \mbox{and} \qquad \mu|_U = \mathrm{d} r \wedge  \mathrm{d}z =  \mathrm{d} r \wedge  \mathrm{d} x_1 \wedge \mathrm{d} y_1 \wedge \dots \wedge  \mathrm{d} x_{n-1} \wedge \mathrm{d} y_{n-1}.
\]
Let us fix a function $k$ supported in $(0,\epsilon)^{2n-2}$ and such that
\[
\int_{(0,\epsilon)^{2n-2}}k(z)\di z=0,\qquad \int_{(0,\epsilon)^{2n-2}}k(z)^2\di z=1\,.
\]
We shall compute the second variation of $\mathrm{length}_V$ along those time-dependent contact vector fields $Y$ that are induced by a contact Hamiltonian $K$ with support in $[0,\pi] \times U$ and of the form
\[
K(t,r,z) = k(z) \bigl( a(r) \sin t + b(r)\sin (2t)\bigr),
\]
where the functions $a$, $b$ are supported in $(0,\epsilon)$ and have vanishing integral. Note that $K(0,\cdot)=K(\pi,\cdot)=0$, and hence $Y$ is an admissible variation. By plugging a function $K$ of this form into \eqref{ildiff2} we obtain the following expression
\[
\begin{split}
\mathrm{d}^2 \mathrm{length}_V &  \bigl(\{\phi^t\}_{t\in[0,\pi]} \bigr) (Y,Y) = -c_2\int_0^{\pi}  \left[ \int_{(0,\epsilon)^{2n-1}} \!\! k(z)^2 \bigl( a(r) \cos t + 2 b(r) \cos (2t) \bigr)^2 \mathrm{d}r \wedge  \mathrm{d}z\, + \right.\\ &+  \left.\int_{(0,\epsilon)^{2n-1}} \!\! k(z)^2 \bigl( a'(r) \sin t + b'(r) \sin (2t) \bigr) (a(r) \cos t + 2 b(r) \cos (2t)\bigr)  \,  \mathrm{d}r \wedge \mathrm{d} z  \right]  \mathrm{d} t,
\end{split}
\]
where we have used the fact that $\partial_t K$ has vanishing integral, since $k,a,b$ have vanishing integral.
By switching the integrals and using the identities
\[
\begin{split}
\int_0^{\pi} \cos^2 t \, \di t = \int_0^{\pi} \cos^2 (2t) \, \di t = \frac{\pi}{2}, \quad \int_0^{\pi} \sin t \cos (2t) \, \di t = - \frac{2}{3},
\quad\int_0^{\pi} \sin (2t) \cos t \, \di t = \frac{4}{3},\\ \int_0^{\pi} \cos t \cos (2t)\, \di t = \int_0^{\pi} \sin t \cos t\, \di t = \int_0^{\pi} \sin(2t) \cos (2t)\, \di t = 0,
\end{split}
\]  
the above expression simplifies to
\[
\mathrm{d}^2 \mathrm{length}_V  \bigl( \{\phi^t\}_{t\in[0,\pi]} \bigr) (Y,Y) = - c_2  \int_0^{\epsilon} \Bigl( \frac{\pi}{2} a(r)^2 + 2\pi b(r)^2 - \frac{4}{3} a'(r) b(r) + \frac{4}{3} a(r) b'(r) \Bigr)\, \mathrm{d} r\,.
\]
Finally, an integration by parts gives
\[
\mathrm{d}^2 \mathrm{length}_V  \bigl( \{\phi^t\}_{t\in[0,\pi]} \bigr) (Y,Y) =  - c_2   \int_0^{\epsilon} \Bigl( \frac{\pi}{2} a(r)^2 + 2\pi b(r)^2 - \frac{8}{3} a'(r) b(r)  \Bigr)\, \mathrm{d} r\,.
\]
By choosing $b=0$, we find that the second variation is negative definite on the infinite dimensional space $W^-$ of time-dependent contact vector fields which are induced by Hamiltonians which are supported in $[0,\pi] \times U$ and there have the form
\[
K(t,r,z) = k(z) a(r) \sin t,
\]
where $a$ is any smooth function with compact support in $(0,\epsilon)$ and vanishing integral.

By choosing
\[
b(r) = \frac{2}{3\pi} a'(r),
\]
which has vanishing integral because $a$ is compactly supported, we obtain
\[
\mathrm{d}^2 \mathrm{length}_V  \bigl( \{\phi^t\}_{t\in[0,\pi]} \bigr) (Y,Y) = \frac{8}{9\pi} c_2  \int_0^{\epsilon} \Bigl( a'(r)^2 - \frac{9}{16} \pi^2 a(r)^2  \Bigr)\, \mathrm{d} r .
\]
This is a quadratic form of finite Morse index on the space of compactly supported smooth functions on $(0,\epsilon)$ with vanishing integral, and hence we can find an infinite dimensional vector space $A$ of functions $a$ as above on which this quadratic form is positive definite. The second  variation of $\mathrm{length}_V$ is then positive definite on the infinite dimensional space $W^+$ of time-dependent contact vector fields which are induced by Hamiltonians which are supported in $[0,\pi] \times U$ and there have the form
\[
K(t,r,z) = k(z) \left( a(r) \sin t + \frac{2}{3\pi} a'(r) \sin (2t) \right),
\]
where $a$ belongs to $A$.
\end{proof}

Now we study the conjugate instants for the geodesic $\phi^t$ in $\mathrm{Cont}(M,\xi)$, where as before $\phi^t$ is the Reeb flow of some contact form $\alpha$ defining $\phi$. As explained in Appendix \ref{liegroups}, the equation for Jacobi vector fields $Y$ along $\phi^t$ is
\[
\partial_{tt} Y = [R_{\alpha},\partial_t Y].
\]
Denoting as before by $K=\alpha(Y)$ the contact Hamiltonian associated to $Y$, and using the bracket induced on $C^\infty(M)$ the above equations reads
\[
\partial_{tt} K = \{1, \partial_t K\}
\]
and, plugging in the formula \eqref{e:bracketHK} for the bracket, we arrive at
\begin{equation}
\label{jacobicont}
\partial_{tt} K + \bigl(\di\,\partial_t K\bigr)(R_{\alpha}) = 0.
\end{equation}
Recalling that the general solution of the first order linear PDE
\[
\partial_{t}F+ \mathrm{d} F(R_{\alpha}) = 0
\]
is given by
\[
F(t,x) = f(\phi^{-t}(x)),
\]
where $f$ is any function on $M$, we obtain that the general solution $K$ of (\ref{jacobicont}) vanishing for $t=0$ is of the form
\[
K(t,x) = \int_0^t f(\phi^{-s}(x))\, \mathrm{d}s,
\]
for some arbitrary function $f$. Therefore, the positive number $t^*$ is a conjugate instant for the geodesic $\phi^t$ if and only if there exists a non-vanishing smooth function $f$ on $M$ such that
\[
\int_0^{t^*} f(\phi^{-s}(x))\, \mathrm{d}s= 0 \qquad \forall x\in M.
\]
By the change of variable $t=t^*-s$ and the fact that $\phi^t$ is a flow, this is equivalent to the condition
\[
\int_0^{t^*}  f(\phi^t(\phi^{-t^*}(x)))\, \mathrm{d}t  = 0 \qquad \forall x\in M,
\]
and hence to the condition
\begin{equation}
\label{conjinstcont}
\int_0^{t^*} f(\phi^t(x))\, \mathrm{d}t  = 0 \qquad \forall x\in M.
\end{equation}
We now determine all solutions of the above equation in the two simple cases that appear in the Introduction as Examples \ref{exintro1} and \ref{exintro2}.

\begin{ex}
{\rm Let $c>0$ and consider the geodesic $\phi^t$ in $\mathrm{Cont}_0(\T,\xi_0)= \mathrm{Diff}_0(\T)$ which is generated by the vector field
\begin{equation}
\label{vectfield}
X(x) = c \frac{\partial}{\partial x},
\end{equation}
where $c>0$. Then $V(X)=c$ and
 $\phi^t(x) = x + ct$, so equation (\ref{conjinstcont}) reads
\[
\int_0^{t^*} f(x+ct) \, \mathrm{d}t = 0 \qquad \forall x\in \T.
\]
By the change of variable $x+ct=s$, the above condition is seen to be equivalent to
\begin{equation}
\label{ppcc}
\int_x^{x+ct^*} f(s)\, \mathrm{d}s = 0 \qquad \forall x\in \T.
\end{equation}
If 
\begin{equation}
\label{rational}
t^* = \frac{p}{q} \cdot \frac{1}{c} =  \frac{p}{q} \cdot \frac{1}{V(X)}
\end{equation}
for some pair of natural numbers $p,q$, then any smooth function on $\R$ which is $\frac{1}{q}$-periodic and has vanishing integral on its period interval satisfies (\ref{ppcc}). This shows that all the positive numbers $t^*$ of the form (\ref{rational}) are conjugate instants and have infinite multiplicity. In order to show that these are the only conjugate instants, we need to show that (\ref{ppcc}) has no non-trivial smooth solution $f$ if $ct^*$ is an irrational number. By writing $f$ in Fourier series as
\[
f(x) = \sum_{k\in \Z} \hat{f}_k\, e^{2\pi i k x},
\]
we rewrite (\ref{ppcc}) as
\[
 ct^*\hat{f}_0+ \sum_{k\in \Z\setminus \{0\}} \frac{\hat{f}_k}{2\pi i k} \bigl(e^{2\pi i k c t^*} - 1\bigr) e^{2\pi i k x} = 0   \qquad \forall x\in \T.
\]
Since a 1-periodic function is identically zero if and only if all its Fourier coefficients vanish, the above condition is equivalent to
\[
\hat{f}_0=0 \quad \mbox{and} \quad \hat{f}_k \bigl( e^{2\pi i k c t^*} - 1\bigr) = 0 \qquad \forall k\in \Z\setminus \{0\}.
\]
The fact that $ct^*$ is irrational implies that all the Fourier coefficients of $f$ vanish, and hence (\ref{ppcc}) has only the trivial solution.

An arbitrary vector field
\[
X(x) = H(x)\, \frac{\partial}{\partial x}, \qquad x\in \T,
\]
with $H>0$ is conjugate to the vector field (\ref{vectfield}) with $c=V(X)$ (see the proof of Proposition \ref{uniqueness2}). Therefore, what we proved for conjugate instants in the special case (\ref{vectfield}) holds in general, as stated in Example \ref{exintro1} from the Introduction.
}
\end{ex}

\begin{rem}
{\rm The above example extends immediately to all geodesics in $\mathrm{Cont}(M,\xi)$ that are given by a Zoll contact form $\alpha$: If $T>0$ is the minimal period of all the orbits of the flow $\phi^t$ of $R_{\alpha}$, then the positive number $t^*$ is a conjugate instant for the geodesic $\phi^t$ if and only if it is of the form $T$ times a rational number. All conjugate instants have infinite multiplicity.}
\end{rem}

\begin{ex}
{\rm 
We now consider the setting of Example \ref{exintro2} of the Introduction: $\xi$ is the contact structure which is induced by the contact form
\[
\alpha(x,y,z) = \cos (2\pi z) \, \mathrm{d}x + \sin (2\pi z) \, \mathrm{d}y
\]
on $\T^3=\R^3/\Z^3$, and $\phi^t$ is the flow of $R_{\alpha}$, which is readily seen to have the form
\[
\phi^t(x,y,z) = \bigl( x + t  \cos (2\pi z) , y + t \sin (2\pi z), z).
\]
Then the positive number $t^*$ is a conjugate instant for the geodesic $\phi^t$ if and only if there is a non identically vanishing smooth function $f$ on $\T^3$ such that
\[
\int_0^{t^*} f(x + t  \cos (2\pi z) , y + t \sin (2\pi z), z)\, \di t = 0 \qquad \forall (x,y,z)\in \T^3.
\]
Using the Fourier representation
\[
f(x,y,z) = \sum_{(h,k)\in \Z^2} \hat{f}_{h,k} (z) e^{2\pi i (hx+ky)},
\]
for suitable smooth functions $\hat{f}_{h,k}: \T \rightarrow \C$, the above condition can be rewritten as
\begin{equation}
\label{conda}
\sum_{(h,k)\in \Z^2} \hat{f}_{h,k} (z) e^{2\pi i (hx+ky)} \int_0^{t^*} e^{2\pi i t ( h \cos (2\pi z) + k \sin (2\pi z))} \, \di t = 0 \qquad \forall (x,y,z)\in \T^3.
\end{equation}
Consider the analytic function $\varphi: \R \rightarrow \C$ 
\[
\varphi(s) = \frac{e^{2\pi i s}-1}{2\pi is} = \sum_{n=0}^{\infty} \frac{(2 \pi i s)^n}{(n+1)!}.
\]
Then 
\[
\int_0^{t^*} e^{2\pi i t ( h \cos (2\pi z) + k \sin (2\pi z))} \, \di t  = t^* \varphi \bigl( t^*  ( h \cos (2\pi z) + k \sin (2\pi z))\bigr),
\]
and hence (\ref{conda}) is equivalent to
\[
 \hat{f}_{h,k} (z) \, \varphi \bigl( t^*  ( h \cos (2\pi z) + k \sin (2\pi z))\bigr) = 0 \qquad \forall (h,k)\in \Z^2, \; \forall z\in \T.
 \]
 For every $(h,k)\in \Z^2$, the analytic function $z\mapsto \varphi \bigl( t^*  ( h \cos (2\pi z) + k \sin (2\pi z))\bigr)$ is not identically zero and hence has at most finitely many zeroes in $\T$. Therefore, the above condition implies that each function $\hat{f}_{h,k}$ is identically zero and hence $f=0$. This shows that the geodesic $\phi^t$ has no conjugate instants.
}
\end{ex}

\section{Proof of Theorem \ref{systolic}}
\label{systolicsec} 

In this section, we show how Theorem \ref{systolic} can be deduced from the fact that the Morse co-index of every timelike geodesic in $\mathrm{Cont}(M,\xi)$ is positive (and actually infinite), thanks to Proposition \ref{infmorind} from the Introduction.

Let $\phi^t_{\alpha_0}$ be the Reeb flow of the Zoll contact form $\alpha_0$ on $(M,\xi)$ and denote by $T_0$ the minimal period of its orbits. Then the positive path $\{\phi^t_{\alpha_0}\}_{t\in [0,T_0]}$ is a geodesic arc in $\mathrm{Cont}(M,\xi)$ and by Proposition \ref{infmorind} we can find a time-dependent contact vector field $Y=\{Y_t\}_{t\in [0,T]}$ vanishing for $t=0$ and $t=T$ and such that
\[
\mathrm{d}^2 \mathrm{length}_V (\{\phi^t_{\alpha_0}\}_{t\in [0,T_0]} \cdot (Y,Y) > 0.
\]
By (\ref{ildiff2}), the second differential of the length functional is continuous on the space of contact vector fields corresponding to contact Hamiltonians in the Hilbert space
\[
H^1_0((0,T),L^2(M,\mu)).
\]
Since smooth functions which are compactly supported in $(0,T)\times M$ are dense in this space, we can assume that the above time-dependent vector field $Y$ vanishes for $t$ in a neighborhood of $0$ and $T$.

As in the previous section, we denote by $\{\theta_{s,t}\}_{(s,t)\in \R \times [0,T_0]}$ the smooth family of contactomorphisms of $(M,\xi)$ which is defined by
\[
\frac{\partial}{\partial s} \theta_{s,t} = Y_t (\theta_{s,t}), \qquad \theta_{0,t} = \mathrm{id}, \qquad \forall  (s,t)\in \R \times [0,T_0].
\]
Then the smooth family 
\[
\psi_{s,t} = \theta_{s,t} \circ \phi^t_{\alpha_0}
\]
satisfies
\[
\psi_{0,t} = \phi^t_{\alpha_0}, \qquad \psi_{s,0} = \mathrm{id}, \qquad \psi_{s,T_0} = \phi^{T_0}_{\alpha_0} = \mathrm{id}, \qquad  \forall  (s,t)\in \R \times [0,T_0],
\]
so there exists $\epsilon>0$ such that if $|s|<\epsilon$ then the path of contactomorphisms $\{\psi_{s,t}\}_{t\in [0,T_0]}$ is positive and joins the identity with itself. 

It $t$ is close enough to $0$ or $T_0$, then $Y_t$ vanishes and hence $\psi_{s,t} = \phi^t_{\alpha_0}$. From the fact that $\phi^t_{\alpha_0}$ is far away from the identity for $t\in (0,T_0)$ far away from $0$ and $T_0$ and from the fact that $\psi_{0,t} = \phi^t_{\alpha_0}$ we deduce that, up to reducing $\epsilon$, we may assume that the diffeomorphisms $\psi_{s,t}$ do not have any fixed point when $(s,t)\in (-\epsilon,\epsilon) \times (0,T_0)$. 

The smooth function 
\[
f: (-\epsilon,\epsilon) \times [0,T_0] \rightarrow \R, \qquad f(s,r) := \mathrm{length}_V ( \{ \psi_{s,t} \}_{t\in [0,r]} ),
\]
satisfies
\[
\begin{split}
f(0,r) =  \mathrm{length}_V & ( \{ \phi^t_{\alpha_0} \}_{t\in [0,r]} ) , \qquad \frac{\partial f}{\partial s} (0,r) = 0, \qquad \forall r\in [0,T_0], \\ \frac{\partial^2 f}{\partial s^2} (0,T_0) &= \mathrm{d}^2 \mathrm{length}_V (\{\phi^t_{\alpha_0}\}_{t\in [0,T_0]} \cdot (Y,Y) > 0.
\end{split}
\]
Therefore, the function $s\mapsto f(s,T_0)$ has a non-degenerate local minimum at $s=0$. Together with the fact that $\frac{\partial f}{\partial r}(0,T_0)$ is positive, from the implicit function theorem we deduce that, up to a further reduction of $\epsilon$, there is a smooth function
\[
\tau: (-\epsilon,\epsilon) \rightarrow (0,T_0]
\]
such that $\tau(0)=T_0$, $\tau(s)< T_0$ for every $s\in (-\epsilon,\epsilon) \setminus \{0\}$ and
\[
f(s,\tau(s)) = f(0,T_0) = \mathrm{length}_V ( \{ \phi^t_{\alpha_0} \}_{t\in [0,T_0} ).
\]
We now consider the smooth family $\{\phi_s\}_{s\in (-\epsilon,\epsilon)}$ of positive paths 
\[
\phi_s: [0,T_0] \rightarrow \mathrm{Cont}(M,\xi) 
\]
which is  given by
\[
\phi_s(t) := \psi_{s, \frac{\tau(s)}{T_0} t}, \qquad \forall (s,t) \in (-\epsilon,\epsilon)\times [0,T_0].
\]
This family has the required properties: $\phi_0 = \phi_{\alpha_0}$, $\phi_s(0)=\mathrm{id}$,
\[
\mathrm{length}_V(\phi_s) = \mathrm{length}_V ( \{ \phi^t_{\alpha_0} \}_{t\in [0,T_0} ),
\]
for every $s\in (-\epsilon,\epsilon)$, and $\phi_s(t)$ has no fixed points if $s\neq 0$ and $t\in (0,T_0]$. This concludes the proof of Theorem \ref{systolic}.

\section{Krein theory, Maslov quasimorphism and proof of Theorem \ref{time}}
\label{timesec}

In this section, we recall some basic facts about Krein theory and about the homogeneous Maslov quasimorphism which we will need in the following sections. 

Every endomorphism $A$ of $\R^{2n}$ can be seen as linearly acting on $\C^{2n}$ in the usual way:
\[
A(u+iv) := Au +i Av \qquad \forall u,v\in \R^{2n}.
\]
An endomorphism $A$ of $\C^{2n}$ arises in this way if and only if it is real, i.e., $A(\R^{2n})\subset\R^{2n}$.
Similarly, the symplectic form $\omega_0$ extends by sesquilinearity uniquely to a skew-Hermitian form on $\C^{2n}$ by setting
\[
\omega_0(u+iv,u'+iv') := \omega_0(u,u') + \omega_0(v,v') + i (\omega_0(v,u')-\omega_0(u,v')), \qquad \forall u,v,u',v'\in \R^{2n}.
\]
 The {\it Krein-form} $\kappa: \C^{2n} \times \C^{2n} \rightarrow \C$ is defined as
\[
\kappa(w,w'):= \omega_0(-i w,w') = \langle - i J_0 w,w' \rangle \qquad \forall w,w'\in \C^{2n},
\]
where $\langle\cdot ,\cdot \rangle$ is the standard Hermitian product on $\C^{2n}$. The form $\kappa$ is easily seen to be Hermitian with signature $(n,n)$. An automorphism of $\C^{2n}$ corresponds to an element of $\mathrm{Sp}(2n)$ if and only if it is real and $\kappa$-unitary. Similarly, an endomorphism of $\C^{2n}$ corresponds to an element of the Lie algebra $\mathrm{sp}(2n)$ if and only if it is real and $\kappa$-skew-Hermitian. 

If $\lambda\in \U:= \{z\in \C \mid |z|=1\}$ is an eigenvalue of $W\in \mathrm{Sp}(2n)$, then $\kappa$ is non-degenerate on the corresponding algebraic eigenspace
\[
E(\lambda):= \{w\in\C^{2n} \mid w\in \Ker (\lambda I - W)^k \mbox{ for some } k\geq 1 \},
\]
and the signature of $\kappa$ on $E(\lambda)$ is called Krein-signature of $\lambda$. The eigenvalue $\lambda\in \U$ is said to be Krein-positive (resp.\ Krein-negative, resp.\ Krein-definite) if $\kappa$ is positive definite (resp.\ negative definite, resp.\ definite) on $E(\lambda)$. If $\lambda\in \U$ has Krein-signature $(p,q)$, then the conjugate eigenvalue $\overline{\lambda}$ has Krein-signature $(q,p)$ (see \cite[Lemma I.2.9]{eke90}). In particular, the eigenvalues 1 and -1 have signature $(p,p)$ and cannot be Krein-definite. 

The following well known lemma explains the role of the Krein-signature in the behaviour of eigenvalues of paths in $\mathrm{Sp}(2n)$.

\begin{lem}
\label{fundlem}
Let $W:(a,b) \rightarrow \mathrm{Sp}(2n)$ be a differentiable curve such that
\begin{equation}
\label{eigen}
W u = e^{i\theta} u,
\end{equation}
for some differentiable function $\theta:(a,b) \rightarrow \R$ and some differentiable curve of eigenvectors $u:(a,b) \rightarrow \C^{2n}$. Then
\begin{equation}
\label{fund0}
\kappa(u,u) \, \theta' = \omega_0 ( u, W'W^{-1} u).
\end{equation}
\end{lem}

\begin{proof}
By differentiating (\ref{eigen}) we find
\[
W' u + W u' = i \theta' e^{i\theta}  u + e^{i\theta} u',
\]
and by taking the $\kappa$-product with $Wu$ we obtain
\[
\kappa(W'u,Wu) + \kappa(Wu',Wu) = i\theta' e^{i\theta} \kappa(u,Wu) + e^{i\theta} \kappa(u',Wu).
\] 
Using the fact that $W$ is $\kappa$-unitary and that $u$ is an eigenvector of $W$ with eigenvalue $e^{i\theta}$, we can rewrite the above expression as
\[
e^{-i\theta}\kappa(W' u,u) + \kappa(u',u) = i \theta' \kappa(u,u) + \kappa(u',u),
\]
and (\ref{fund0}) follows since $e^{-i\theta}\kappa(W' u,u)=\kappa(W' e^{-i\theta}u,u)=\kappa(W' W^{-1}u,u)$.
\end{proof}

Note that the right-hand side of (\ref{fund0}) is positive (resp.\ non-negative) when the curve $W$ is timelike (resp.\ causal). Recalling that simple eigenvalues and the corresponding eigenvectors of differentiable paths are differentiable, the above lemma implies that if $W$ is a timelike (resp.\ causal) curve, then the argument of any Krein-positive simple eigenvalue of $W(t)$ on $\U$ is a strictly increasing (resp.\ non-decreasing) function of $t$. The same is true for Krein-definite eigenvalues of higher multiplicity, see \cite[Proposition I.3.2 and Corollary I.3]{eke90}.

We now consider the function
\[
\upsilon: \mathrm{Sp}(2n) \rightarrow \U, \qquad \upsilon(W) := (-1)^{m} \prod_{\lambda\in \sigma(W) \cap \U \setminus \{\pm 1\}} \lambda^{p(\lambda)},
\]
where $2m$ denotes the total algebraic multiplicity of real negative eigenvalues of $W$ and $(p(\lambda),q(\lambda))$ is the Krein-signature of the eigenvalue $\lambda \in \U$. This function, which was introduced by Gel'fand and Lidski\v{i} in \cite{gl58}, is continuous, invariant under symplectic conjugacy, homogeneous, i.e., $\upsilon(W^k) = \upsilon(W)^k$ for every $k\in \Z$, and coincides with the complex determinant on the subgroup $\mathrm{O}(2n) \cap \mathrm{Sp}(2n)\cong \mathrm{U}(n)$. Moreover, it induces an isomorphism of fundamental groups
\[
\upsilon_* : \pi_1(\mathrm{Sp}(2n)) \rightarrow \pi_1(\U) = \Z.
\]
See also \cite{sz92} and \cite[Section 1.3.4]{abb01} for the proof of these properties. We now consider the universal cover
\[
\pi : \widetilde{\mathrm{Sp}}(2n) \rightarrow \mathrm{Sp}(2n)
\]
and define the {\it homogeneous Maslov quasi-morphism}
\[
\mu: \widetilde{\mathrm{Sp}}(2n) \rightarrow \R
\]
as lift of $\upsilon$, and more precisely as the unique continuous function satisfying
\[
\upsilon(\pi(w)) = e^{2\pi i \mu(w)} \quad \forall w\in \widetilde{\mathrm{Sp}}(2n), \qquad \mu(\mathrm{id}) = 0.
\]
As recalled in  Section  \ref{secI} of the Introduction, this function is the unique homogeneous real quasi-morphism on $\widetilde{\mathrm{Sp}}(2n)$ whose restriction to $\pi^{-1}(\mathrm{U}(n))$ agrees with the lift of the complex determinant. See \cite{bss10} for more on this and for the proof of uniqueness. Alternative definitions of $\mu$ not requiring Krein theory are possible, but one advantage of the above definition is that it immediately yields the following result.

\begin{prop}
The homogeneous Maslov quasi-morphism $\mu$ is non-decreasing along every causal curve in $\widetilde{\mathrm{Sp}}(2n)$.
\end{prop}

Indeed, by the form of $\upsilon$ the function $\mu$ can change only when some eigenvalue on the unit circle $\U$ moves. Since, as explained above, for a causal curve all Krein-positive eigenvalues on $\U$ cannot move clockwise, the function $\mu$ is non-decreasing.

Now we prove Theorem \ref{time} from the Introduction: $\widetilde{\mathrm{Sp}}(2n)$ admits a time function, i.e., a continuous real function which increases strictly along each causal curve.

The homogeneous Maslov quasi-morphism is surely not a time function because it is locally constant on the open set consisting of elements $w\in \widetilde{\mathrm{Sp}}(2n)$ whose projection $\pi(w)$ has no eigenvalues on $\U$.  However, it is strictly increasing along causal curves which are contained in the set $\pi^{-1}(\Gamma)$, where $\Gamma$ denotes the open subset of $\mathrm{Sp}(2n)$ consisting of matrices having $2n$ distinct eigenvalues on $\U$. Indeed, if $w:(a,b) \rightarrow \pi^{-1}(\Gamma)$ is a causal curve and $W:=\pi\circ w$, then we can find a basis $u_1(t),\dots,u_n(t),\overline{u}_1(t), \dots, \overline{u}_n(t)$ of $\C^{2n}$ which is $\kappa$-unitary, i.e.,
\[
\kappa(u_j,u_j) = 1 = - \kappa(\overline{u}_j,\overline{u}_j) \; \forall j, \quad \kappa(u_j,u_h) = \kappa(\overline{u}_j,\overline{u}_h) = 0 \; \forall j\neq h, \quad \kappa(u_j,\overline{u}_h) = 0 \; \forall j,h,
\] 
and satisfies 
\[
W(t) u_j(t) = e^{i \theta_j(t)} u_j(t), \qquad W(t) \overline{u}_j(t)  = e^{-i \theta_j(t)} \overline{u}_j(t), \qquad  \forall t\in (a,b),
\]
where the real functions $\theta_j$ satisfy
\begin{equation}
\label{laespr}
\mu(w(t)) = \frac{1}{2\pi} \sum_{j=0}^n \theta_j(t) \qquad \forall t\in (a,b).
\end{equation}
Here, the curves $u_j$ and the real functions $\theta_j$ are differentiable. By Lemma \ref{fundlem}, we have
\[
\theta_j'  = \omega_0 ( u_j ,  W'W^{-1}u_j)  \qquad \mbox{on } (a,b).
\]
The fact that $W$ is a causal curve tells us that the Hermitian form $u\mapsto \omega_0 ( u , W' (t)W^{-1}(t) u)$ is positive semi-definite but not zero for every $t\in (a,b)$. Therefore, for every $t\in (a,b)$ the numbers $\theta_j'(t)$ are non-negative and at least one of them is positive. We conclude that the derivative of the sum of the $\theta_j$'s is strictly positive on $(a,b)$ and by (\ref{laespr}) the function $\mu\circ w$ is strictly increasing on $(a,b)$, as we claimed.

Now we choose a countable subset $\{w_j \mid j\in \N\}$ of $\widetilde{\mathrm{Sp}}(2n)$ such that
\begin{equation}
\label{totale}
\widetilde{\mathrm{Sp}}(2n) = \bigcup_{j\in \N} \pi^{-1}(\Gamma)w_j,
\end{equation}
and for $\epsilon>0$ we consider the function
\[
f: \widetilde{\mathrm{Sp}}(2n) \rightarrow \R, \qquad f(w) := \mu(w) + \epsilon \sum_{j=1}^{\infty}
 2^{-j} \arctan \mu  (ww_j^{-1}).
\]
This function is continuous because the above series converges uniformly. Thanks to the bi-invariance of the cone distribution defined by $\mathrm{sp}^+(2n)$, each summand is non-decreasing along each causal curve $t\mapsto w(t)$ and by (\ref{totale}) and what we have seen above about the behaviour of $\mu$ on $\pi^{-1}(\Gamma)$, for every $t_0$ at least one of the summands is strictly increasing for $t$ close to $t_0$. This proves that $f$ is strictly increasing on each causal curve in $\widetilde{\mathrm{Sp}}(2n)$ and hence is a time function. Since
\[
\|f-\mu\|_{\infty} < \frac{\pi}{2} \epsilon,
\]
$f$ is at bounded distance from $\mu$ and hence is also a quasimorphism. Moreover, the distance between $f$ and $\mu$ can be chosen to be arbitrarily small, as stated in the Introduction. This concludes the proof of Theorem \ref{time}.

\medskip

The homogenous Maslov quasimorphism is invariant under conjugacy in $\widetilde{\mathrm{Sp}}(2n)$, but the time function we have constructed above loses this invariance property. Actually, no time function on $\widetilde{\mathrm{Sp}}(2n)$ can be conjugacy invariant. Indeed, the following stronger statement holds:

\begin{prop}
\label{nonconju}
There are no continuous real functions on $\widetilde{\mathrm{Sp}}(2n)$ which are strictly increasing along timelike curves and invariant under conjugacy.
\end{prop}

\begin{proof}
It is enough to consider the case $n=1$. Assume that $\tilde{g}: \widetilde{\mathrm{Sp}}(2) \rightarrow \R$ is continuous, strictly increasing on timelike curves, and conjugacy invariant. Let $\Omega$ be the open subset of $\mathrm{Sp}(2)$ consisting of automorphisms without negative real eigenvalues. Since $\Omega$ is simply connected, the function $\tilde{g}$ descends to a continuous function $g: \Omega\rightarrow \R$ which is still strictly increasing on timelike curves and conjugacy invariant. The discriminant $\Sigma_1$, i.e., the subset of $\mathrm{Sp}(2)$ consisting of all automorphisms having the eigenvalue 1, is contained in $\Omega$ and is a two-dimensional double cone with vertex at the identity, see the left-hand picture in Figure \ref{figura} in the Introduction. Elements of each of the two components of $\Sigma_1\setminus \{\mathrm{id}\}$ are pairwise conjugate and since the identity is in the closure of both components, the function $g$ must be constant on the whole $\Sigma_1$. But there are timelike curves in $\Omega$ that go from one component of $\Sigma_1\setminus \{\mathrm{id}\}$ to the other one. An example is given by the timelike geodesic
\[
W(t) = e^{tJ_0} A \qquad \mbox{with} \quad A:= \left( \begin{array}{cc} 1 & 1 \\ 0 & 1 \end{array} \right).
\]
This curve passes through $A\in \Sigma_1$ at $t=0$, takes values in the space of positively hyperbolic automorphisms for $t\in (0,t^*)$ with $t^*=\arctan \frac{4}{3}$, and
\[
W(t^*) = \left( \begin{array}{cc} \frac{3}{5} & -\frac{1}{5} \\ \frac{4}{5} & \frac{7}{5} \end{array} \right)
\]
belongs again to $\Sigma_1$. The function $g\circ W$ is strictly increasing and this contradicts the fact that it takes identical values at $t=0$ and $t=t^*$.
\end{proof}

\section{Causality, Lorentz distance and proof of Theorem \ref{Lordist} on $\widetilde{\mathrm{Sp}}(2)$}
\label{Sp(2)sec}

In this section, we study some properties of causal curves in the universal cover of $\mathrm{Sp}(2)=\mathrm{SL}(2,\R)$. As discussed in  Section \ref{secC} of the Introduction, this space can be identified with the universal cover of the three-dimensional anti-de Sitter space $\mathrm{AdS}_3$, which is well studied. Therefore, we will be rather sketchy on the facts which are well known, see e.g., \cite[p.\ 131-134]{he10}, and add more details about statements which we could not find in the literature.

First, one can check that the symplectic group $\mathrm{Sp}(2)$ can be parametrized by the following diffeomorphism 
\begin{equation}
\label{parametrization}
S^2_+ \times (\R/2\pi\Z) \rightarrow \mathrm{Sp}(2), \quad (\varphi,\theta,t) \mapsto \frac{1}{\cos \varphi} {\scriptsize \left( \begin{array}{cc} \cos t & - \sin t \\ \sin t & \cos t \end{array} \right) } + \tan \varphi {\scriptsize  \left( \begin{array}{cc} \cos \theta &  \sin \theta \\ \sin \theta & - \cos \theta \end{array} \right)} .
\end{equation}
Here, $S^2_+$ denotes the open upper half-sphere in $\R^3$ with spherical coordinates $(\varphi,\theta)$ given by the latitude $\varphi\in [0,\frac{\pi}{2})$, where $\varphi=0$ corresponds to the north pole, and the longitude $\theta\in \R/2\pi \Z$. Moreover, the pull-back of the Lorentzian metric of $\mathrm{Sp}(2)$ by this diffeomorphism is the Lorentzian metric
\begin{equation}
\label{metric}
\frac{1}{\cos^2 \varphi} (\di s^2 - \di t^2),
\end{equation}
where $\di s^2$ refers to the round metric on $S^2_+$. Therefore, $\widetilde{\mathrm{Sp}}(2)$ is diffeomorphic to $S^2_+\times\R$ and, as Lorentz manifold, $\widetilde{\mathrm{Sp}}(2)$ is conformally equivalent to a portion of the {\it Einstein space} $(S^2\times \R,\di s^2-\di t^2)$, in which timelike and lightlike curves are easy to visualize.

\begin{figure}[h]
\hspace{3.8cm}
\begin{footnotesize}
  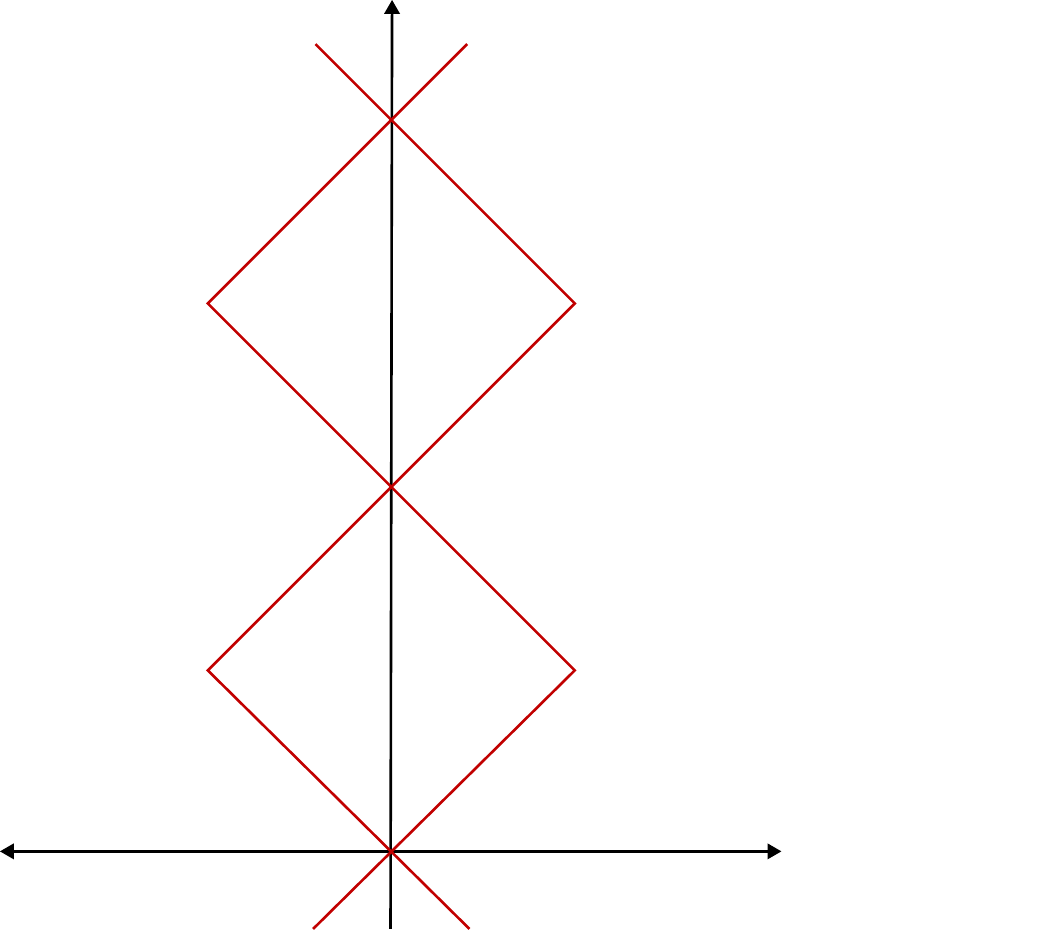
\end{footnotesize}
\caption{The identification of $\widetilde{\mathrm{Sp}}(2)$ with $S^2_+\times \R$.}
\label{AdS}
\end{figure}

Figure \ref{AdS} schematically represents $\widetilde{\mathrm{Sp}}(2)$ after the above identification with $S^2_+\times \R$. The universal cover of the group of rotations $\mathrm{U}(1)\subset \mathrm{Sp}(2)$ corresponds to $\varphi=0$, i.e., to the $t$-axis. The identity of $\widetilde{\mathrm{Sp}}(2)$ sits at $t=0$, the first lift of minus the identity, i.e., the element given by the homotopy class of the path $\{e^{\pi t J_0} \}_{t\in [0,1]}$, sits at $t=\pi$, and the lift of the identity in $\mathrm{Sp}(2)$ given by the homotopy class of the path $\{e^{2\pi t J_0} \}_{t\in [0,1]}$ sits at $t=2\pi$.

Lightlike curves (in red) have slope 1, while timelike curves (in green a timelike geodesic segment and in blue another timelike curve) have slope larger than 1. The set of lightlike curves emanating from the identity spans the cone $\widetilde{\Sigma}_1$ given by one component of the inverse image of the set of elements of $\mathrm{Sp}(2)$ having the eigenvalue 1 under the covering map $\pi: \widetilde{\mathrm{Sp}}(2) \rightarrow \mathrm{Sp}(2)$. The regions enclosed by the red diamonds correspond to lifts of elliptic elements in $\mathrm{Sp}(2)$, while the region outside of them to the lifts of hyperbolic elements. The positively elliptic region $\widetilde{\mathrm{Sp}}_{\mathrm{ell}}^+(2)$ is the shaded region enclosed by the lower diamond. It can be characterized as the set of points $w\in \widetilde{\mathrm{Sp}}(2)$ such that there is a timelike curve from $\mathrm{id}$ to $w$ and a timelike curve from $w$ to the first lift of $-\mathrm{id}$.

The submanifold $\widetilde{\mathrm{Sp}}_{\mathrm{ell}}^+(2)$ is {\it globally hyperbolic}, meaning that for every pair of points $w_0$, $w_1$ in it the set of causal curves from $w_0$ to $w_1$ spans a compact subset. Equivalently, $\widetilde{\mathrm{Sp}}_{\mathrm{ell}}^+(2)$ admits a {\it Cauchy hypersurface}, i.e., a hypersurface which is met exactly once by any inextensible causal curve in $\widetilde{\mathrm{Sp}}_{\mathrm{ell}}^+(2)$ (see \cite[Section 3.11]{ms08}). Here, a natural Cauchy hypersurface in $\widetilde{\mathrm{Sp}}_{\mathrm{ell}}^+(2)$ is given by $w$ such that $\pi(w)$ has spectrum $\{i,-i\}$, that is, by the set $S^2_+ \times \{\pi\}$ in the above identification with a portion of the Einstein space.

The Lorentzian distance $\mathrm{dist}_G$ on $\widetilde{\mathrm{Sp}}(2)$ is completely described by the following result.

\begin{prop}
\label{lordisSp2}
Let $w\in \widetilde{\mathrm{Sp}}(2)$ be such that $w\geq \mathrm{id}$. Then:
\begin{enumerate}[(i)]
\item If $w$ belongs to the closure of $\widetilde{\mathrm{Sp}}_{\mathrm{ell}}^+(2)$ then
\[
\mathrm{dist}_G(\mathrm{id},w) = \theta,
\]
where $\theta\in [0,\pi]$ is such that $e^{\pm i \theta}$ are the eigenvalues of $\pi(w)$. If moreover $w\in \widetilde{\mathrm{Sp}}_{\mathrm{ell}}^+(2)$, then $\mathrm{dist}_G(\mathrm{id},w)$ is achieved by the unique timelike geodesic segment from $\mathrm{id}$ to $w$.
\item If $w$ does not belong to the closure of $\widetilde{\mathrm{Sp}}_{\mathrm{ell}}^+(2)$, then  there are 
arbitrarily long timelike curves from $\mathrm{id}$ to $w$ and hence $\mathrm{dist}_G(\mathrm{id},w)=+\infty$.
\end{enumerate}
\end{prop}

Statement (i) is proven at the end of this section. As for (ii): if $w\geq \mathrm{id}$ is not in  the closure of $\widetilde{\mathrm{Sp}}_{\mathrm{ell}}^+(2)$, then an arbitrarily long timelike curve from $\mathrm{id}$ to $w$  can be obtained by first following a timelike curve which is very close to a lightlike curve and gets close to the boundary of $S^2_+\times \R$, then move along a short segment in the $t$-direction, thus acquiring large length because of the factor $\frac{1}{\cos \varphi}$ in the expression (\ref{metric}) for the Lorentzian metric, and then reach $w$ by following a path which is close to a lightlike one. See the blue curve in Figure \ref{AdS}. Actually, the last statement in (ii) is not specific of the Lorentzian distance induced by $G$ and holds in the following more general form.

\begin{prop}
\label{trivlordisSp}
Let $d: \widetilde{\mathrm{Sp}}(2) \times \widetilde{\mathrm{Sp}}(2) \rightarrow [0,+\infty]$ be a function such that:
\begin{enumerate}[(i)]
\item $d(w_0,w_1)>0$ if there is a timelike curve from $w_0$ to $w_1$;
\item $d(w_0,w_2) \geq d(w_0,w_1) + d(w_1,w_2)$ if $w_0\leq w_1 \leq w_2$;
\item $d$ is conjugacy invariant.
\end{enumerate}
If $w\geq \mathrm{id}$ and $w$ does not belong to the closure of $\widetilde{\mathrm{Sp}}_{\mathrm{ell}}^+ (2)$, then $d(\mathrm{id},w) = +\infty$.
\end{prop}

The proof of this proposition uses the following algebraic lemma.

\begin{lem}
\label{algebraic}
Let $W_0\in \mathrm{Sp}(2)$ be hyperbolic. Then there exists $H\in \mathrm{sp}(2)$ such that the curve
\[
W: \R \rightarrow \mathrm{Sp}(2), \qquad W(t):= e^{-tH} W_0 e^{tH}
\]
is timelike.
\end{lem}

\begin{proof}
Without loss of generality, we have
\[
W_0 = \left( \begin{array}{cc} \lambda & 0 \\ 0 & \lambda^{-1} \end{array} \right),
\]
for some real number $\lambda$ with $0< |\lambda| < 1$. If $W$ is as above, then
\[
W' W^{-1} = e^{-tH}  [W_0,H]W_0^{-1}  e^{tH},
\]
so the curve $W$ is timelike if and only if the element $[W_0,H]W_0^{-1}$ belongs to $\mathrm{sp}^+(2)$. By choosing
\[
H = \left( \begin{array}{cc} 0 & 1 \\ 1 &0 \end{array} \right),
\]
we compute
\[
[W_0,H]W_0^{-1} = \left( \begin{array}{cc} 0 & \lambda^{2}-1 \\ \lambda^{-2}-1 &0 \end{array} \right) = J_0 S,
\]
where the symmetric matrix
\[
S=  \left( \begin{array}{cc} \lambda^{-2}-1 & 0 \\ 0 &  1-\lambda^{2}  \end{array} \right) 
\]
is positive definite because $0<\lambda^2<1$.
\end{proof}

\begin{proof}[Proof of Proposition \ref{trivlordisSp}]
Since $w\geq \mathrm{id}$ is not in the closure of $\widetilde{\mathrm{Sp}}_{\mathrm{ell}}^+ (2)$, we can find a $w_0\in \widetilde{\mathrm{Sp}}(2)$ such that $\mathrm{id} \leq w_0 \leq w$ and $\pi(w_0)$ is hyperbolic. See Figure \ref{AdS}. By (ii), it is enough to prove that $d(\mathrm{id},w_0) = +\infty$. By Lemma \ref{algebraic}, there exists $w_1$ in the conjugacy class of $w_0$ such that there is a timelike curve from $w_0$ to $w_1$. By (iii) and (ii), we have
\[
d(\mathrm{id},w_0) = d(\mathrm{id},w_1) \geq d(\mathrm{id},w_0) + d(w_0,w_1).
\]
By (i), $d(w_0,w_1)>0$ and the above inequality forces $d(\mathrm{id},w_0)=+\infty$.
\end{proof}

We conclude this section by proving statement (i) of Proposition \ref{lordisSp2}. It is enough to consider the case of some $w\in \widetilde{\mathrm{Sp}}^+_{\mathrm{ell}}(2)$ with
\[
\sigma(\pi(w_1)) = \{ e^{\pm i \theta} \}, \qquad \theta\in (0,\pi),
\]
and prove that any timelike curve from $\mathrm{id}$ to $w$ has length at most $\theta$. Indeed, the timelike geodesic segment from $\mathrm{id}$ to $w$ has length $\theta$, being of the form $\{e^{t\theta J}\}_{t\in [0,1]}$ for some $\omega_0$-compatible complex structure $J$. Moreover, the bound on the length of causal curves which are not timelike and the case of a $w$ in the closure of $\widetilde{\mathrm{Sp}}^+_{\mathrm{ell}}(2)$ follow by an easy perturbation argument. 

Therefore, it is enough to consider a timelike curve
\[
W:[0,1] \rightarrow \mathrm{Sp}(2)
\]
such that $W(0)=\mathrm{id}$ and $W(t) \in \mathrm{Sp}^+_{\mathrm{ell}}(2)$ for every $t\in (0,1]$. Such a curve has the form
\begin{equation}
\label{form}
W(t) = A(t)^{-1} e^{\theta(t) J_0} A(t), \qquad \forall t\in (0,1],
\end{equation}
where $\theta(t)\in [0,\pi)$ and $A(t)\in \mathrm{Sp}(2)$, and we must prove the bound
\[
\mathrm{length}_G(W) \leq \theta(1).
\] 
The differentiability of $W$, together with the fact that the elements of $\mathrm{Sp}_{\mathrm{ell}}^+(2)$ have two distinct eigenvalues, implies that both $\theta$ and $A$ depend differentiably on $t$ in $(0,1]$, see \cite[Theorem II.5.4]{kat80}. The function $\theta$ is continuous at $t=0$ with $\theta(0)=0$, whereas $A$ needs not extend continuously at $t=0$.
For paths of the above form, we have the following simple lemma, in which we use the fact that elements of $\mathrm{sp}(2)$ have vanishing trace.

\begin{lem}
\label{det}
Let $W: (0,1] \rightarrow \mathrm{Sp}(2)$ be a path of the form $(\ref{form})$ for some differentiable functions $\theta:(0,1] \rightarrow \R$ and $A:(0,1] \rightarrow \mathrm{Sp}(2)$. If we denote the coefficients of the path of matrices $A'A^{-1} \in \mathrm{sp}(2)$ by
\[
A' A^{-1} = \left( \begin{array}{cc} a & b \\ c & - a \end{array} \right),
\]
we have the identity
\[
\det W' = {\theta'}^2  - \bigl( 4a^2 + (b+c)^2 \bigr) \sin^2 \theta \qquad \mbox{on } (0,1].
\]
\end{lem}  

\begin{proof}
The derivative of $W$ is
\[
W' = A^{-1} \bigl( [ e^{\theta J_0}, A' A^{-1} ] + \theta' J_0 e^{\theta J_0} \bigr) A.
\]
From the identities
\[
[ e^{\theta J_0}, A' A^{-1} ] = \sin \theta \left( \begin{array}{cc} - b -c & 2a \\ 2a & b+c  \end{array} \right), \qquad J_0 e^{\theta J_0} = \left( \begin{array}{cc} - \sin \theta & - \cos \theta \\ \cos \theta & - \sin \theta \end{array} \right),
\]
we obtain
\[
\det W' = \det \bigl( [ e^{\theta J_0}, A' A^{-1} ] + \theta' J_0 e^{\theta J_0} \bigr) =  {\theta'}^2  - \bigl( 4a^2 + (b+c)^2 \bigr) \sin^2 \theta.\qedhere
\]
\end{proof}

By (\ref{form}), the Krein-positive eigenvalue of $W(t)$ is $e^{i \theta(t)}$ and we can find a differentiable curve $u:(0,1) \rightarrow \C^2$ of eigenvectors of $W$ corresponding to this eigenvalue such that $\kappa(u,u)=1$. Since $W$ is timelike, the quantity $\omega_0(u,W'W^{-1}u)$ is positive and the formula of Lemma \ref{fundlem} implies that $\theta'> 0$ on $(0,1)$. We can then use the identity from Lemma \ref{det} and obtain the desired upper bound for $\mathrm{length}_G(W)$:
\begin{equation}
\label{laine}
\begin{split}
\mathrm{length}_G(W) &= \int_0^1 (\det W'(t))^{\frac{1}{2}}  \, \di t = \int_0^1 \bigl( {\theta'}^2  - \bigl( 4a^2 + (b+c)^2 \bigr) \sin^2 \theta \bigr)^{\frac{1}{2}} \, \di t \\ &\leq \int_0^1 |{\theta'}| \, \di t = \int_0^1 \theta'\, \di t = \theta(1).
\end{split}
\end{equation}
This concludes the proof of statement (i) in Proposition \ref{lordisSp2}. 

\section{Proof of Theorem \ref{Lordist}} 
\label{dimLordist}

Let $W:[0,1] \rightarrow \mathrm{Sp}(2n)$ be a causal curve such that $W(t)$ is elliptic for every $t\in [0,1]$, meaning that all the eigenvalues of $W(t)$ belong to the unit circle $\U$. By the continuous dependence of the spectrum, see \cite[Theorem II.5.1]{kat80}, the spectrum of $W(t)$ is given by
\[
\sigma(W(t)) = \{ e^{\pm i \theta_1(t)}, \dots, e^{\pm i \theta_n(t)} \}
\]
for some continuous functions $\theta_j: [0,1] \rightarrow \R$, $j=1,\dots,n$. Moreover, we can assume that the Krein-positive eigenvalues of $W(t)$ are the eigenvalues
\[
e^{i \theta_1(t)}, \dots, e^{ i \theta_n(t)}.
\]
Here, we are counting eigenvalues according to their algebraic multiplicity and are seeing an eigenvalue on $\U$ of Krein signature $(p,q)$ as $p$ Krein-positive eigenvalues and $q$ Krein-negative ones. The fact that the curve $W$ is causal implies that all the functions $\theta_j$ are non-decreasing, see Lemma \ref{fundlem} above for the case of simple eigenvalues and \cite[Corollary I.3.5]{eke90} for the general case.

The proof of statement (i) in Theorem \ref{Lordist} from the Introduction is based on the following result.

\begin{prop}
Let $W:[0,1] \rightarrow \mathrm{Sp}(2n) $ be a causal curve taking values in the closure of $\mathrm{Sp}_{\mathrm{ell}}^+(2n)$ and denote by
\[
\sigma(W(t)) = \{ e^{\pm i \theta_1(t)}, \dots, e^{\pm i \theta_n(t)} \}
\]
the spectrum of $W(t)$, where $\theta_j :[0,1] \rightarrow [0,\pi]$ are continuous functions. Then
\[
\mathrm{length}_G(W) \leq \frac{1}{n} \sum_{j=1}^n \bigl( \theta_j(1) - \theta_j(0) \bigr).
\]
\end{prop}

\begin{proof}
The boundary of $\mathrm{Sp}_{\mathrm{ell}}^+(2n)$ consists of endomorphisms having the eigenvalue 1 or -1, so $W(t)$ belongs to this boundary if and only if $\theta_j(t)$ has the value $0$ or $\pi$ for at least one $j$. Since each function $\theta_j$ is non-decreasing and takes values in $[0,\pi]$, $W$ can take values in the boundary of $\mathrm{Sp}_{\mathrm{ell}}^+(2n)$ only on two closed subintervals of the form $[0,t_0]$ and $[t_1,1]$. When one of these intervals has positive length, the restriction of $W$ to it is lightlike and hence this interval gives no contribution to $\mathrm{length}_G(W)$. Therefore,
\[
\mathrm{length}_G(W) = \mathrm{length}_G(W|_{[t_0,t_1]}) = \sup_{\epsilon>0} \mathrm{length}_G(W|_{[t_0+\epsilon,t_1-\epsilon]}).
\]
By this observation and by the monotonicity of the functions $\theta_j$, it is enough to prove the desired upper bound on $\mathrm{length}_G(W)$ in the case of a causal curve $W:[0,1] \rightarrow \mathrm{Sp}(2n)$ which takes values in $\mathrm{Sp}_{\mathrm{ell}}^+(2n)$. Moreover, any such curve can be $C^1$-approximated by a curve which is timelike and analytic. Therefore, in the following we can assume that the path $W$ is timelike, analytic and takes values in the open set $\mathrm{Sp}_{\mathrm{ell}}^+(2n)$. 

By the last condition, $W(t)$ is diagonalizable for every $t\in [0,1]$, so we can find a basis $u_1(t),\ldots,u_n(t),\overline{u}_1(t),\ldots,\overline{u}_n(t)$ of $\C^{2n}$ which is $\kappa$-unitary and satisfies
\begin{equation}
\label{autovec}
W(t) u_j(t)  = e^{i \theta_j(t)} u_j(t), \qquad W(t) \overline{u}_j(t)  = e^{-i \theta_j(t)} \overline{u}_j(t), \qquad \forall t\in [0,1].
\end{equation}
The functions $\theta_j$ are differentiable (see \cite[Theorems II.5.4 and II.5.6]{kat80}), but in general the eigenvectors $u_j(t)$ need not depend continuously on $t$. Indeed, one may lose continuity of the eigenvectors when different eigenvalues collide. However, by the analyticity of $W$ the set $\mathcal{T}$ of exceptional instants at which some eigenvalues change their multiplicity is finite, and the functions $u_j$ are analytic on $[0,1] \setminus \mathcal{T}$ (see \cite[Section II.1.4]{kat80}). Then Lemma \ref{fundlem} gives us the identity
\begin{equation}
\label{fund}
\theta_j' = \omega_0(u_j,W'W^{-1}  u_j) \qquad \mbox{on } [0,1] \setminus \mathcal{T}.
\end{equation}
The fact that $W$ is timelike tells us that the endomorphism $W'(t)W^{-1}(t)$ is in $\mathrm{sp}^+(2n)$ for every $t\in [0,1]$ and hence
\[
\omega_0( u_j,W'W^{-1}  u_j) > 0  \qquad \mbox{on } [0,1].
\]
Therefore, (\ref{fund}) implies that the functions $\theta_j$ are strictly increasing on $[0,1]$. 

By Proposition \ref{propA1}, the spectrum of $W'(t)W^{-1}(t) \in  \mathrm{sp}^+(2n)$ has the form
\[
\{ \pm i \lambda_1(t), \ldots , \pm i \lambda_n(t) \}
\]
for some continuous positive functions $\lambda_j : [0,1] \rightarrow \R$, and we can find a $\kappa$-unitary basis $v_1(t),\dots,v_n(t),\overline{v}_1(t), \dots,\overline{v}_n(t)$ of $\C^{2n}$ such that
\[
W'(t)W^{-1}(t) v_j(t) = i \lambda_j(t) v_j(t), \qquad W'(t)W^{-1}(t) \overline{v}_j(t) = -i \lambda_j(t) \overline{v}_j(t) \qquad \forall t\in [0,1].
\]
We now express each vector $u_j\in \C^{2n}$ as a linear combination of the latter basis:
\[
u_j = \sum_{h=1}^n ( \alpha_{jh} v_h + \beta_{jh} \overline{v}_h) \qquad \forall j,
\]
for suitable complex numbers $\alpha_{jh},\beta_{jh}$. The fact that both bases $u_1,\dots,u_n,\overline{u}_1,\dots,\overline{u}_n$ and  $v_1,\dots,v_n,\overline{v}_1,\dots,\overline{v}_n$ are $\kappa$-unitary implies the identities
\[
\delta_{jk} = \kappa(u_j,u_k) = \sum_{h=1}^n ( \alpha_{jh} \overline{\alpha}_{kh} - \beta_{jh} \overline{\beta}_{kh}) \qquad \forall j,k.
\]
If $A:=(\alpha_{jh})$ and $B:=(\beta_{jh})$ are the $n\times n$ complex matrices given by these coefficients, the above identity can be rewritten more compactly as
\begin{equation}
\label{matide}
AA^* - BB^* = I.
\end{equation}
In particular, the self-adjoint matrix $AA^*$ satisfies $AA^*\geq I$. But then we also have $A^*A\geq I$, because this inequality can be read from the spectrum, and $A^*A$ and $AA^*$ have the same spectrum. In particular, the $h$-th diagonal element of the matrix $A^*A$ is in absolute value not smaller than 1:
\begin{equation}
\label{unitary}
(A^*A)_{hh} = \sum_{j=1}^n (A^*)_{hj} (A)_{jh} = \sum_{j=1}^n \overline{\alpha}_{jh} \alpha_{jh} =  \sum_{j=1}^n |\alpha_{jh}|^2 \geq 1.
\end{equation}
We now rewrite (\ref{fund}) in terms of the $\kappa$-unitary basis $v_1(t),\dots,v_n(t),\overline{v}_1(t), \dots,\overline{v}_n(t)$ and get
\[
\begin{split}
\theta_j' &=  - i \kappa\Bigl( W^{-1} W'   \sum_{h=1}^n ( \alpha_{jh} v_h + \beta_{jh} \overline{v}_h),  \sum_{h=1}^n ( \alpha_{jh} v_h + \beta_{jh} \overline{v}_h) \Bigr) \\ &=  - i \kappa\Bigl( i \sum_{h=1}^n ( \lambda_h \alpha_{jh} v_h - \lambda_h \beta_{jh} \overline{v}_h),  \sum_{h=1}^n ( \alpha_{jh} v_h + \beta_{jh} \overline{v}_h) \Bigr) \\ &= \sum_{h=1}^n ( |\alpha_{jh}|^2 + |\beta_{jh}|^2) \lambda_h \geq \sum_{h=1}^n  \lambda_h |\alpha_{jh}|^2,
\end{split}
\]
on $[0,1]\setminus \mathcal{T}$. By adding over $j$ and using (\ref{unitary}) we obtain
\[
\sum_{j=1}^n \theta_j' \geq \sum_{j=1}^n \sum_{h=1}^n \lambda_h  |\alpha_{jh}|^2 = \sum_{h=1}^n \lambda_h \sum_{j=1}^n |\alpha_{jh}|^2 \geq \sum_{h=1}^n \lambda_h  \qquad \mbox{on } [0,1] \setminus \mathcal{T}.
\]
Integration over $[0,1]$ gives us the chain of inequalities
\[
\begin{split}
 \sum_{j=1}^n \bigl(\theta_j(1) - \theta_j(0)\bigr) &= \sum_{j=1}^n \int_0^1 \theta_j'(t)\, \di t \geq \int_0^1 \sum_{h=1}^n \lambda_h(t)\, \di t \geq n \int_0^1 \left( \prod_{h=1}^n \lambda_h(t) \right)^{\frac{1}{n}} \, \di t \\ &= n \int_0^1 \bigl( \det (W'(t)W^{-1}(t) )\bigr)^{\frac{1}{2n}} \, \di t = n\, \mathrm{length}_G(W).
\end{split}
\]
where at the end of the first line we have used the inequality between the arithmetic and geometric means. This proves the bound
\[
 \mathrm{length}_G(W) \leq \frac{1}{n}  \sum_{j=1}^n \bigl(\theta_j(1) - \theta_j(0)\bigr)
 \]
and concludes the proof.
\end{proof}

Statement (i) of Theorem  \ref{Lordist} from the Introduction is an immediate consequence of the above proposition. Statement (ii) of that theorem follows immediately from the next proposition.

\begin{prop} 
\label{verylong}
Let $W:[0,1] \rightarrow \mathrm{Sp}(2n)$ be a timelike curve such that $W(0)=\mathrm{id}$ and $W(1)$ is not in the closure of $\mathrm{Sp}_{\mathrm{ell}}^+(2n)$. Then there are timelike curves from $\mathrm{id}$ to $W(1)$ which are homotopic with fixed ends to $W$ and have arbitrarily large Lorentz--Finsler length $\mathrm{length}_G$.
\end{prop}

\begin{proof} 
Denote by $t^*\in (0,1)$ the first positive instant at which $W$ reaches the boundary of $\mathrm{Sp}_{\mathrm{ell}}^+(2n)$. Equivalently, the number $t^*$ is the first instant at which $W$ hits the singular hypersurface consisting of elements of $\mathrm{Sp}(2n)$ having the eigenvalue $-1$. Up to a $C^1$-small perturbation keeping $W$ timelike and within the same homotopy class, we may assume that the eigenvalue $-1$ of $W(t^*)$ has algebraic multiplicity 2. Then $W(t^*)$ splits as a symplectic automorphism of a symplectic plane and a symplectic automorphism of a symplectic $(2n-2)$-dimensional subspace, and the same is true for $t$ close to $t^*$. More precisely, there exists an interval $[t^-,t^+]\subset [0,1]$ containing $t^*$ in its interior such that
	\[
	W(t) = \Phi(t)^{-1} \bigl( W_0(t) \oplus W_1(t) \bigr) \Phi(t) \qquad \forall t\in  [t^-,t^+],
	\]
	for some smooth curves $W_0:[t^-,t^+] \rightarrow \mathrm{Sp}(2)$, $W_1: [t^-,t^+] \rightarrow \mathrm{Sp}(2n-2)$, $\Phi: [t^-,t^+] \rightarrow \mathrm{Sp}(2n)$. Here, $W_1$ takes values in $\mathrm{Sp}_{\mathrm{ell}}^+(2n-2)$ and $W_0(t)$ is in $\mathrm{Sp}_{\mathrm{ell}}^+(2)$ for $t<t^*$, on the boundary of this set for $t=t^*$, and outside of the closure of $\mathrm{Sp}_{\mathrm{ell}}^+(2)$ for $t>t^*$.

For all $k\in\N$, the set $U_k:=\{\mathrm{id}\}\cup\mathrm{Sp}_{\mathrm{ell}}^+(2k)$ is contractible as it is the homeomorphic image under the exponential map of $\{0\}\cup \mathrm{sp}_{\mathrm{ell}}^+(2k)$ (see Proposition \ref{propA2} in Appendix \ref{linear-app}). Therefore we can find timelike curves
\[
\widetilde{W}_0: [0,t^-] \rightarrow U_2, \qquad \widetilde{W}_1: [0,t^-] \rightarrow U_{2n-2},
\]
such that $\widetilde{W}_0(0)=\mathrm{id}$, $\widetilde{W}_0(t^-)=W_0(t^-)$, $\widetilde{W}_1(0)=\mathrm{id}$, $\widetilde{W}_1(t^-)=W_1(t^-)$ and the image of the timelike curve
\[
t\mapsto \Phi(t^-)^{-1} \bigl( \widetilde{W}_0(t) \oplus \widetilde{W}_1(t) \bigr) \Phi(t^-), \qquad t\in [0,t^-],
\]
is contained in $U_{2n}$. Since $U_{2n}$ is contractible, the above curve is homotopic with fixed ends to the restriction $W|_{[0,t^-]}$. Homotoping also the conjugacy from $\Phi(t^-)$ to $\Phi(t^+)$, we see that $W|_{[0,t^+]}$ is homotopic with fixed ends to
\[
\Phi(t^+)^{-1} \bigl( (\widetilde{W}_0\# W_0|_{[t^-,t^+]}) \oplus (\widetilde{W}_1\# W_1|_{[t^-,t^+]}) \bigr) \Phi(t^+),
\]
where $\#$ denotes concatenation of paths. Since $\widetilde{W}_1\# W_1|_{[t^-,t^+]}$ is contained in $U_{2n-2}$ and $W_1(t^+)\neq\mathrm{id}$, this path is homotopic with fixed ends to a timelike curve $\widehat W_1\colon [0,t^+]\to U_{2n-2}$ for which we can assume that $G(\widehat{W}_1')=(\det \widehat{W}_1')^\frac{1}{2n-2}$ is constant. If $\mu$ is the Maslov quasi-morphism, then $\mu(\widetilde{W}_0\# W_0|_{[t^-,t^+]})\geq\tfrac12$ and $W_0(t^+)$ is not in the closure of $\mathrm{Sp}_{\mathrm{ell}}^+(2)$. By statement (ii) in Proposition \ref{lordisSp2}, for every $c>0$ there is a timelike curve $\widehat{W}_0: [0,t^+]\rightarrow \mathrm{Sp}(2)$ which is homotopic with fixed ends to $\widetilde{W}_0\# W_0|_{[t^-,t^+]}$ and satisfies
\[
\mathrm{length}_G(\widehat{W}_0) \geq c.
\]
Without loss of generality we can assume that $G(\widehat{W}_0')=(\det \widehat{W}_0')^{\frac12}$ is constant and equals $c/t^+$. Then the timelike curve
\[
\widehat{W}\colon[0,t^+]\to \mathrm{Sp}(2n),\qquad \widehat{W}:= \Phi(t^+)^{-1}(\widehat{W}_0 \oplus \widehat{W}_1)\Phi(t^+)
\] is homotopic with fixed ends to $W|_{[0,t^+]}$ and, by the bi-invariance of $G$, has length
\[
\mathrm{length}_G(\widehat{W}) = t^+\bigl(\det \widehat{W}_0' \oplus \widehat{W}_1' \bigr)^{\frac{1}{2n}} = t^+\Big(\frac{c}{t^+}\Big)^{\frac1n}\bigl(\det \widehat{W}_1' \bigr)^{\frac{n-1}{n}}
\] 
which can be made arbitrarily large for $c$ arbitrarily large. Hence the concatenation $\widehat{W}\# (W|_{[t^+,1]})$ is a timelike curve homotopic with fixed ends to $W$ and with arbitrarily large Lorentz--Finsler length. 
\end{proof}

\section{Proof of Theorem \ref{nolordisCont}} 
\label{nolordisContsec}

In order to prove Theorem \ref{nolordisCont}, it is convenient to work with the Hamiltonian formalism instead of the contact one. By seeing $\R \mathrm{P}^{2n-1}$ as the manifold of lines through the origin in $\R^{2n}$, the contactomorphism group of $\R \mathrm{P}^{2n-1}$ can be identified with the group $\mathcal{G}_n$ of all 1-homogeneous symplectomorphisms of $\R^{2n}\setminus \{0\}$, i.e., the group of diffeomorphisms of $\R^{2n}\setminus \{0\}$ which preserve the Liouville 1-form $\lambda_0$ and are equivariant with respect to the antipodal $\Z_2$-action $z\mapsto -z$. These maps extend to homeomorphisms of $\R^{2n}$, but they are not differentiable at the origin, unless they are linear. The Lie algebra $\mathrm{cont}(\R \mathrm{P}^{2n-1},\xi_{\mathrm{st}})$ is then identified with the space of 1-homogeneous Hamiltonian vector fields on $\R^{2n}$. These are precisely the Hamiltonian vector fields that are induced by Hamiltonian functions in the space of 2-homogeneous functions
\[
\mathcal{H}_n := \{ H: \R^{2n} \rightarrow \R \mid H \mbox{ smooth on } \R^{2n}\setminus \{0\}, \; H(tz) = t^2 H(z) \; \forall t\in \R, \; \forall z\in \R^{2n} \}.
\] 
Hamiltonians in $\mathcal{H}_n$ are differentiable at $0$ with differential which is Lipschitz-continuous on $\R^{2n}$, but in general are not twice differentiable at $0$. 
Given $H\in \mathcal{H}_n$, we denote by $X_H$ the induced \textit{Hamiltonian vector field} on $\R^{2n}$, which in the symplectic setting is defined by
\[
\omega_0(X_H(z),v) = -\mathrm{d}H(z)\cdot v \qquad \forall z,v\in \R^{2n}.
\] 
This vector field is Lipschitz-continuous on $\R^{2n}$ and smooth on $\R^{2n}\setminus \{0\}$. The path in $\mathcal{G}_n$ which is generated by the Hamiltonian vector field of a (in general time-dependent) Hamiltonian $H$ is denoted by $\phi_H^t$.

The bijection between elements of $\mathcal{H}_n$ and contact Hamiltonians on $(\mathbb{R}^{2n-1},\xi_{\mathrm{st}})$ associated to the standard contact form $\alpha_0$ is given just by the restriction of the even function $H\in \mathcal{H}_n$ to $S^{2n-1}$. Therefore, a path in $\mathrm{Cont}(\mathbb{R}^{2n-1},\xi_{\mathrm{st}})$ is positive if and only if the corresponding path in $\mathcal{G}_n$ is generated by a time-dependent Hamiltonian $\{H_t\}_{t\in [0,1]}$ such that $H_t>0$ on $\R^{2n}\setminus \{0\}$ for every $t\in [0,1]$. In this case, the corresponding path in $\mathcal{G}_n$ is also called {\it positive}.

\begin{lem}
\label{posconju}
Let $\psi\in \mathcal{G}_n$ and $L\in \mathcal{H}_n$. Then the path $t\mapsto \phi_L^{-t}\circ \psi \circ \phi_L^t$ in $\mathcal{G}_n$ is positive if and only if
\[
L - L\circ \psi > 0 
\]
on $\R^{2n}\setminus \{0\}$.
\end{lem}

\begin{proof}
Differentiating this path we obtain
\[
\frac{\mathrm{d}}{\mathrm{d}t}  \phi_L^{-t}\circ \psi \circ \phi_L^t = - X_L + (\phi_L^{-t}\circ \psi)_* X_L = X_{-L+(\phi_L^{-t}\circ \psi)_* L},
\]
where all vector fields are evaluated at $ \phi_L^{-t}\circ \psi \circ \phi_L^t$. Therefore, the path is generated by the time-dependent Hamiltonian
\[
H_t = -L + L\circ \psi^{-1}\circ \phi_L^t = (-L+L\circ \psi^{-1})\circ \phi_L^t,
\]
where in the last identity we have used the fact that $L$ is invariant under the flow $\phi^t_L$. This $t$-dependent Hamiltonian is positive  on $\R^{2n}\setminus \{0\}$ if and only if  $L\circ \psi^{-1}-L$ is, which is equivalent to the condition we stated.
\end{proof}

\begin{rem}
{\rm The above lemma can be restated in the setting of contact Hamiltonians as follows: If $\phi_L^t$ is the flow of the autonomous contact vector field on $(M,\xi)$ defined by the contact Hamiltonian $L:M \rightarrow \R$ with respect to some contact form $\alpha$ defining $\xi$ and $\psi$ is any contactomorphism of $(M,\xi)$, then the path $t\mapsto \phi_L^{-t}\circ \psi \circ \phi_L^t$ is positive in $\mathrm{Cont}(M,\xi)$ if and only if
\[
f_{\psi} L - L\circ \psi > 0 
\]
on $M$, where $f_{\psi}$ is the conformal factor defined by $\psi^* \alpha = f_{\psi} \alpha$.}
\end{rem}

Our proof of Theorem \ref{nolordisCont} is based on the following result.

\begin{prop}
\label{exposconju}
For every $\epsilon>0$ there exists a smoothly time-dependent Hamiltonian $\{H_t\}_{t\in [0,1]}\subset \mathcal{H}_n$ with
\[
0 < H_t(z) \leq \epsilon |z|^2 \qquad \forall t\in [0,1], \; \forall z\in \R^{2n}\setminus \{0\},
\]  
and an autonomous Hamiltonian $L\in \mathcal{H}_n$ such that
\[
L - L\circ \phi_H^1 >0
\]
on $\R^{2n}\setminus \{0\}$.
\end{prop}

Before discussing the proof of this proposition, we show how it implies Theorem \ref{nolordisCont}. The argument is analogous to the proof of Proposition \ref{trivlordisSp}.

\begin{proof}[Proof of Theorem \ref{nolordisCont}]
Let $\phi_0$ and $\phi_1$ be elements of $\widetilde{\mathrm{Cont}}_0(\R \mathrm{P}^{2n-1},\xi_{\mathrm{st}})$. If there is no non-negative and somewhere positive path from $\phi_0$ to $\phi_1$, then $d(\phi_0,\phi_1)=0$ by assumption (i). Therefore, we must prove that if there is a non-negative and somewhere positive path from $\phi_0$ to $\phi_1$, then $d(\phi_0,\phi_1)=+\infty$. 

By the bi-invariance of $d$ (assumption (iii)), we may assume that $\phi_0=\mathrm{id}$. Let $0\leq t_0 < t_1 \leq 1$ be such that the path $\{\phi_t\}_{t\in [t_0,t_1]}$ is positive. Then $\phi_{t_1}\circ \phi_{t_0}^{-1}$ is generated by a positive Hamiltonian and we denote by $\epsilon>0$ its minimum. By Lemma \ref{posconju} and Proposition \ref{exposconju}, there exists a smoothly time-dependent contact Hamiltonian $\{H_t\}_{t\in [0,1]}\subset C^{\infty} (\R \mathrm{P}^{2n-1})$ such that
\[
0 < H \leq \epsilon \qquad \mbox{on } [0,1]\times \R \mathrm{P}^{2n-1}
\]
and an autonomous contact Hamiltonian $L\in C^{\infty} (\R \mathrm{P}^{2n-1})$ such that, denoting by $\phi_H^t$ and $\phi_L^t$ the generated paths in $\widetilde{\mathrm{Cont}}_0(\R \mathrm{P}^{2n-1},\xi_{\mathrm{st}})$, we have that 
\begin{equation}
\label{thepath}
t\mapsto \phi_L^{-t}\circ \phi_H^1\circ  \phi_L^t
\end{equation}
is a positive path. By the above bounds on $H$, we have
\[
\mathrm{id} \leq \phi_H^1 \leq \phi_{t_1}\circ \phi_{t_0}^{-1} \leq \phi_1,
\]
so by the reverse triangular inequality (ii) it is enough to show that $d(\mathrm{id}, \phi_H^1)=+\infty$. By the properties (iii) and (ii) of $d$ we have
\[
d(\mathrm{id}, \phi_H^1) = d(\mathrm{id}, \phi_L^{-1}\circ \phi_H^1\circ  \phi_L^1) \geq d(\mathrm{id}, \phi_H^1) + d(\phi_H^1, \phi_L^{-1}\circ \phi_H^1\circ  \phi_L^1).
\]
The latter term $d(\phi_H^1, \phi_L^{-1}\circ \phi_H^1\circ  \phi_L^1)$ is positive because of condition (i) and the positivity of the path (\ref{thepath}). Then the above inequality forces
\[
d(\mathrm{id}, \phi_H^1) = +\infty,
\]
concluding the proof.
\end{proof}

We conclude this section by proving Proposition \ref{exposconju}. By identifying $\mathrm{Diff}_0(\T) = \mathrm{Cont}_0(\R \mathrm{P}^1,\xi_{\mathrm{st}})$ with $\mathcal{G}_1$, we have for every $k\in \N$ the Lie group homomorphism
\[
j_k : \mathrm{PSp}_k(2) \rightarrow \mathcal{G}_1.
\]
We identify the Lie algebra $\mathrm{sp}(1)$ with the space of quadratic forms on $\R^2$ by mapping $J_0S$, $S\in \mathrm{Sym}(2)$, to the quadratic form
\[
h(z) = \frac{1}{2} S z \cdot z, \qquad \forall z\in \R^2.
\]
We write the quadratic form $h$ in polar coordinates as
\[
h(re^{i\theta}) = \frac{r^2}{2} \widehat{h} (2\theta), \qquad \forall r\geq 0, \; \forall \theta\in \R/2\pi \Z,
\]
where $\widehat{h}$ is a trigonometric polynomial of degree 1. By the results of Section \ref{sottogruppisec}, the linear mapping
\[
\di j_k( \mathrm{id}) : \mathrm{sp}(2) \rightarrow \mathcal{H}_1
\]
maps the quadratic form $h$ to the 2-homogeneous function
\[
h_k \in \mathcal{H}_1, \qquad h_k(re^{i\theta}) := \frac{r^2}{2k} \, \widehat{h} (2k\theta), \qquad \forall r\geq 0, \; \forall \theta\in \R/2\pi \Z.
\]

Fix some timelike curve $w: [0,1] \rightarrow \mathrm{PSp}(2)$ such that $w(0)=\mathrm{id}$, $w(1)$ is hyperbolic, and let $h$ be the smooth path of quadratic forms on $\R^2$ generating $w$. Given $k\in\N$, we denote by $w_k:[0,1] \rightarrow \mathrm{PSp}_k(2)$ the lift of $w$ such that $w_k(0)=\mathrm{id}$. Then the positive path $j_k\circ w_k : [0,1] \rightarrow \mathcal{G}_1$ is generated by the positive time-dependent Hamiltonian $h_k: [0,1] \rightarrow \mathcal{H}_1$ which satisfies
\begin{equation}
\label{ntd1}
h_k(t,z) \leq \frac{c}{k} |z|^2, \qquad \mbox{where } c:= \max_{(t,z)\in [0,1]\times S^1} h(t,z).
\end{equation}
Since $w(1)$ is hyperbolic, by Lemma \ref{algebraic} there exists $S\in \mathrm{Sym}(2)$ such that the curve
\begin{equation}
\label{ntd2}
t\mapsto e^{-t J_0 S} w(1) e^{tJ_0 S}
\end{equation}
is timelike in $\mathrm{PSp}(2)$. The autonomous path $t\mapsto e^{tJ_0 S}$ is generated by the quadratic Hamiltonian
\[
\ell(z) := \frac{1}{2} S z\cdot z, \qquad \forall z\in \R^2.
\]
Then the image by $j_k$ of the lift to $\mathrm{PSp}_k(2)$ of this path is the autonomous path in $\mathcal{G}_1$ which is generated by the Hamiltonian $\ell_k\in \mathcal{H}_1$. The fact that the curve (\ref{ntd2}) is timelike implies that the path
\[
t\mapsto j_k \bigl( e^{-tJ_0 S} w_k(1) e^{tJ_0 S} \bigr) = \phi_{\ell_k}^{-t} \circ \phi_{h_k}^1 \circ \phi_{\ell_k}^t
\]
is positive in $\mathcal{G}_1$. By Lemma \ref{posconju}, this implies that
\begin{equation}
\label{ntd3}
\ell_k - \ell_k \circ \phi_{h_k}^1 >0 \qquad \mbox{on } \R^2\setminus \{0\}.
\end{equation}
Consider the 2-homogeneous functions
\[
\begin{split}
H_k: [0,1]\times \R^{2n} \rightarrow \R, \qquad & H_k(t,z_1,\dots,z_n) := \sum_{j=1}^n h_k(t,z_j), \\
L_k: \R^{2n} \rightarrow \R, \qquad & L_k(z_1,\dots,z_n) := \sum_{j=1}^n \ell_k(z_j).
\end{split}
\]
These functions are smooth on the open set where all $z_j$ do not vanish, but globally they are just of class $C^{1,1}$, which as usual denotes the class of continuously differentiable functions whose differential is Lipschitz-continuous. The Hamiltonian vector field $X_{H_k}$ is Lipschitz-continuous and its non-autonomous flow is the following path of bi-Lipschitz-continuous homeomorphisms
\[
\phi_{H_k}^t(z_1,\dots,z_n) = \bigl( \phi_{h_k}^t(z_1), \dots, \phi_{h_k}^t(z_n) \bigr).
\]
Therefore, (\ref{ntd3}) implies
\begin{equation}
\label{ntd4}
L_k - L_k \circ \phi_{H_k}^1 >0 \qquad \mbox{on } \R^{2n}\setminus \{0\}.
\end{equation}
Moreover, from the positivity of $h_k$ and (\ref{ntd1}) we deduce
\begin{eqnarray}
\label{ntd5}
H_k > 0 \qquad \mbox{on } [0,1] \times \bigl( \R^{2n} \setminus \{0\} \bigr), \\
\label{ntd6}
H_k(t,z) \leq \frac{c}{k} |z|^2 \qquad \forall (t,z) \in [0,1]\times \R^{2n}.
\end{eqnarray}
We now fix $k\in \N$ large enough, so that
\begin{equation}
\label{ntd7}
\frac{c}{k} \leq \frac{\epsilon}{2},
\end{equation} 
where $\epsilon$ is the positive number appearing in the statement of Proposition \ref{exposconju}. The last step is to approximate $H_k$ and $L_k$ by 2-homogeneous Hamiltonians which are smooth on the whole $\R^{2n} \setminus \{0\}$. Here is the standard argument for such an approximation. By (\ref{ntd4}),  (\ref{ntd5}),  (\ref{ntd6}) and  (\ref{ntd7}) we have
\[
\begin{split}
\eta \leq H_k(t,z) &\leq \frac{\epsilon}{2} \qquad \forall (t,z) \in [0,1]\times S^{2n-1}, \\
L_k(z) - L_k\circ \phi_{H_k}^1(z) &\geq \eta \qquad \forall z\in S^{2n-1},
\end{split}
\]
for some $\eta>0$. For every $\delta>0$ we can find a smoothly time-dependent function $\{H(t,\cdot)\}_{t\in [0,1]}\subset \mathcal{H}_n$ and $L\in \mathcal{H}_n$ such that
\[
\|H(t,\cdot) - H_k(t,\cdot)\|_{C^{1,1}(S^{2n-1})} < \delta, \qquad \|L-L_k\|_{C^0(S^{2n-1})} < \delta.
\]
The first bound implies that $X_H$ is $C^{0,1}$-close to $X_{H_k}$ on compact subsets of $[0,1]\times \R^{2n}$ and hence $\phi^1_H$ is $C^0$-close to $\phi^1_{H_k}$ on compact subsets of $\R^{2n}$. By choosing $\delta$ small enough, we can then ensure the bounds
\[
\frac{\eta}{2} \leq H(t,z) \leq \epsilon \qquad \forall (t,z) \in [0,1]\times S^{2n-1},
\]
and
\[
\begin{split}
L - L \circ \phi_H^1 &= L_k -  L_k \circ \phi_{H_k}^1 + L - L_k + L_k \circ  \phi_{H_k}^1 - L_k \circ \phi_H^1 + (L_k-L)\circ \phi_H^1 \\ &\geq \eta - |L-L_k| - | L_k \circ  \phi_{H_k}^1 - L_k \circ \phi_H^1| - | (L_k-L)\circ \phi_H^1| \geq \frac{\eta}{2},
\end{split}
\]
on $S^{2n-1}$, where we are using the fact that $L_k$ is uniformly continuous on compact subsets of $\R^{2n}$. The functions $H$ and $L$ satisfy the requests of Proposition \ref{exposconju}, which is then proven. 

\section{Proof of Theorem \ref{quantum}}

In this section, we prove Theorem \ref{quantum} from the Introduction and discuss a related example. The proof of statement (i) in Theorem \ref{quantum} uses the following result of Nazarov \cite{naz21}.

\begin{thm} 
\label{bernstein}
Let $k\in \N$. If the trigonometric polynomial
\[
p(x) = \sum_{j=0}^k \bigl( a_j \cos (2\pi j x) + b_j \sin (2\pi j x) \bigr), \qquad a_j,b_j\in \R,
\]
is non-negative, then
\[
\int_0^1 |p'(x)| \, \mathrm{d} x \leq 4k \int_0^1 p(x) \, \mathrm{d} x.
\]
\end{thm}

Note that without the non-negativity assumption, the $L^1$-norm of the derivative of the trigonometric polynomial $p$ on $[0,1]$ would be bounded by $2\pi k$ times the $L^1$-norm of $p$, by the classical Bernstein inequality. Thanks to the assumption $p\geq 0$, the constant $2\pi$ can be replaced by $4$, which is optimal, as the example of
\[
p(x) = 1 + \cos (2\pi k x)
\]
shows.

In order to prove statement (i), we consider for a given number $k\in \N$ and positive Hamiltonian 
$H\in \mathcal{P}_k$ the path $\{\phi^t\}_{t\in [0,1]}$ in $\mathrm{Diff}_1(\R)$ which is defined by
\begin{equation}
\label{nonaut}
\frac{\mathrm{d}}{\mathrm{d}t} \phi^t = H(t,\phi^t), \qquad \phi^0 = \mathrm{id}.
\end{equation}
Assuming that
\begin{equation}
\label{ass}
\phi^1(x)\leq x + \frac{s}{4k}
\end{equation}
for some $s<1$, we must prove the bound
\[
\mathrm{length}_V \bigl( \{\phi^t\}_{t\in [0,1]} \bigr) \leq \frac{1}{1-s} \int_0^1 \bigl( \phi^1(x) - x \bigr)\, \mathrm{d}x.
\]
By multiplying both sides of \eqref{nonaut} by $\partial_x \phi^t$ and integrating over $[0,1]$ after a change of variable, we find
\begin{equation}
\label{prima}
\int_0^1 H(t,x)\, \mathrm{d}x = \int_0^1 H(t,\phi^t(x)) \partial_x \phi^t(x) \, \mathrm{d}x = \int_0^1  \partial_t \phi^t(x) \partial_x \phi^t(x)\, \mathrm{d}x.
\end{equation}
Writing
\begin{equation}
\label{lapsi}
\phi(t,x) = \phi^t(x) = x + \psi(t,x),
\end{equation}
where the function $\psi$ is 1-periodic in $x$, non-negative and monotonically increasing in $t$, we manipulate the right-hand side of \eqref{prima} as follows:
\begin{equation}\label{unobis}
\begin{split}
\int_0^1  \partial_t \phi^t(x) \partial_x \phi^t(x)\, \mathrm{d}x &=  \int_0^1 \partial_t \phi(t,x) \, \mathrm{d}x +  \int_0^1 \partial_t \phi(t,x)  \partial_x \psi(t,x)\, \mathrm{d}x \\
&= \int_0^1 \partial_t \phi(t,x) \, \mathrm{d}x + \int_0^1  H(t, \phi(t,x))  \partial_x \psi(t,x)\, \mathrm{d}x \\ &=  \int_0^1 \partial_t \phi(t,x) \, \mathrm{d}x - \int_0^1 \partial_x \bigl( H(t,\phi(t,x)) \bigr)  \, \psi(t,x)\, \mathrm{d}x \\ &=  \int_0^1 \partial_t \phi(t,x) \, \mathrm{d}x - \int_0^1 \partial_x H(t,\phi^t(x))  \partial_x \phi^t(x)  \, \psi(t,x)\, \mathrm{d}x \\ &=  \int_0^1 \partial_t \phi(t,x) \, \mathrm{d}x - \int_0^1 \partial_x H(t,x)   \, \psi(t,(\phi^t)^{-1}(x))\, \mathrm{d}x.
\end{split}
\end{equation}
We estimate the latter term using Theorem \ref{bernstein} and obtain
\[
\begin{split}
\Big|\int_0^1 \partial_x H(t,x)   \, \psi(t,(\phi^t)^{-1}(x))\, \mathrm{d}x \Big|&\leq \sup_{x\in [0,1]} |\psi(t,x)|\int_0^1 |\partial_x H(t,x)|\, \mathrm{d}x\\
&\leq 4k  \sup_{x \in [0,1]} |\psi(1,x)|\int_0^1 H(t,x)\, \mathrm{d}x\\
&=4k  \sup_{x\in [0,1]} (\phi^1(x)-x )\int_0^{1}H(t,x)\, \mathrm{d}x\\
&\leq s\int_0^1 H(t,x)\, \mathrm{d}x,
\end{split}
\]
where we have used that $\psi(t,x)$ is increasing in the variable $t$, that both $\phi^1(x)-x$ and $H(t,x)$ are non-negative for all $t$ and $x$, and that assumption \eqref{ass} holds.

Plugging this estimate back in \eqref{unobis} and using \eqref{prima} we deduce the bound
\[
\int_0^1 H(t,x)\, \mathrm{d}x \leq  \int_0^1 \partial_t \phi(t,x) \, \mathrm{d}x + s\int_0^1 H(t,x)\, \mathrm{d}x. 
\]
which is equivalent to
\[
\int_0^1 H(t,x)\, \mathrm{d}x  \leq  \frac{1}{1-s} \int_0^1 \partial_t \phi(t,x) \, \mathrm{d}x
\]
since $s<1$. Combining this bound with the Jensen inequality applied to the convex function $r\mapsto 1/r$, $r\in (0,+\infty)$, we obtain
\[
\left(\int_0^1 \frac{\mathrm{d}x}{H(t,x)} \right)^{-1}\leq \int_0^1 H(t,x)\, \mathrm{d}x \leq\frac{1}{1-s} \int_0^1 \partial_t \phi(t,x) \, \mathrm{d}x.
\]
Integrating this last inequality in $t$ produces the desired upper bound
\[
\mathrm{length}_V \bigl( \{\phi^t\}_{t\in [0,1]} \bigr) = \int_0^1 \left(\int_0^1 \frac{\mathrm{d}x}{H(t,x)} \right)^{-1} \, \mathrm{d}t \leq \frac{1}{1-s} \int_0^1 \bigl( \phi^1(x) - x \bigr) \, \mathrm{d}x .\qedhere
\]
There remains to prove statement (ii). Given $\phi\in \mathrm{Diff}_1(\R)$ such that
\[
\phi(x) - x > \frac{1}{k} \qquad \forall x\in \R,
\]
we fix $\epsilon>0$ such that
\[
\phi(x) > \psi(x):= x + \frac{1}{k} + \epsilon   \qquad \forall x\in \R.
\]
The diffeomorphism $\psi$ belongs to the image of the lift 
\[
\tilde{\jmath}_k : \widetilde{\mathrm{Sp}}(2) \rightarrow \mathrm{Diff}_1(\R)
\]
of the homomorphism
\[
j_k: \mathrm{PSp}_k(2) \rightarrow \mathrm{Diff}_0(\T).
\]
More precisely, it is the image by $\tilde{\jmath}_k$ of the element $w\in \widetilde{\mathrm{Sp}}(2)$ which is given by the homotopy class of the path
\[
\{e^{\pi(1+k\epsilon)tJ_0} \}_{t\in [0,1]} \subset \mathrm{Sp}(2).
\] 
Since $\pi(1+k\epsilon)> \pi$, $w\geq \mathrm{id}$ does not belong to the closure of the positively elliptic region $\widetilde{\mathrm{Sp}}_{\mathrm{ell}}^+(2)$. By Proposition \ref{lordisSp2} (ii), there are timelike curves in $ \widetilde{\mathrm{Sp}}(2)$ from the identity to $w$ having arbitrarily large $\mathrm{length}_G$. By Proposition \ref{sottogruppi}, the image of these curves by $\tilde{\jmath}_k$ are positive paths in $\mathrm{Diff}_1(\R)$ from $\mathrm{id}$ to $\psi$ with arbitrarily large $\mathrm{length}_V$. By Remark \ref{poltrig}, these curves are generated by Hamiltonians which belong to $\mathcal{P}_k$. This concludes the proof of Theorem \ref{quantum}.

\begin{ex}
\label{leonid}
{\rm The above proof of statement (ii) shows how to use the Lorentz geometry of $\mathrm{Sp}(2)$ in order to construct arbitrarily long positive paths in $\mathrm{Diff}_1(\R)$ starting at the identity and ending at some given translation. Here is another more direct construction of an arbitrarily long positive path starting at $\mathrm{id}$ and remaining smaller than the translation $x\mapsto x+1+\epsilon$, for some fixed $\epsilon>0$. Choose a 1-periodic smooth function $h: \R \rightarrow \R$ with minimum $\epsilon$ and constant value $c$ on a subset of $[0,1]$ of measure $1-\delta$. By choosing $c$ large and $\delta$ small with respect to $\epsilon$, we can make the quantity
\[
V\bigl(h {\textstyle \frac{\partial}{\partial x}} \bigr) = \left( \int_0^1 \frac{\mathrm{d}x}{h(x)} \right)^{-1}
\]
as large as we wish. Let $\{\phi^t\}_{t\in [0,1]}$ be the positive path in $\mathrm{Diff}_1(\R)$ which is generated by the non-autonomous Hamiltonian $H_t(x):= h(x-\epsilon t)$. Its length is the arbitrarily large quantity
\[
\mathrm{length}_V\bigl( \{\phi^t\}_{t\in [0,1]} \bigr) = \int_0^1 V\bigl( H_t {\textstyle \frac{\partial}{\partial x}} \bigr)\, \mathrm{d}t = V\bigl(h {\textstyle \frac{\partial}{\partial x}} \bigr),
\]
but we claim that $\phi^1(x) \leq x + 1 + \epsilon$ for every $x\in \R$. Indeed, if $x:[0,1] \rightarrow \R$ solves the ODE $x'=H_t(x)$, then $y(t):= x(t) - \epsilon t$ solves the ODE $y' = h(y) - \epsilon$. The fact that the 1-periodic function $h-\epsilon$ has zeroes implies that $y(1)\leq y(0)+1=x(0)+1$, and hence $x(1) = y(1)+\epsilon \leq x(0)+ 1 + \epsilon$, as claimed. By replacing the function $h$ with $x\mapsto h(kx)$, we obtain arbitrarily long positive paths starting at $\mathrm{id}$ and remaining smaller than the translation $x\mapsto x+\frac{1}{k}+\epsilon$. 
}
\end{ex}

\section{Proof of Theorem \ref{lengthboundn}}
\label{lengthboundnnsec}

In this section, we prove Theorem \ref{lengthboundn} from the Introduction. Let $\phi\in \widetilde{\mathrm{Cont}}_0(\R \mathrm{P}^{2n-1},\xi_{\mathrm{st}})$ be such that
\begin{equation}
\label{ipotesi2}
\nu(\phi) \leq \frac{1}{2},
\end{equation}
with the strict inequality in the case $n=1$, and let $\{\phi_t\}_{t\in [0,1]}$ be the path generated by a Hamiltonian $H\in \mathcal{H}_c$ such that $\phi_0=\mathrm{id}$ and $\phi_1=\phi$. Let 
\[
\widetilde{Q}(t,z) = \frac{1}{2} S(t)z\cdot z
\]
be a smooth path of positive definite quadratic forms on $\R^{2n}$ such that the corresponding path $Q$ of positive quadratic Hamiltonians on $\R \mathrm{P}^{2n-1}$ satisfies
\begin{equation}
\label{pinch}
Q \leq H \leq c \, Q.
\end{equation}
Denote by $w:[0,1] \rightarrow \widetilde{\mathrm{Sp}}(2n)$ the timelike curve induced by the positive definite linear Hamiltonian system
\[
W'(t) = J_0 S(t) W(t), \qquad W(0) = \mathrm{id}.
\]
The induced path $\tilde{\jmath}\circ w$ in $\widetilde{\mathrm{Cont}}_0(\R \mathrm{P}^{2n-1},\xi_{\mathrm{st}})$ is generated by the Hamiltonian $Q$, so the first inequality in (\ref{pinch}) implies
\begin{equation}
\label{1s}
\mu(w(1)) = \nu(\tilde{\jmath}(w(1))) \leq \nu(\phi) \leq \frac{1}{2},
\end{equation}
with strict inequality in the case $n=1$. This upper bound on $\mu(w(1))$ implies that $w(1)$ belongs to the positively elliptic region $\widetilde{\mathrm{Sp}}^+_{\mathrm{ell}}(2n)$. Indeed, the non-decreasing function $t\mapsto \mu(w(t))$ has the value $0$ at $t=0$, and for $t>0$ and as long as $w(t)$ remains in $\widetilde{\mathrm{Sp}}^+_{\mathrm{ell}}(2n)$ has the value
\[
\mu(w(t)) = \frac{1}{2\pi} \sum_{j=1}^n \theta_j(t)
\]
where $\theta_j(t)\in (0,\pi)$ and $e^{i\theta_j(t)}$ are the Krein-positive eigenvalues of $\pi(w(t))$. The timelike curve $w$ can exit from $\widetilde{\mathrm{Sp}}^+_{\mathrm{ell}}(2n)$ only at some $t^*>0$ for which some $\theta_j$ takes the value $\pi$ and, in the case $n>1$, for such $t^*$ we would have $\mu(w(t^*))> \frac{1}{2}$, contradicting (\ref{1s}). In the case $n=1$, we would get $w(t^*)=\frac{1}{2}$, which contradicts the strict inequality in (\ref{1s}). Therefore, $w(1)$ belongs to  $\widetilde{\mathrm{Sp}}^+_{\mathrm{ell}}(2n)$ , as claimed.

By Theorem \ref{Lordist} (i), the Lorentz--Finsler length of $w$ has the upper bound
\[
\mathrm{length}_G(w) \leq \frac{2\pi}{n} \mu(w(1)) \leq \frac{2\pi}{n} \nu(\phi),
\]
where we have used (\ref{1s}) again. Combining this with the identity 
\[
j^* V = \frac{2^{\frac{1}{n}}}{2\pi} G
\]
of Proposition \ref{sottogruppi}, we obtain the bound
\[
\mathrm{length}_V(\tilde{\jmath}(w)) = \frac{2^{\frac{1}{n}}}{2\pi}  \mathrm{length}_G(w) \leq  \frac{2^{\frac{1}{n}}}{n}  \nu(\phi).
\]
Then the second inequality in (\ref{pinch}), together with the 1-homogeneity of the Lorentz--Finsler metric $V$, implies
\[
\mathrm{length}_V\bigl( \{\phi_t\}_{t\in [0,1]} \bigr) \leq c \, \mathrm{length}_V(\tilde{\jmath}\circ w) \leq  c\frac{2^{\frac{1}{n}}}{n}  \nu(\phi),
\]
concluding the proof of Theorem \ref{lengthboundn}.

\begin{rem}
\label{optimal}
{\rm The bound $\frac{1}{2}$ in assumption (\ref{ipotesi2}) is optimal. Indeed, let $s>\frac{1}{2}$ in the case $n>1$ and $s=\frac{1}{2}$ in the case $n=1$. In this case, we can find an element $w_1\in \widetilde{\mathrm{Sp}}(2n)$ with $\mu(w_1)=s$ which is the end-point of a timelike curve starting at the identity and which is not in the closure of $\widetilde{\mathrm{Sp}}^+_{\mathrm{ell}}(2n)$: in the case $n>1$, we consider a positive rotation of an angle $2\pi s - \epsilon> \pi$ in one symplectic plane and positive rotations of angle $\frac{\epsilon}{n-1}$ in the remaining $n-1$ symplectic planes; in the case $n=1$, we consider an element $w_1$ such that $\mu(w_1)=\frac{1}{2}$ and $\pi(w_1)$ is hyperbolic (necessarily with negative eigenvalues). By Theorem \ref{Lordist} (ii), there are arbitrarily long timelike curves $w$ from $\mathrm{id}$ to $w_1$. Setting $\phi=\tilde{\jmath}(w_1)$, we deduce that there are arbitrarily long positive paths $\tilde{\jmath}\circ w$ which connect the identity to $\phi$ in $\widetilde{\mathrm{Cont}}(\R \mathrm{P}^{2n-1},\xi_{\mathrm{st}})$ and are generated by positive quadratic Hamiltonians, i.e., by Hamiltonians in $\mathcal{H}_1$.
}
\end{rem}

\begin{rem}
\label{alternative}
{\rm
Assume that the answer to Question \ref{converse} is positive: if $\tilde{\jmath}(w_0) \leq \tilde{\jmath}(w_1)$ then $w_0\leq w_1$. Then in the above proof we can replace assumption (\ref{ipotesi2}) by the assumption $\phi\leq e$, where $e\in \widetilde{\mathrm{Cont}}_0(\R \mathrm{P}^{2n-1},\xi_{\mathrm{st}})$ is generated by the constant Hamiltonian $H=1$, and get the same conclusion. Indeed, let $w$ be the timelike curve in $\widetilde{\mathrm{Sp}}(2n)$ as in the above proof. The first inequality in (\ref{pinch}) implies that 
\[
\tilde{\jmath}(w(1))\leq \phi \leq e = \tilde{\jmath}(v),
\]
where $v$ is the element in $\widetilde{\mathrm{Sp}}(2n)$ corresponding to the homotopy class of the path $\{e^{\pi t J_0}\}_{t\in [0,1]}$. Then the fact that Question \ref{converse} is assumed to have a positive answer implies that $w(1)\leq v$. Since we also have $w(1)\geq \mathrm{id}$, $w(1)$ belongs to the closure of $\widetilde{\mathrm{Sp}}^+_{\mathrm{ell}}(2n)$ and the proof proceeds as before.
}
\end{rem}

We conclude this section by giving a short proof of the fact that Question \ref{converse} has a positive answer for $n=1$.

\begin{prop}
\label{converse1}
Consider the lift $\tilde{\jmath} : \widetilde{\mathrm{Sp}}(2) \rightarrow \mathrm{Diff}_1(\R)$ of the homomorphism $j: \mathrm{PSp}(2) \rightarrow \mathrm{Diff}_0(\T)$. Then $w_0 \leq w_1$ in $\widetilde{\mathrm{Sp}}(2)$ if and only if $\tilde{\jmath}(w_0) \leq  \tilde{\jmath}(w_1)$ in $\mathrm{Diff}_1(\R)$.
\end{prop}

\begin{proof}
The non-trivial implication we need to prove here is: $\tilde{\jmath}(w_0) \leq  \tilde{\jmath}(w_1)$ $\Rightarrow$ $w_0 \leq w_1$. By bi-invariance, we may assume $w_0=\mathrm{id}$. Since $\tilde{\jmath}(w_1)\geq \mathrm{id}$, we have
\[
\mu(w_1) = \frac{1}{2} \rho(\tilde{\jmath}(w_1)) \geq \frac{1}{2} \rho(\mathrm{id}) = 0.
\]
Assume by contradiction that $w_1\geq \mathrm{id}$ does not hold.  Then the non-negativity of $\mu(w_1)$ implies that $\mu(w_1)$ is actually zero and $w_1\in \widetilde{\mathrm{Sp}}(2)$ is either hyperbolic or parabolic and, up to conjugacy, given by the homotopy class of  
\[
\mbox{either} \qquad \left\{ \Bigr( \begin{array}{cc} e^{at} & 0 \\ 0 & e^{-at} \end{array} \Bigr) \right\}_{t\in [0,1]} \qquad  \mbox{or} \qquad \left\{ \Bigr( \begin{array}{cc} 1 & at \\ 0 & 1 \end{array} \Bigr) \right\}_{t\in [0,1]},
\]
for some $a>0$. See Section \ref{Sp(2)sec}, in particular Figure \ref{AdS} and identity (\ref{parametrization}). The image $\tilde{\jmath}(w_1)$ of the above hyperbolic element is a diffeomorphism $\varphi\in \mathrm{Diff}_1(\R)$ with fixed points at $\frac{1}{2} \Z$ and such that $\varphi(x)<x$ for every $x\in (0,\frac{1}{2}) + \Z$. The image $\tilde{\jmath}(w_1)$ of the above parabolic element is a diffeomorphism $\varphi\in \mathrm{Diff}_1(\R)$ with fixed points at $\Z$ and such that $\varphi(x)<x$ for every $x\in \R\setminus \Z$. In both cases, it is not true that $\varphi = \tilde{\jmath}(w_1) \geq \mathrm{id}$, which is the desired contradiction.
\end{proof}

\section{Proofs of the results of Section N}
\label{Mproofs}

The aim of this section is to prove the results stated in  Section \ref{secN} of the Introduction. We start by showing the strong concavity of the function $\mathcal{G}: \mathcal{H}^+(\overline{\Omega}) \rightarrow \R$.

\begin{prop}
\label{Mconcave}
For every $H\in \mathcal{H}^+(\overline{\Omega})$ and $K\in  \mathcal{H}(\overline{\Omega})$ we have 
\begin{equation}
\label{iinnee}
\mathrm{d}^2 \mathcal{G}(H)\cdot (K,K) \leq 0,
\end{equation}
with equality if and only if $K$ belongs to the line $\R H$.
\end{prop}

\begin{proof}
The inequality (\ref{iinnee}) is an immediate consequence of Lemma \ref{l:hessdet}. By the same lemma, the equality holds in (\ref{iinnee}) if and only if
\[
\nabla^2 K = u \,\nabla^2 H
\]
for some function $u:\overline{\Omega} \rightarrow \R$. From the fact that $\nabla^2 H$ is everywhere invertible, we deduce that $u$ is smooth. We rewrite the above identity using partial derivatives as
\[
\partial_i \nabla K = u \,\partial_i \nabla H \qquad \forall i=1,\dots,2n.
\]
By differentiating once more we obtain
\[
\partial_j \partial_i \nabla K = \partial_j u \,\partial_i \nabla H + u \, \partial_j \partial_i \nabla H.
\]
By subtracting from the above identity the analogous one which is obtained by exchanging $i$ and $j$ we find, by the symmetry of third derivatives of $H$ and $K$,
\[
\partial_j u \, \partial_i \nabla H - \partial_i u \, \partial_j \nabla H = 0.
\]
For every $i\neq j$ the vectors $\partial_i \nabla H$ and $\partial_j \nabla H$ are linearly independent, because the Hessian matrix of $H$ is invertible. Therefore, $\partial_i u = \partial_j u = 0$ for all $i\neq j$. We conclude that all partial derivatives of $u$ vanish, so $u$ is a constant function and
\[
\nabla^2 K = c\, \nabla^2 H
\]
for some real number $c$. It follows that
\begin{equation}
\label{intermedio}
K=c H + v
\end{equation}
for some affine function $v$. As a linear combination of functions representing elements of the vector space $\mathcal{H}(\overline{\Omega})$, $v$ is constant on each leaf of the characteristic foliation of $\partial \Omega$. Any affine function with this property is necessarily constant. Indeed, if $v$ is not constant then its level sets are given by a family of parallel hyperplanes. Choose a hyperplane $V$ in this family such that $V\cap \overline{\Omega}$ contains a point $z\in \partial \Omega$ and $\Omega$ is contained in one of the two open half-spaces determined by $V$. It follows that $V$ is tangent to $\partial \Omega$ at all points in $\overline{\Omega}\cap V$ and hence the characteristic foliation of $\partial \Omega$ is linear on the compact set $\partial \Omega\cap V$ (which possibly reduces just to the singleton $\{z\}$). In particular, the connected component containing $z$ of the intersection of the leaf through $z$ with $V$ is a closed segment. This implies that the leaf through $z$  cannot be fully contained in $V$, which contradicts the fact that $V$ is a level set of the function $v$.   

We conclude that $v$ is a constant function and hence (\ref{intermedio}) implies that $K$ agrees with $cH$ up to an additive constant. Therefore, $K$ belongs to the line $\R H$ in the quotient space $\mathcal{H}(\overline{\Omega})$, as we wished to prove.
\end{proof}

In the next proposition we prove the identity (\ref{GGversusG}) for the length of positive paths in $\mathrm{Symp}_0(\overline{\Omega})$.

\begin{prop}
\label{propGGversusG}
Let $\phi=\{\phi^t\}_{t\in [0,1]}$ be a positive path in $\mathrm{Symp}_0(\overline{\Omega})$. Then
\begin{equation}
\label{GGversusG2}
\mathrm{length}_{\mathcal{G}}(\phi) = \frac{1}{\mathrm{vol}(\Omega)} \int_{\Omega} \mathrm{length}_G(\{t\mapsto \mathrm{d} \phi^t(z)\})\, \mathrm{d}z.
\end{equation}
\end{prop} 

\begin{proof}
Let $H\in C^{\infty}([0,1]\times \overline{\Omega})$ be a uniformly convex Hamiltonian generating the path $\phi$. By differentiating
\[
\frac{\mathrm{d}}{\mathrm{d} t} \phi^t (z) = J_0 \nabla H_t(\phi^t(z)),
\]
we obtain the linearized equation
\[
\frac{\mathrm{d}}{\mathrm{d} t} \mathrm{d} \phi^t (z) = J_0 \nabla^2 H_t(\phi^t(z)) \mathrm{d} \phi^t (z),
\]
and hence for every $z\in \overline{\Omega}$ we have
\[
\mathrm{length}_G(\{t\mapsto \mathrm{d} \phi^t(z)\}) = \int_0^1 G( J_0 \nabla^2 H_t(\phi^t(z)))\, \di t = \int_0^1 \bigl( \det \nabla^2 H_t(\phi^t(z)) \bigr)^{\frac{1}{2n}}\, \mathrm{d}t.
\]
By integrating the above identity over $z\in \Omega$, switching the integrals and using the fact that $\phi^t$ is a volume-preserving diffeomorphism of $\Omega$, we find
\[
\begin{split}
\int_{\Omega}  \mathrm{length}_G&(\{t\mapsto \mathrm{d} \phi^t(z)\}) \, \mathrm{d}z = \int_{\Omega} \int_0^1  \bigl( \det \nabla^2 H_t(\phi^t(z)) \bigr)^{\frac{1}{2n}}\, \mathrm{d}t \, \mathrm{d}z \\ &= \int_0^1 \int_{\Omega} \bigl( \det \nabla^2 H_t(\phi^t(z)) \bigr)^{\frac{1}{2n}}\, \mathrm{d}z \,  \mathrm{d}t =\int_0^1 \int_{\Omega} \bigl( \det \nabla^2 H_t(z) \bigr)^{\frac{1}{2n}}\, \mathrm{d}z \,  \mathrm{d}t.
\end{split}
\]
On the other hand, by the definition of $\mathcal{G}$ we have
\[
\mathrm{length}_{\mathcal{G}}(\phi) = \int_0^1 \mathcal{G}(H_t) \, \di t = \int_0^1 \frac{1}{\mathrm{vol}(\Omega)} \int_{\Omega}  \bigl( \det \nabla^2 H_t(z) \bigr)^{\frac{1}{2n}}\, \mathrm{d}z \,  \mathrm{d}t,
\]
and (\ref{GGversusG2}) follows.
\end{proof}

We can now prove Theorems \ref{corti} and \ref{varthm} and Proposition \ref{smooth-uniqueness} from the Introduction.

\begin{proof}[Proof of Theorem \ref{corti}]
Let $z\in \overline{\Omega}$. The paths $\mathrm{d}\phi(z) = \{\mathrm{d}\phi^t(z)\}$ and $\mathrm{d}\psi(z) =\{\mathrm{d}\psi^t(z)\}$ are positive in $\mathrm{Sp}(2n)$ and start at the identity. The assumption that $\mathrm{d}\phi^t(z)$ never has the eigenvalue $-1$ guarantees that $\mathrm{d}\phi^t(z)$ belongs to the positive elliptic region $\mathrm{Sp}_{\mathrm{ell}}^+(2n)$ for every $t\in (0,1]$ and the path $\mathrm{d}\phi(z)$ has Conley--Zehnder index $n$ (for the definition of the Conley--Zehnder index and the properties which are needed here, see \cite{rs95}). Being a positive path, $\mathrm{d}\psi^t(z)$ belongs to  $\mathrm{Sp}_{\mathrm{ell}}^+(2n)$ for $t>0$ and small enough. Being homotopic to $\mathrm{d}\phi(z)$, the path $\mathrm{d}\psi(z)$ also has Conley--Zehnder index $n$. This implies that $\mathrm{d}\psi^t(z)$ remains in $\mathrm{Sp}_{\mathrm{ell}}^+(2n)$ for every $t\in (0,1]$. Indeed, the Conley--Zehnder index of any positive path starting at the identity which exits from the positively elliptic region and enters again in it is larger than $n$, as this path must meet the discriminant, i.e.\ the set of elements in $\mathrm{Sp}(2n)$ having the eigenvalue 1, and each intersection of a positive path with the discriminant contributes positively to the Conley--Zehnder index. 

We conclude that $\mathrm{d}\psi^t(z)$ belongs to $\mathrm{Sp}_{\mathrm{ell}}^+(2n)$ for every $t\in (0,1]$, and in particular does not have the eigenvalue $-1$ for every $t\in [0,1]$. By the former fact together with Theorem \ref{Lordist} (i) we obtain the bound
\[
\mathrm{length}_G(\mathrm{d}\psi(z)) \leq \frac{2\pi}{n} \mu\bigl([\mathrm{d}\psi(z)]\bigr) = \frac{2\pi}{n}\mu\bigl([\mathrm{d}\phi(z)]\bigr).
\]
Integration over $z\in \Omega$ and Proposition \ref{propGGversusG} yield the desired bound
\[
\mathrm{length}_{\mathcal{G}}(\psi) \leq \frac{2\pi}{n\, \mathrm{vol}(\Omega)} \int_{\Omega} \mu\bigl([\mathrm{d}\phi(z)]\bigr)\, \mathrm{d}z =   \frac{2\pi}{n} \mathcal{M}(\tilde\phi).\qedhere
\]
\end{proof}

\begin{proof}[Proof of Theorem \ref{varthm}]
Let $\phi$ be a positive path of symplectomorphisms extending the path $\psi$ and let $H\in C^{\infty}([0,1]\times \overline{\Omega})$ be a uniformly convex Hamiltonian generating $\phi$. By the H\"older inequality we obtain the upper bound
\[
\begin{split}
\mathrm{length}_{\mathcal{G}} (\phi) &= \frac{1}{\mathrm{vol}(\Omega)} \int_0^1 \int_{\Omega} \bigl( \det \nabla^2 H_t(z) \bigr)^{\frac{1}{2n}} \, \mathrm{d}z \, \mathrm{d} t \\ &\leq \mathrm{vol} (\Omega)^{-\frac{1}{2n}} \int_0^1 \left( \int_{\Omega} \det \nabla^2 H_t(z)\, \di z \right)^{\frac{1}{2n}}\, \mathrm{d}t,
\end{split}
\]
with equality if and only if 
\begin{equation}
\label{equa}
\det \nabla^2 H_t(z) = c(t) \qquad \forall t\in [0,1].
\end{equation} 
Since the map $\nabla H_t: \Omega \rightarrow \R^{2n}$ is a diffeomorphism onto its image, the change of variable formula gives us the identity
\[
 \int_{\Omega} \det \nabla^2 H_t(z)\, \di z  =   \mathrm{vol} \bigl( \nabla H_t(\Omega) \bigr),
 \]
and the bound
\[
\mathrm{length}_{\mathcal{G}} (\phi) \leq \mathrm{vol} (\Omega)^{-\frac{1}{2n}} \int_0^1 \mathrm{vol} \bigl( \nabla H_t(\Omega) \bigr)^{\frac{1}{2n}}\, \mathrm{d} t = \mathcal{V}(\psi)
\]
follows, with equality if and only if (\ref{equa}) holds.  
\end{proof}

\begin{proof}[Proof of Proposition \ref{smooth-uniqueness}]
The uniqueness of maximizers of the optimal extension problem (\ref{varprob}) follows from the strong concavity of $\mathcal{G}$. Indeed, let $\phi$ and $\phi'$ be two positive paths in $\mathrm{Symp}_0(\overline{\Omega})$ extending the same path of diffeomorphisms $\psi$ and having maximal length. If $H$ and $H'$ are uniformly convex smooth Hamiltonians generating $\phi$ and $\phi'$, then by Proposition \ref{Mconcave} we have for every $t\in [0,1]$
\begin{equation}
\label{primadiint}
\mathcal{G}\Bigr(\frac{1}{2} (H_t + H_t') \Bigl) \geq \frac{1}{2} \mathcal{G}(H_t) + \frac{1}{2} \mathcal{G}(H_t')
\end{equation}
with equality if and only if $H_t'$ is, up to an additive constant, a constant multiple of $H_t$. The Hamiltonian $\frac{1}{2} (H_t + H_t')$ generates a positive path $\phi''$ which still extends the path $\psi$ and integrating the above inequality in $t$ we obtain
\[
\mathrm{length}_{\mathcal{G}} (\phi'') \geq \frac{1}{2} \mathrm{length}_{\mathcal{G}} (\phi) + \frac{1}{2} \mathrm{length}_{\mathcal{G}} (\phi') =  \mathrm{length}_{\mathcal{G}} (\phi).
\]
The fact that $\phi$ is a maximizer implies that the above inequality is an equality, and we deduce that the inequality in (\ref{primadiint}) is an equality for every $t\in [0,1]$. Therefore,
\[
H_t'(z) = a(t) H_t(z) + b(t) \qquad \forall (t,z)\in [0,1]\times \overline{\Omega},
\]
for suitable numbers $a(t)>0$ and $b(t)$. Since $\phi$ and $\phi'$ extend the same path of diffeomorphisms of $\partial \Omega$, we have $\nabla H_t'(z) = \nabla H_t(z)$ for every $(t,z)\in [0,1]\times \partial \Omega$, so in the above identity we must have $a(t)=1$. Then $H'$ and $H$ differ by a function of $t$ and hence define the same positive path: $\phi'=\phi$.
\end{proof}

We conclude this section by describing a relaxation of the problem of maximizing the functional
\begin{equation}
\label{functional}
\mathcal{F}(H) = \int_{\Omega} \bigl(\det \nabla^2 H(z) \bigr)^{\frac{1}{2n}}\, \di z
\end{equation}
over the set of all uniformly convex smooth functions $H:\overline{\Omega} \rightarrow \R$ such that
\begin{equation}
\label{boundary}
H = K \mbox{ and } \nabla H = \nabla K \mbox{ on } \partial \Omega,
\end{equation}
for some fixed uniformly convex smooth functions $K:\overline{\Omega} \rightarrow \R$. This relaxation and the arguments for dealing with it are analogous to the study of the first variational problem for the affine area functional from \cite{tw05} and \cite{tw08}. 

Given a convex function $H:\Omega \rightarrow \R$, we denote by $\partial H(z_0)$ the set of its subgradients at $z_0\in \Omega$, i.e.\
\[
\partial H(z_0) := \{ p \in \R^{2n} \mid  H(z) \geq H(z_0) + p\cdot (z-z_0) \; \forall z\in \Omega \},
\]
which in the case of a point of differentiability reduces to the singleton $\{\nabla H(z_0)\}$. We consider the set of functions
\[
\mathcal{C}(\Omega,K) := \{ H : \overline{\Omega} \rightarrow \R \mid H \mbox{ is convex, } H|_{\partial \Omega} = K|_{\partial \Omega}, \; \partial H(\Omega) \subset \nabla K(\overline{\Omega}) \},
\]
 which is readily seen to be compact with respect to the $C^0(\overline{\Omega})$-topology. By a theorem of Aleksandrov, convex functions are almost everywhere twice differentiable and hence the functional $\mathcal{F}$ extends to a functional $\mathcal{F}: \mathcal{C}(\Omega,K) \rightarrow [0,+\infty]$. Actually, the measure
\[
\bigl(\det \nabla^2 H(z) \bigr)^{\frac{1}{2n}}\, \di z
\]
can be shown to be the absolutely continuous part of the \textit{Monge--Amp\`ere measure} of $H$, which associates to any Borel subset $E\subset \Omega$ the Lebesgue measure of $\partial H(E)$, see \cite[Lemma 2.3]{tw08}. The functional $\mathcal{F}$ is upper semicontinuous on $\mathcal{C}(\Omega,K)$ with respect to the $C^0(\overline{\Omega})$-topology. This is proven in \cite[Lemma 6.4]{tw08}, building on the weak continuity of the Monge--Amp\`ere measure as a function of $H$ (the proof of  \cite[Lemma 6.4]{tw08} is for a different exponent, but the modification for the exponent $\frac{1}{2n}$ is straightforward). We conclude that the functional $\mathcal{F}$ is finite and has maximizers in $\mathcal{C}(\Omega,K)$.

It can also be proven that the Monge--Amp\`ere measure of any maximizer is absolutely continuous with respect to the Lebesgue measure, by a straightforward modification of the proof of \cite[Lemma 6.5]{tw08}. Further regularity results for maximizers of the affine area functional, under the assumption that these maximizers are uniformly convex, are proven in \cite[Theorem 6.6]{tw08}, but we do not know whether analogous results hold also for the functional $\mathcal{F}$.

Maximizers  in $\mathcal{C}(\Omega,K)$ of the affine area functional are unique, see \cite[Theorem 6.5]{tw08}, but we do not know whether uniqueness also holds for maximizers of $\mathcal{F}$. This is due to the fact that, unlike the integrand of the affine area functional, the integrand of $\mathcal{F}$ is not uniformly concave, being linear in the radial direction. As a consequence, if we assume that $H,H'\in \mathcal{C}(\Omega,K)$ are maximizers of $\mathcal{F}$, from Lemma \ref{l:hessdet} we deduce that
\[
u(z) \nabla^2 H(z) + v(z) \nabla^2 H'(z) = 0 \mbox{ for a.e.\ } z\in \Omega,
\]
for suitable real functions $u$ and $v$. If $H$ is uniformly convex and both $H$ and $H'$ are three times differentiable, then the argument of the proof of Proposition \ref{Mconcave} shows that the functions $u$ and $v$ must be constant, and by using the boundary conditions we conclude that $H=H'$. Without uniform convexity and higher differentiability assumptions, we do not know whether uniqueness holds.

\appendix

\renewcommand{\thesection}{\roman{section}}

\section{Bi-invariant Lorentz--Finsler metrics on Lie groups}
\label{liegroups}


The study of bi-invariant Finsler metrics on Lie groups was initiated in the seventies by Grove, Karcher and Ruh, see \cite{gkr74}. In this Appendix, we establish a few facts about bi-invariant Lorentz--Finsler metrics.

Let $G$ denote a Lie group and denote by $\mathfrak g$ the Lie algebra of $G$, namely the tangent space of $G$ at the identity $e$. Given $X\in \mathfrak g$ and $w\in G$, we shall use the notation
\begin{equation}\label{e:Xw}
Xw:= \mathrm{d} R_w(e)\cdot X \in T_w G,
\end{equation}
where $R_w:G \rightarrow G$ denotes the right multiplication by $w$. We call $\mathcal X$ the right-invariant vector field on $G$ extending $X$. In other words $\mathcal X(w)=Xw$ for all $w\in G$. 
We define the Lie bracket of two elements $X,Y$ of $\mathfrak g$ by
\begin{equation}\label{e:brabra}
	[X,Y]=[\mathcal X,\mathcal Y](e),
\end{equation}
where $\mathcal X$ and $\mathcal Y$ are the right-invariant extensions of $X$ and $Y$. The bracket on the right-hand side is the bracket of two vector fields on the manifold $G$, which is given by the convention
\begin{equation}\label{e:oppositesign}
	[\mathcal X,\mathcal Y]=-\mathcal L_{\mathcal X}\mathcal Y,
\end{equation}
where $\mathcal L$ denotes the Lie derivative of vector fields, see the discussion in Section \ref{secG} of the Introduction.

Let $(K,F)$ be a Lorentz--Finsler structure on a Lie group $G$, in the sense of Definition \ref{deflorfin} in the Introduction, and assume that $(K,F)$ is bi-invariant. In this Appendix, we actually do not need $\overline{K} \cap \overline{-K}$ to coincide with the zero-section and $F$ to extend continuously to the closure of $K$ by the zero extension on the boundary, because we are going to consider only timelike curves.

The pair $(K,F)$ is uniquely determined by its restriction to $\mathfrak{g}$, namely by the open convex cone $\kappa:=K\cap \mathfrak g$ and by the smooth function $f:=F|_{\kappa}$. The cone $\kappa$ and the function $f$ are invariant under the adjoint action. This implies that
\begin{equation}
\label{B0}
\mathrm{d}f(X)\cdot[X,Y]=0,\qquad\forall\,X\in \kappa, \; \forall\, Y\in\mathfrak g.
\end{equation}

\begin{rem}
\label{injective} 
{\rm Let $c>0$. By our assumptions on $f$, the restriction of the map $\mathrm{d}f: \kappa \rightarrow \mathfrak g^*$ to the subset $\{f=c\}$ is injective. Indeed, if $X$ and $Y$ are distinct elements in $\{f=c\}$, then the segment joining them is contained in the convex set $\mathcal{\kappa}$ and, by the positive 1-homogeneity of $f$ and the condition $f(X)=f(Y)$, the vector $Y-X$ is not collinear to $X+t (Y-X)$ for any $t\in \R$. Therefore, the smooth function
\[
g: [0,1] \rightarrow \R, \qquad g(t) = f(X+t(Y-X))
\]
satisfies
\[
g''(t) = \mathrm{d}^2 f (X+t (Y-X)) \cdot (Y-X,Y-X) < 0,
\]
and hence 
\[
\mathrm{d}f(X)\cdot (Y-X) = g'(0) > g'(1) = \mathrm{d} f(Y)\cdot (Y-X).
\]
In particular, $\mathrm{d} f(X)\neq \mathrm{d} f(Y)$, proving that the restriction of $\mathrm{d}f$ to $\{f=c\}$ is indeed 
injective.}
\end{rem}

 To any continuously differentiable curve $w: [0,1] \rightarrow G$ we can associate the continuous curve $X: [0,1] \rightarrow \mathfrak g$ that  is defined by $w'=Xw$, see \eqref{e:Xw}. Then $w$ is timelike if and only if $X$ takes values in $\kappa$. The Lorentz--Finsler length of the curve $w$ has then the form
\[
\mathrm{length}_F(w):=\int_0^1 F(w'(t))\,\mathrm{d}t=\int_0^1f(X(t))\,\mathrm{d}t.
\]
In order to compute the first variation of the functional $\mathrm{length}_F$, we need the following lemma, whose proof is adapted from \cite[Proposition I.1.1]{ban78}. 
\begin{lem}\label{l:inter2}
Let $w_s:[0,1]\to G$ be a smooth one-parameter family of continuously differentiable paths with $s\in(-\epsilon,\epsilon)$. Define paths $X_s,Y_s$ in $\mathfrak{g}$  by $\partial_tw_s=X_sw_s$ and $\partial_sw_s=Y_sw_s$. Then,
	\[
	\partial_s X_s=\partial_tY_s+ [Y_s,X_s].
	\]
\end{lem}

\begin{proof}
Define $\tilde G:=G\times (-\epsilon,\epsilon)\times[0,1]$ and consider the standard projection
\[
p\colon \tilde G\to(-\epsilon,\epsilon)\times[0,1],\qquad p(w,s,t)=(s,t).
\]
The path $w_s$ yields naturally a section of $p$ 
\[
\tilde w\colon (-\epsilon,\epsilon)\times[0,1]\to \tilde G,\qquad \tilde w(s,t)=(w_s(t),s,t),
\]
which is an embedded surface in $\tilde G$ transverse to $p$. For every $(s,t)\in (-\epsilon,\epsilon)\times[0,1]$, let $\mathcal X_{(s,t)}$ and $\mathcal Y_{(s,t)}$ be the extensions of $X_s(t)$ and $Y_s(t)$ as right-invariant vectors fields on $G$. Thus, the lifted vector fields on $\tilde G$
\[
\tilde{\mathcal X}:=\mathcal X+\partial_t \qquad \mbox{and} \qquad  \tilde{\mathcal Y}:=\mathcal Y+\partial_s
\]
are tangent to $\tilde{w}$. It follows that the Lie bracket of vector fields $[\tilde {\mathcal X},\tilde{\mathcal Y}]$  is also tangent to $\tilde{w}$ and hence transverse to $p$. On the other hand, 
\[
p_*[\tilde{\mathcal X},\tilde{\mathcal Y}]=[p_*\tilde{\mathcal X},p_*\tilde{\mathcal Y}]=[\partial_t,\partial_s]=0,
\]
namely $[\tilde{\mathcal X},\tilde{\mathcal Y}]$ is tangent to the fibers of $p$. We conclude that
\[
[\tilde{\mathcal X},\tilde{\mathcal Y}]=0\,.
\]
On the other hand, using the convention \eqref{e:oppositesign} we get
\[
[\tilde{\mathcal X},\tilde{\mathcal Y}]=[\mathcal X+\partial_t,\mathcal Y+\partial_s]=[\mathcal X,\mathcal Y]-\partial_t\mathcal Y+\partial_s\mathcal X\,. 
\]
Therefore, for all $(s,t)\in (-\epsilon,\epsilon)\times[0,1]$ we have the following equality of vector fields on $G$
\[
\partial_s\mathcal X_{(s,t)}=\partial_t\mathcal Y_{(s,t)}+[\mathcal Y_{(s,t)},\mathcal X_{(s,t)}]\,.
\]
Evaluating this equality at the identity $e\in G$ and using \eqref{e:brabra}, we arrive at the desired formula. 
\end{proof}

The differential of the functional $\mathrm{length}_F$ at some timelike curve $w: [0,1] \rightarrow G$ is defined on the space sections of the vector bundle $w^*(TG)$. Using the group structure, we identify these sections with curves in the Lie algebra $\mathfrak{g}$: a curve $Y: [0,1] \rightarrow \mathfrak{g}$ defines the section $Yw$. The differential of $\mathrm{length}_F$ at $w$ can hence be seen as a linear functional on curves in $\mathfrak{g}$ and we use the notation
\[
\mathrm{d}\, \mathrm{length}_F(w)\cdot Y := \frac{\mathrm{d}}{\mathrm{d}s} \Bigr|_{s=0} \mathrm{length}_F(e^{sY} w),
\]
where $Y: [0,1] \rightarrow \mathfrak{g}$.

\begin{prop}[First variation]
	\label{lemvar1b}
	Let $w: [0,1] \rightarrow G$ be a timelike curve with tangent vector field $w'=Xw$, where $X:[0,1]\to\kappa$. Then the first variation of $\mathrm{length}_F$ at $w$ in the direction $Y: [0,1] \rightarrow  \mathfrak g$ has the form
	\[
	\mathrm{d}\, \mathrm{length}_F(w)\cdot Y = \int_0^1 \mathrm{d} f(X(t))\cdot Y'(t) \, \mathrm{d}t.
	\]
	In particular, $\mathrm{d}\, \mathrm{length}_F(w) \cdot Y=0$ for every smooth $Y$ with compact support in $(0,1)$ if and only if $w$ is the reparametrization of an autonomous path, i.e.,
	\[
	w(t) = e^{\tau(t) A}w(0)
	\]
	for some $A\in \kappa$ and for some continuously differentiable function $\tau$ with $\tau'>0$ and $\tau(0)=0$.
\end{prop}

\begin{proof}
Set 
\[
w_s(t):= e^{sY(t)} w(t).
\]
Then
\[
\partial_s w_s = Y w_s.
\]
We define the paths $X_s$ by
\[
\partial_t w_s = X_s w_s, 
\]
so that $X_0 = X$.
Using Lemma \ref{l:inter2} and the invariance of $f$ by the adjoint action, see (\ref{B0}), we compute
	\begin{equation}
		\label{var1b}
		\begin{split}
			\frac{\mathrm{d}}{\mathrm{d} s} \mathrm{length}_F(w_s) &= \frac{\mathrm{d}}{\mathrm{d} s} \int_0^1 f(X_s)\, \mathrm{d}t =\int_0^1 \mathrm{d} f(X_s)\cdot\partial_sX_s\, \mathrm{d}t \\ &= \int_0^1 \Big(\mathrm{d} f(X_s)\cdot Y' +\mathrm{d} f(X_s)\cdot  [Y,X_s] \Big)\, \mathrm{d} t=\int_0^1 \mathrm{d} f(X_s)\cdot Y' \, \mathrm{d}t.
		\end{split}
	\end{equation}
By evaluating at $s=0$ we find the desired formula.

Now suppose that the first variation of $\mathrm{length}_F$ at $w$ in the direction $Y$ vanishes for every smooth $Y$ with compact support in $(0,1)$. By a time reparametrization, which leaves the critical set of $\mathrm{length}_F$ invariant, we assume that $f(X(t))=c>0$ does not depend on $t$. The Du Bois-Reymond lemma tells us that $\mathrm{d} f(X(t))$ is constant in $t$. Thanks to Remark \ref{injective}, we conclude that $X(t)=X_0$ is constant. It follows that the path $t\mapsto e^{-tX_0}w(t)=w_0$ is constant and hence $w(t)=e^{tX}w_0$ as was to be shown.	
\end{proof}

We deduce that timelike geodesics are precisely the curves of the form
\[
w(t) = e^{tX} w_0,
\]
where $X\in \kappa$ and $w_0\in G$. We can now compute the second variation of $\mathrm{length}_F$ at a timelike geodesic, again seen as symmetric bilinear form on the space of curves in $\mathfrak{g}$.

\begin{prop}[Second variation]
	\label{lemvar2b}
	Let $w: [0,1] \rightarrow G$, $w(t)=e^{tX}w_0$, with $X\in\kappa$ and $w_0\in G$, be a timelike geodesic. Then the second variation of $\mathrm{length}_F$ at $w$ is the symmetric bilinear form
\begin{equation}
\label{secvar}
\mathrm{d}^2\, \mathrm{length}_F(w) \cdot (Y_1,Y_2) = \int_0^1 \mathrm{d}^2f(X)\cdot(Y'_1+ [Y_1,X],Y'_2)\, \mathrm{d}t,
\end{equation}
for every pair of curves $Y_1,Y_2:[0,1] \rightarrow \mathfrak{g}$ vanishing at $t=0,1$.
\end{prop}

\begin{proof}
Consider a curve $Y: [0,1] \rightarrow \mathfrak{g}$ vanishing at $t=0,1$, and let $w_s$, $X_s$ and $Y_s$ be as in the proof of Proposition \ref{lemvar1b}. Differentiating \eqref{var1b} with respect to $s$, we have by Lemma \ref{l:inter2}
\[
\begin{split}
\frac{\mathrm{d}^2}{\mathrm{d}s^2}\mathrm{length}_F(w_s)&=\int_0^1 \Big(\mathrm{d}^2f(X_s)\cdot(\partial_sX_s,Y') =\int_0^1 \Big(\mathrm{d}^2f(X_s)\cdot(Y'+ [Y,X_s],Y').
\end{split}
\]
Evaluating at $s=0$, we get the formula
\[
\mathrm{d}^2 \,\mathrm{length}_F(w) \cdot (Y,Y) = \int_0^1 \mathrm{d}^2f(X)\cdot(Y'+ [Y,X],Y'), \mathrm{d}t.
\]	
This is precisely the quadratic form that is induced by the bilinear form (\ref{secvar}). In order to conclude, we must show that the bilinear form (\ref{secvar}) is symmetric. Since the second differential of $f$ is symmetric, we just need to prove the symmetry of the bilinear form
\[
(Y_1,Y_2) \mapsto \int_0^1 \mathrm{d}^2 f(X)\cdot ([Y_1,X],Y_2')\, \di t.
\]
By differentiating \eqref{B0} we find
\[
\mathrm{d}^2 f(X)\cdot (Z,[Y,X]) + \mathrm{d} f(X)\cdot [Y,Z] = 0, \qquad \forall X\,\in \kappa, \; \forall\, Y,Z\in \mathfrak{g},
\]
and hence
\[
\begin{split}
 \int_0^1 \mathrm{d}^2 f(X)\cdot ( [Y_1,X],Y_2')\, \mathrm{d}t &= \int_0^1 \mathrm{d}f(X)\cdot  [ Y_2',Y_1] \, \mathrm{d}t \\ &=  \int_0^1 \mathrm{d} f(X) \cdot [ Y_1',Y_2] \, \mathrm{d}t =  \int_0^1 \mathrm{d}^2 f(X)\cdot ([Y_2,X],Y_1')\, \mathrm{d}t,
 \end{split}
\]
where the middle equality follows from Lemma \ref{lemvar1b} applied to the critical curve $X$, thanks to the fact 
that $[Y_2',Y_1]$ and $[Y_1',Y_2]$ differ by the derivative of the curve $[Y_2,Y_1]$, which vanishes at $0$ and $1$.
\end{proof}

\begin{rem} \label{extensiontoDiff}
{\rm
Propositions \ref{lemvar1b} and \ref{lemvar2b} hold also when $G$ is an infinite dimensional group of diffeomorphisms of a manifold $M$. In this case, $\mathfrak{g}$ is a subspace of the space of smooth vector fields on $M$ and the Lie bracket on $\mathfrak g$ is the Lie bracket of vector fields on $M$ using the sign convention \eqref{e:oppositesign}. Indeed, if $G$ is a group of diffeomorphisms, Lemma \ref{l:inter2} is proved by Banyaga in \cite[Proposition I.1.1]{ban78} (using the sign convention opposite to \eqref{e:oppositesign}). Moreover, the argument in Remark \ref{injective}, which is used in the proof of Proposition \ref{lemvar1b}, applies also when $\mathfrak{g}$ is infinite dimensional and $\mathfrak{g}^*$ is its algebraic dual space.}
\end{rem}

Let $w(t) = e^{tX} w_0$ be a timelike geodesic. Because of the invariance of $\mathrm{length}_F$ by time reparametrizations, $\mathrm{d}^2 \, \mathrm{length}_F(w)$ has an infinite dimensional kernel: indeed, the fact that $\R X$ is in the kernel of $\mathrm{d} ^2f(X)$ implies any curve of the form $Y(t)= u(t) X$, with $u$ a real function vanishing at $t=0,1$, is in the kernel of $\mathrm{d}^2 \mathrm{length}_F(w)$.

In order to get rid of this invariance by reparametrizations, let us consider the linear splitting
\[
\mathfrak{g} = \mathfrak{g}_X \oplus \R X \qquad \mbox{where} \quad \mathfrak g_X:=\ker \mathrm{d} f(X),
\]
and correspondingly
\[
H^1_0((0,1),\mathfrak{g}) = H^1_0((0,1),\mathfrak{g}_X) \oplus H^1_0((0,1),\R X),
\]
where $H^1_0((0,1),V)$ denotes the Sobolev space of  absolutely continuous curves in the vector space $V$ vanishing at the end-points and having square integrable derivative. The symmetric bilinear form $\mathrm{d}^2\, \mathrm{length}_F(w)$ is continuous on $H^1_0((0,1),\mathfrak{g})$ and has an infinite dimensional kernel containing the second space of the above splitting. By restricting it to the first space, 
we obtain the continuous symmetric bilinear form
\[
\begin{split}
\mathcal{H}_X&: H^1_0((0,1),\mathfrak{g}_X) \times H^1_0((0,1),\mathfrak{g}_X) \rightarrow \R, \\
\mathcal{H}_X (Y_1,Y_2) &:= \mathrm{d}^2\, \mathrm{length}_F (w)\cdot (Y_1,Y_2) = \int_0^1 \mathrm{d}^2f(X)\cdot(Y_1'+[Y_1,X],Y_2')\, \mathrm{d}t.
\end{split}
\]
Since $\mathrm{d}^2f(X)$ is negative definite on $\mathfrak{g}_X$, the space $H^1_0((0,1),\mathfrak{g}_X)$ admits the equivalent inner product
\[
(Y_1,Y_2)_{H^1_0} :=-\int_0^1\mathrm{d}^2f(X)\cdot (Y_1',Y_2')\, \mathrm{d}t
\]
and we denote by $H_X$ the linear selfadjoint operator on $H^1_0((0,1),\mathfrak{g}_X)$ representing $\mathcal H_X$ with respect to this inner product. 

Thanks to the compact embedding of $H^{1}_0$ into $H^{1/2}$, $H_X$ is a compact perturbation of minus the identity. Hence, the spectrum of $H_X$ consists of $-1$ and of a sequence of real eigenvalues of finite multiplicity converging to $-1$. Hence, only a finite number of eigenvalues are non-negative and each of these has finite multiplicity. Therefore, the symmetric bilinear form $\mathcal{H}_X$ has finite dimensional kernel and finite co-index. The eigenvectors corresponding to the eigenvalue zero are exactly the Jacobi field vanishing at the boundary of the interval $[0,1]$.

In general, the eigenvalues are precisely the real numbers $\lambda$ such that the symmetric bilinear form
\[
((H_X-\lambda I)Y_1,Y_2)_{H^1_0} =  \int_0^1 \mathrm{d}^2f(X)\cdot \bigl((1+\lambda) Y_1'+[Y_1,X],Y_2'\bigr)\, \mathrm{d} t
\]
has a non-trivial kernel. By a standard regularity argument, integration by parts and (\ref{B0}), we deduce that $Y_\lambda\in H^1_0((0,1),\mathfrak{g}_X)$ lies in the kernel of $H_X-\lambda I$ if and only if it is a smooth solution of the equation
\begin{equation}\label{e:eigenhess}
(1+\lambda)Y''_\lambda=[X,Y'_\lambda].
\end{equation}
The solutions of this equation satisfying $Y_{\lambda}(0)=0$ can be explicitly written as
\begin{equation}\label{e:jack}
Y_{\lambda}(t)=(1+\lambda)Y_Z(\tfrac{1}{1+\lambda}t),\qquad Y_Z(t):=\int_0^t\mathrm{Ad}(e^{\tau X})\cdot Z\,\mathrm{d}\tau,
\end{equation}
for some $Z\in\mathfrak g_X$. Requiring that $Y_\lambda(1)=0$ however implies that the condition $Z\in\mathfrak g_X$ is redundant since from \eqref{e:eigenhess} and the bi-invariance of $f$, see \eqref{B0}, the function $t\mapsto \mathrm{d} f(X)\cdot Y_\lambda(t)$ is affine and therefore it must vanish identically since it vanishes for $t=0,1$. Thus, the eigenspace with eigenvector $\lambda$ is isomorphic to the vector space
\[
V_\lambda:=\big\{Z\in\mathfrak g\ \big|\ Y_Z(\tfrac{1}{1+\lambda})=0\big\}.
\]
The fields $Y_Z$ correspond to the eigenvalue $0$ and yields therefore the Jacobi fields along $W$. As usual, we call $t\in[0,1]$ a conjugate instant along the timelike geodesic $W$ if the space
\[
V_0(t):=\{Z\in\mathfrak g\ |\ Y_Z(t)=0\}
\]
is non-zero and we call $m(t):=\dim V_0(t)$ the multiplicity of the conjugate instant. For $\lambda>0$, $V_\lambda=V_0(\tfrac{1}{1+\lambda})$ by definition, thus $V_\lambda$ is non-trivial if and only if $\tfrac{1}{1+\lambda}$ is a conjugate instant and $\dim V_\lambda=m(t)$. Since the function $\lambda\mapsto \tfrac{1}{1+\lambda}$ is a bijection between $(0,\infty)$ and $(0,1)$, we get the co-index formula
\begin{equation}\label{e:co-ind}
\mbox{co-ind} \, \mathcal{H}_X :=\sum_{\lambda>0}\dim V_\lambda=\sum_{t\in(0,1)}m(t).
\end{equation}

We summarize the above discussion into the following proposition.

\begin{prop}
\label{morse-thm}
Let $w(t)=e^{tX} w_0$, $X\in \kappa$, $w_0\in G$ be a timelike geodesic. Then the symmetric bilinear form $\mathrm{d}^2\, \mathrm{length}_F(w)$ on $H^1_0((0,1),\mathfrak{g})$ has finite co-index and an infinite dimensional kernel containing $H^1_0((0,1),\R X)$. The kernel of its restriction $\mathcal{H}_X$ to $H^1_0((0,1),\mathfrak{g}_X)$ is finite dimensional and coincides with the space of Jacobi vector fields along $w$, i.e., solutions $Y: [0,1] \rightarrow \mathfrak{g}$ of the equation
\[
Y'' = [X,Y']
\]
vanishing at $t=0$ and $t=1$. Moreover
\[
\mbox{\rm co-ind} \, \mathrm{d}^2\mathrm{length}_F (w) = \mbox{\rm co-ind} \, \mathcal{H}_X = \sum_{t^* \in(0,1)}m(t^*),
\]
where $m(t^*)$ denotes the dimension of the space of Jacobi vector fields $Y$ such that $Y(0)=Y(t^*)=0$.
\end{prop}

The {\it Morse co-index} of the timelike geodesic segment $w: [0,1] \rightarrow G$ is defined to be the co-index of the second differential of $\mathrm{length}_F$ at $w$:
\[
\mbox{co-ind} (w) := \mbox{co-ind} \, \mathrm{d}^2\, \mathrm{length}_F(w) = \mbox{ co-ind} \, \mathcal{H}_X.
\]

\section{Some facts about the Lie algebra of the symplectic group}
\label{linear-app}

We recall that $\mathrm{sp}^+(2n)$ is defined as the subset of the Lie algebra $\mathrm{sp}(2n)$ consisting of those endomorphisms $X$ for which the symmetric bilinear form $(u,v) \mapsto \omega_0(u,Xv)$ is positive definite. In this appendix we prove a characterization of the elements of  $\mathrm{sp}^+(2n)$ which is used extensively in this monograph. Here, $\kappa=-i\omega_0$ denotes the Krein form on $\C^{2n}$, see  Section \ref{secJ} in the Introduction, and a basis $u_1,\dots,u_n,v_1,\dots,v_n$ of $\C^{2n}$ is said to be $\kappa$-unitary if 
\[
\kappa(u_j,u_j) = 1 = - \kappa(v_j,v_j) \; \forall j, \quad \kappa(u_j,u_h) = \kappa(v_j,v_h) = 0 \; \forall j\neq h, \quad \kappa(u_j,v_h) = 0 \; \forall j,h.
\]
Moreover, such a basis is said to be real if $v_j = \overline{u}_j$ for every $j$.

\begin{prop}
\label{propA1}
Let $X$ be an endomorphism of $\R^{2n}$. Then the following facts are equivalent:
\begin{enumerate}[(i)]
\item $X$ belongs to $\mathrm{sp}^+(2n)$;
\item there exists a $\kappa$-unitary real basis $w_1,\dots,w_n,\overline{w}_1,\dots,\overline{w}_n$ of $\C^{2n}$ such that
\begin{equation}
\label{autoval1}
X w_j =i \theta_j w_j, \quad  X \overline{w}_j =-i \theta_j \overline{w}_j\qquad \forall j\in \{1,\dots,n\}.
\end{equation}
for some positive numbers $\theta_j$;
\item there is an $X$-invariant symplectic splitting of $\R^{2n}$ into $n$ symplectic planes, i.e.,
\[
\R^{2n} = \bigoplus_{j=1}^n V_j, \qquad \dim V_j=2, \quad \omega_0(u,v) = 0\quad  \forall u\in V_j, \; v\in V_h \mbox{ with } j\neq h,
\]
with respect to which $X$ has the form
\begin{equation}
\label{theform}
X = \bigoplus_{j=1}^n \theta_j J_j,
\end{equation}
where each $\theta_j$ is a positive number and each $J_j:V_j\to V_j$  is an $\omega_0$-compatible complex structure on the symplectic plane $V_j$. 
\end{enumerate}
\end{prop}

\begin{proof}
(i) $\Rightarrow$ (ii). Being an element of $\mathrm{sp}(2n)$, $X$ is real and $\kappa$-skew-Hermitian. The latter fact implies that $-iX$ is $\kappa$-Hermitian. Moreover, the fact that $X$ is in $\mathrm{sp}^+(2n)$ implies that 
\[
\kappa(-iXw,w) = -\omega_0(Xw,w) = \omega_0(w,Xw) > 0 \qquad \forall w\in \C^{2n}\setminus \{0\}.
\]
Therefore, the hypersurface
\[
\Sigma:= \{w\in \C^{2n} \mid \kappa(-iXw,w) = 1 \} 
\]
is compact. Let $z_1\in \Sigma$ be a maximizer for the restriction of the real function $z\mapsto \kappa(z,z)$ to $\Sigma$. Since this real function is somewhere positive on $\Sigma$, we have $\kappa(z_1,z_1) > 0$. By the Lagrange multipliers theorem, we have 
\[
-i X z_1 = \theta_1 z_1
\]
for some $\theta_1\in \R$. By taking the $\kappa$-product with $z_1$ we obtain
\[
1=\kappa(-i X z_1, z_1) = \theta_1 \kappa (z_1,z_1),
\]
so $\theta_1= 1/\kappa(z_1,z_1)$ is positive. Setting $w_1 := \sqrt{\theta_1} z_1$ we have
\[
X w_1 = i \theta_1 w_1, \qquad \kappa(w_1,w_1) = 1.
\]
Since $X$ is real, we also have
\[
X \overline{w}_1 = - i \theta_1 \overline{w}_1.
\]
From the identity
\[
\overline{\kappa(w,w')} = -\kappa(\overline{w},\overline{w}') \qquad \forall w,w'\in \C^{2n},
\]
we obtain
\[
\kappa(\overline{w}_1,\overline{w}_1) = -1
\]
and, together with the fact that $\kappa$ is Hermitian,
\[
\kappa(w_1,\overline{w}_1) = \kappa( \overline{\overline{w}}_1 , \overline{w}_1) = - \overline{\kappa(\overline{w}_1,w_1)}= - \kappa(w_1,\overline{w_1}) ,
\]
which implies
\[
\kappa(w_1,\overline{w}_1) =0.
\]
By considering the $\kappa$-orthogonal complement to the 2-dimensional complex subspace generated by the vectors $w_1$ and $\overline{w}_1$, on which $\kappa$ has signature $(n-1,n-1)$, we can iterate the above argument and produce the desired real unitary basis of eigenvectors of $Y$. 

\medskip

\noindent (ii) $\Rightarrow$ (iii). Denote by $V_j\subset \R^{2n}$ the plane that is obtained from intersecting the 2-dimensional conjugation invariant complex subspace spanned by $w_j$ and $\overline{w}_j$ with $\R^{2n}$. The fact that the real basis $w_1,\dots,w_n,\overline{w}_1,\dots,\overline{w}_n$ is $\kappa$-unitary implies that the $V_j$'s form a symplectic splitting of $\R^{2n}$. By (\ref{autoval1}), this splitting is $X$-invariant and $J_j:= \theta_j^{-1} X|_{V_j}$ is a complex structure on $V_j$. For every $u = \alpha w_j + \overline{\alpha} \overline{w}_j, v= \beta w_j + \overline{\beta} \overline{w}_j \in V_j$ we have
\[
\omega_0(u, J_j v) = \omega_0 (  \alpha w_j + \overline{\alpha} \overline{w}_j, \beta i w_j - \overline{\beta}  i \overline{w}_j) = \kappa (  \alpha w_j + \overline{\alpha} \overline{w}_j, \beta  w_j - \overline{\beta}   \overline{w}_j) = 2\, \re (\alpha \overline{\beta}),
\]
so the bilinear form $u,v\mapsto \omega_0(u,J_j v)$ is symmetric and positive definite on $V_j$. We conclude that $J_j$ is an $\omega_0$-compatible complex structure on $V_j$.

\medskip

\noindent (iii) $\Rightarrow$ (i). Writing any $u,v\in \R^{2n}$ as
\[
u = \sum_{j=1}^n u_j, \qquad v = \sum_{j=1}^n v_j,
\]
with $u_j,v_j\in V_j$ for every $j$, we have
\[
\omega_0(u,X v) = \sum_{j=1}^n \theta_j \, \omega_0(u_j,J_j v_j).
\]
The fact that each $\theta_j$ is positive and each bilinear form $(u,v)\mapsto \omega_0(u,J_j v)$ is symmetric and positive definite on $V_j$ implies that the bilinear form $(u,v)\mapsto \omega_0(u,X v)$ is symmetric and positive definite on $\R^{2n}$. Therefore, $X$ belongs to $\mathrm{sp}^+(2n)$.
\end{proof}

We conclude this section by stating the analogous characterization for the elements of $\mathrm{Sp}^+_{\mathrm{ell}}(2n)$, the open subset of $\mathrm{Sp}(2n)$ consisting of elliptic automorphisms all of whose eigenvalues are Krein-definite and such that the Krein-positive ones have positive imaginary part. 

\begin{prop} 
\label{propA2}
The set
\[
\mathrm{sp}^+_{\mathrm{ell}}(2n):= \{ X \in \mathrm{sp}^+(2n) \mid \sigma(X) \subset (-\pi,\pi)i \}
\]
is contractible and $\exp: \mathrm{sp}(2n) \rightarrow \mathrm{Sp}(2n)$ maps it diffeomorphically onto $\mathrm{Sp}^+_{\mathrm{ell}}(2n)$. For an automorphism $W$ of $\R^{2n}$, the following facts are equivalent:
\begin{enumerate}[(i)]
\item $W$ belongs to $\mathrm{Sp}^+_{\mathrm{ell}}(2n)$;
\item there exists a $\kappa$-unitary real basis $w_1,\dots,w_n,\overline{w}_1,\dots,\overline{w}_n$ of $\C^{2n}$ such that
\[
W w_j= e^{i \theta_j} w_j, \quad  W \overline{w}_j =e^{-i \theta_j} \overline{w}_j\qquad \forall j\in \{1,\dots,n\}.
\]
for some numbers $\theta_j\in (0,\pi)$;
\item there is a $W$-invariant symplectic splitting of $\R^{2n}$ into $n$ symplectic planes, i.e.,
\[
\R^{2n} = \bigoplus_{j=1}^n V_j, \qquad \dim V_j=2, \quad \omega_0(u,v) = 0\quad  \forall u\in V_j, \; v\in V_h \mbox{ with } j\neq h,
\]
with respect to which $W$ has the form
\[
W = \bigoplus_{j=1}^n e^{\theta_j J_j},
\]
where each $\theta_j$ is in the interval $(0,\pi)$ and each $J_j:V_j\to V_j$ is an $\omega_0$-compatible complex structure on the symplectic plane $V_j$. 
\end{enumerate}
\end{prop}

\begin{proof}
By Proposition \ref{propA1}, the space $\mathrm{sp}^+_{\mathrm{ell}}(2n)$ is precisely the set of all endomorphisms $X: \R^{2n} \rightarrow \R^{2n}$ of the form
\[
X = \bigoplus_{j=1}^n \theta_j J_j,
\]
where the direct sum refers to a symplectic splitting of $\R^{2n}$ into planes $V_1,\dots,V_n$, each $\theta_j$ belongs to the interval $(0,\pi)$, and each $J_j:V_j\to V_j$  is an $\omega_0$-compatible complex structure on the symplectic plane $V_j$. By moving the numbers $\theta_j$ so that they all become equal to $\pi/2$, this space is readily seen to be homotopically equivalent to the space of $\omega_0$-compatible complex structures on $\R^{2n}$, which is well-known to be contractible, see e.g., \cite[Lemma 2.5.5]{ms95}. This shows that the space $\mathrm{sp}^+_{\mathrm{ell}}(2n)$ is contractible. 

The exponential map is a local diffeomorphism on it thanks to Theorem \ref{Morsegen}. Let $X$ and $X'$ be elements of $\mathrm{sp}^+_{\mathrm{ell}}(2n)$ such that $e^X=e^{X'}$. We wish to show that $X=X'$. Thanks to Proposition \ref{propA1}, by considering the spectral decomposition of $X$ and $X'$ we are reduced to the case in which $X=\theta J$ and $X=\theta J'$, where $\theta\in (0,\pi)$ and $J, J'$ are two $\omega_0$-compatible complex structures on a symplectic vector subspace $V$ of $(\R^{2n},\omega_0)$. The identities
\[
e^{\theta J} = (\cos \theta) I + (\sin \theta)J, \qquad e^{\theta J'} = (\cos \theta) I + (\sin \theta)J',
\]
and the fact that $\sin \theta\neq 0$ imply that $J=J'$ and hence $X=X'$. This proves that the exponential map restricts to a diffeomorphism on the contractible open set $\mathrm{sp}^+_{\mathrm{ell}}(2n)$.

If $X$ is in $\mathrm{sp}^+_{\mathrm{ell}}(2n)$, then $e^X$ belongs to the set $\mathrm{Sp}^+_{\mathrm{ell}}(2n)$ as defined in  Section \ref{secJ} of the Introduction. Conversely, the fact that Krein-definite eigenvalues are semisimple (see \cite[Chapter I, Proposition 7]{eke90}), implies that any $W$ in $\mathrm{Sp}^+_{\mathrm{ell}}(2n)$ is semisimple. The equivalence of (i), (ii) and (iii) can now be deduced from the normal form of semisimple symplectic matrices, see e.g.\ \cite[Section 1.3.2]{abb01}. In particular, (iii) shows that any $W$ in $\mathrm{Sp}^+_{\mathrm{ell}}(2n)$ is the exponential of some $X$ in $\mathrm{sp}^+_{\mathrm{ell}}(2n)$. This concludes the proof.
\end{proof}

\section{Contact Hamiltonians}
\label{contham}

In this appendix, we collect for the reader's convenience some basic facts about the identification between contact vector fields and Hamiltonian functions. 

Let $\xi$ be a co-oriented contact structure on the closed manifold $M$ and let $\alpha$ be a contact form on $M$ defining $\xi$. The map
\begin{equation}
\label{the-ide-app}
\mathrm{cont}(M,\xi) \rightarrow C^{\infty}(M), \qquad X \mapsto \imath_X \alpha,
\end{equation}
is invertible. Indeed, its inverse is the map
\[
C^{\infty}(M) \rightarrow \mathrm{cont}(M,\xi)  , \qquad H \mapsto X_H,
\]
where $X_H$ is the unique vector field satisfying the identities 
\[
\imath_{X_H} \alpha = H, \qquad 
\imath_{X_H} \mathrm{d} \alpha = - \mathrm{d} H + \bigl(\imath_{R_{\alpha}} \mathrm{d} H\bigr) \alpha.
\]
Here, $R_{\alpha}$ denotes the Reeb vector field of $\alpha$. See \cite[Theorem 2.3.1]{gei08}. The function $H$ is the \textit{contact Hamiltonian} defining the contact vector field $X_H$.

\begin{rem}
\label{no-tangency}
{\rm The injectivity of the map (\ref{the-ide-app}) implies that the only contact vector field which is a section of the contact structure $\xi$ is the zero vector field.}
\end{rem}

If $H\in C^{\infty}(M)$ is positive, then $H^{-1} \alpha$ is a contact form defining $\xi$. Therefore, its Reeb vector field is a contact vector field and from the identity
\[
\imath_{R_{H^{-1}\alpha}} \alpha = H  \imath_{R_{H^{-1}\alpha}} (H^{-1} \alpha) = H 
\]
and the injectivity of (\ref{the-ide-app}) we deduce the identity
\begin{equation}
\label{Reeb-contact}
X_H = R_{H^{-1} \alpha}.
\end{equation}
This identity implies that the elements of $\mathrm{cont}^+(M,\xi)$ are precisely the Reeb vector fields of contact forms defining $\xi$. 

The adjoint action of $\mathrm{Cont}(M,\xi)$ on $\mathrm{cont}(M,\xi)$ is given by the push-forward
\[
\mathrm{Ad}_{\phi} X = \phi_* X.
\]
In terms of the contact Hamiltonians, this action reads
\begin{equation}
\label{adjoint-action}
\phi_* X_H= X_K \qquad \mbox{where}\quad K:= f^{-1} \phi_* H,
\end{equation}
and the function $f\in C^{\infty}(M)$ is defined by
\[
\phi_* \alpha = f \alpha.
\]
Indeed, this follows from the chain of identities
\[
f K = f \imath_{X_K} \alpha = \imath_{X_K} (f\alpha) = \imath_{\phi_* X_H} ( \phi_* \alpha) = \phi_* \bigl( \imath_{X_H} \alpha \bigr) = \phi_* H.
\]
In this monograph, the Lie bracket of two vector fields is defined by the non-standard sign convention
\[
[X,Y] = - \mathcal{L}_X Y,
\]
see the discussion in Section \ref{secG} from the Introduction. The Lie bracket of two contact vector fields is a contact vector field, and the \textit{contact Poisson bracket} $\{H,K\}\in C^\infty(M)$ of two functions $H,K\in C^{\infty}(M)$ is defined by the identity
\[
X_{\{H,K\}} = [X_H,X_K].
\]
One can show that
\begin{equation}\label{e:bracketHK}
\{H,K\} =\di H(X_K)-\di K(R_\alpha)H\, ,
\end{equation}
see \cite[Remark 3.5.18]{ms95}.

 \addcontentsline{toc}{section}{References}

 
 \providecommand{\bysame}{\leavevmode\hbox to3em{\hrulefill}\thinspace}
\providecommand{\MR}{\relax\ifhmode\unskip\space\fi MR }
\providecommand{\MRhref}[2]{%
  \href{http://www.ams.org/mathscinet-getitem?mr=#1}{#2}
}
\providecommand{\href}[2]{#2}

\end{document}